\documentclass[11pt]{article}

\usepackage{geometry}       
\usepackage{authblk}  

\usepackage{rotating}

\usepackage{amssymb}
\usepackage{amsmath}

\usepackage{multirow}
\usepackage{booktabs}

\usepackage[hidelinks]{hyperref} 
\usepackage{mathrsfs,xspace,xcolor,graphicx,subcaption,mathtools,empheq,
algorithm,algorithmic}
\usepackage{dirtytalk}
\newcommand{\lifex}{\texttt{life}$^{\color{red}\texttt{x}}$\xspace}

\DeclareMathOperator*{\argmin}{arg\,min}
\graphicspath{{./Images/}}

\title{Cardiocirculatory Computational Models for the Study of Hypertension}

\author[1]{Simone Celora}
\author[2]{Andrea Tonini}
\author[2]{Francesco Regazzoni}
\author[2]{Luca Dede'}
\author[3,4]{Gianfranco~Parati}
\author[2,5]{Alfio Quarteroni}

\affil[1]{Politecnico di Milano, Piazza Leonardo da Vinci~32, 20133 Milano, Italy}
\affil[2]{MOX, Department of Mathematics, Politecnico di Milano,\\ Piazza Leonardo da Vinci~32, 20133 Milano, Italy}
\affil[3]{Department of Cardiology, S.Luca Hospital, IRCCS, Istituto Auxologico Italiano,\\ Piazzale Brescia 20, 20149 Milano, Italy}
\affil[4]{Department of Medicine and Surgery, University of Milano-Bicocca,\\ via Vizzola 5, 20126 Milano, Italy}
\affil[5]{Institute of Mathematics, École Polytechnique Fédérale de Lausanne,\\ Station 8, Av.~Piccard, CH-1015 Lausanne, Switzerland}

\date{\today}

\begin{document}
\maketitle

\begin{abstract}
In this work, we develop patient-specific cardiocirculatory models with the aim of 
building Digital Twins for hypertension. In particular, in our pathophysiology-based 
framework, we consider both 0D cardiocirculatory models and a 3D-0D electromechanical 
model.
  The 0D model, which consists of an RLC circuit, is studied in two 
  variants, with and without capillaries. The 3D--0D model consists of a 
  three-dimensional electromechanical model of the left ventricle, coupled with a 0D 
  model for the external
  blood circulation: this representation enables the assessment of additional 
  quantities related to ventricular deformation and stress, 
  and offers 
  a more detailed representation compared to a fully 0D model.
  Sensitivity analysis is performed on the 0D model, with both a mono- and a 
  multi-parametric approach, in order to identify the parameters that most influence 
  the model outputs and guide the calibration process.
We studied three different scenarios, corresponding to systemic, pulmonary and 
renovascular hypertension, each in three nuances of severity.
  To maintain a fair comparison among the models, a parameter calibration 
  strategy is developed;
  the outputs of the 0D model with capillaries are utilized to 
  enhance the 3D--0D model. 
  The results 
  demonstrate that the 3D--0D model yields an accurate representation of 
  cardiocirculatory dynamics in the presence of hypertension; this model represents a 
  powerful step toward digital twins for real-time hypertension control, 
  providing refined and clinically meaningful insights beyond those achievable with 
  0D models alone.
\end{abstract}




\section{Introduction}
\label{intro}%
Hypertension, along with higher than normal or elevated blood 
pressure levels, is one of the world's leading risk factors for death and disability 
\cite{who2023, who2025}, since it may cause stroke, ischaemic heart disease,
altered vascular conditions and kidney 
disorders \cite{intro:1, intro:2}. 
High systolic blood 
pressure is responsible for more deaths than any other health risk. According to \cite{who2023}, it accounted for over 
10~million deaths annually and contributed to over half of all 
cardiocirculatory-related deaths worldwide in 2019 \cite{intro:5}; furthermore, in less than 30 years, the number of hypertensive individuals doubled: from 
650 million in 1990 to 1.3 billion in 2019 \cite{intro:3}. According to the most recent World Health Organization report~\cite{who2025} ``in 2024, an estimated 1.4 billion people aged 30--79 years were affected worldwide, yet fewer than 1 in 5 (320 million) had the condition adequately controlled".

Hypertension is often referred to as the ``silent killer'' because it is typically 
asymptomatic: unless 
blood pressure is measured, most individuals remain unaware 
of their condition until they experience a clinical complication, i.e. heart 
attack, stroke or kidney failure. Proper management of hypertension requires
accurate diagnosis, prompt initiation 
of treatment in severe cases and continuous monitoring 
of progress. Timelines are crucial: the longer a person suffers from 
undiagnosed or poorly managed hypertension, the 
worse the health consequences 
are likely to be \cite{who2023, who2025}.

Due to hypertension profound impact on global health, a comprehensive understanding of 
this condition is crucial.
Mathematical models are increasingly recognized as indispensable tools for 
studying complex physiological phenomena, offering insights that traditional approaches 
alone may not provide. 
Within this context, the concept of Digital Twins (DTs) is gaining ground in cardiovascular 
research \cite{dt2, dt3}. 
Unlike a simple patient-specific model, a digital twin establishes a continuous 
interaction between the physical and digital entities: patient data drive 
the calibration of the computational model, while the predictions obtained from the 
model can guide clinical choices, such as the selection of therapeutic strategies 
\cite{dt1}.
By integrating clinical and sensor data into mathematical models, DTs allow simulations 
tailored to the individual patient, enabling accurate reproduction of haemodynamics 
and prediction of disease progression.
Hypertension has already been investigated using machine 
learning models, as evidenced 
by works such as \cite{es:mach_learn1,es:mach_learn2}; nonetheless, 
there has also been an increasing use 
of differential equation (DE)-based models to simulate the cardiocirculatory 
system and study hypertension (e.g. \cite{model_hyper:es,es:hyp_study}). 
The use of DE-based models, combined with real 
patient data, enables the creation 
of personalized models tailored to individual patients: by incorporating 
patient-specific physiological information, these models are able to simulate 
cardiocirculatory haemodynamics with greater accuracy and reflect the unique characteristics 
of each patient's condition.
When combined with real patient data, such DE-based models naturally lend themselves 
to the construction of DTs, as they can reproduce patient-specific cardiocirculatory 
dynamics and reflect the unique characteristics of individual conditions.

Albeit advanced and complex models are increasingly 
being used as complementary tools to the study of pathological conditions, hypertension 
modeling has 
been little explored in literature. This paper is proposed as an innovative 
contribution in this direction, 
providing a framework that not only simulates different forms of hypertension but also 
provides a foundation for future developments toward cardiovascular DTs.
The adopted approach, which combines 
lumped-parameter models and detailed electromechanical representations, makes it 
possible to explore the dynamics of the cardiocirculatory 
system with a hitherto little studied level of detail: the current work
stands as a forerunner in the modeling of hypertension, helps to fill a gap in the 
existing literature and 
supports the advancement of patient-specific DTs for predictive and 
personalized analysis.

Several differential mathematical and numerical models are employed to 
analyze the dynamics of hypertension, 
simulate its effects under various conditions and compare the simulated outputs with 
clinically observed trends and reference values from the literature:
in particular, a lumped-parameter cardiocirculatory model 
(from now on referred to as $(\mathscr{C})$ or 0D model), based on a 0D 
representation by an RLC circuit, 
and an electromechanical model (\cite{FEDELE2023115983}) coupled with a 
lumped-parameter model (from now on referred to as 3D--0D model), 
where the left ventricle is represented in 3D through 
sophisticated electromechanical modeling, coupled with a 0D model of the blood 
circulation.
The integrated approach not only 
enhances the understanding of 
hypertension, but also contributes to the development of improved diagnostic, 
therapeutic and preventive strategies,
within a DT framework.
Through the proposed modeling 
framework, the simulations reproduce the 
haemodynamic alterations associated with various forms and grades of hypertension; 
the results show the effect of increased vascular resistance and reduced compliance 
on both systemic and pulmonary circulation, demonstrating how these changes deeply affect 
pressure dynamics and cardiac workload. Furthermore, the comparison between 0D and 3D--0D 
models yields valuable insights into advantages and limitations of each approach, 
emphasizing the importance of patient-specific modeling for an accurate 
representation of hypertensive conditions.

The current paper is organized as follows: Section~\ref{hyper} provides a clinical overview of 
hypertension, with a brief description of its causes and effects;
in Section~\ref{models}, the employed mathematical models are presented, 
along with an analysis of the relevant parameters involved;
Section~\ref{hyper-model} illustrates the implementation of 
different kinds of hypertension within 
the mathematical framework, with a description of the necessary modifications and 
assumptions for correct modeling;
Section~\ref{results} is devoted to numerical results obtained 
from the simulations, comparison of the model outputs with real-world clinical data 
and evaluation of their consistency; finally, conclusions 
follow in Section~\ref{conclusions}.

\section{Pathophysiology of Hypertension}
\label{hyper}
In accordance with the 2023 European and current international guidelines 
\cite{guidelines:eu}, hypertension is diagnosed when office
blood pressure 
measurements (i.e. taken in a 
clinical or healthcare setting) is equal or higher than 140~mmHg for systolic blood pressure or 90~mmHg 
for diastolic blood pressure, although recent USA guidelines have adopted a lower threshold, i.e. equal or higher than 130/80~mmHg.
Alternative definitions are considered to identify related conditions 
such as pulmonary hypertension, as described below.
Hypertension can be classified in several ways
considering aspects such as the underlying causes, the severity of blood 
pressure elevation, the timing and consistency of measurements,  
the presence of associated complications, or other risk factors. 
An important distinction is made between primary and secondary hypertension: while primary 
hypertension has no clearly identifiable cause and accounts for the vast majority of 
cases, secondary hypertension arises from specific, diagnosable medical conditions (e.g. 
renal, endocrine or vascular disorders) \cite{percent_esshyper}.
Hypertension can be further classified based on its origin and the vascular district 
involved. Here, the forms that are more related to this work are briefly described:
\begin{enumerate}
    \item systemic hypertension refers to elevated pressure in the systemic circulation 
    and is the most common form \cite{percent_esshyper}. It is typically associated with increased peripheral 
    resistance and large artery stiffening, affecting left ventricular function and 
    contributing to long-term cardiovascular damage;
    \item pulmonary hypertension by contrast, is defined as a mean pulmonary artery 
    pressure above 20~mmHg at rest \cite{pulhyper:guidelines}. It affects the pulmonary circulation and can be 
    caused by a variety of mechanisms, including left heart dysfunction, lung diseases, hypoxia exposure,  
    or vascular abnormalities. It leads to increased right ventricular workload and is 
    associated with a high morbidity and mortality risk;
    \item renovascular hypertension is a form of secondary hypertension caused by 
    reduced renal perfusion, most commonly due to atherosclerotic renal artery stenosis 
    or fibromuscular dysplasia \cite{reno-def}. The decrease in renal blood flow promotes renin secretion,
    vasoconstriction, sympathetic 
    stimulation and fluid retention, 
    all of which 
    contribute to a sustained elevation 
    in systemic blood pressure.
\end{enumerate}

\subsection{Etiology and Clinical Manifestations}
\label{sec:cause_sympt}%
Hypertension is a multifactorial condition with both genetic and environmental influences. 
It may be asymptomatic in its early stages but can lead to progressive damage to key 
organs such as the heart, brain, vasculature and kidneys as blood pressure rises 
\cite{guidelines:2024esc}. 
Hypertensive individuals often exhibit enhanced vasoconstriction in both large 
and small arteries. This phenomenon, alongside the increase in blood pressure, appears 
to be favored by systemic autoregulation, 
triggered by the preliminary increase of blood 
volume and cardiac output; this mechanism refers to the ability of blood vessels to maintain a 
constant blood flow, despite changes in blood pressure \cite{autoregulation}. 
Hypertension is very commonly associated with impaired vasodilation, due mostly to 
endothelial dysfunction and structural changes in blood vessels \cite{vasc5}. 
In small vessels, remodeling of arterioles results in smooth muscle 
layer thickening and lumen narrowing, 
increasing overall vascular resistance \cite{vasc11}.
Other studies documented a two-way relation between aortic stiffness and 
hypertension, as one condition increases due to the other one
in a 
self-promoting circle. A known impact of stiffness is the transmission of pressure into 
the smaller vessels, which may cause harm to the organs \cite{vasc20}.

Chronic high blood pressure puts stress on the left ventricular wall by 
increasing the afterload, promoting left ventricular remodeling and hypertrophy; 
the latter is a well-established risk factor for adverse cardiocirculatory outcomes, 
such as myocardial infarction, heart failure, atrial fibrillation, stroke and 
cardiocirculatory mortality \cite{guidelines:2024esc}. 
Long-term elevated pressure in the left ventricle may lead to the 
dysfunction of the left atrium through dilation.
There is ample evidence linking hypertension with the etiology of valvular heart 
disease \cite{symptoms:heart14}: an increased systolic blood pressure has been 
associated with an elevated risk of conditions such as aortic stenosis and aortic 
regurgitation. 
The effects of hypertension involve endothelial dysfunction, vascular remodeling and 
enhanced vascular stiffness, which not only contribute to the development of the 
disease, but are also results of it in turn \cite{symptoms:heart17}. 
The kidney is also a key actor in the development and sustainment of hypertension: 
regulating systemic vasoconstriction and blood pressure through 
both direct hormonal actions and feedback mechanisms,
including fluid balance, sodium handling and sympathetic activity 
\cite{kidney6,kidney8}.
Hypertension and chronic kidney disease are closely linked in a complex relationship: 
the former can be considered both a cause and a consequence of the latter, 
promoting its progression \cite{symptoms:kidney2}.

\section{Differential Mathematical Models for the Cardiocirculatory System}
\label{models}
This section provides an overview of the models adopted for this paper, describing 
their mathematical formulation, the key assumptions, the functionality and the 
involved parameters \cite{3d-0d_ventr,model-whole-heart1,lumped:capillary}.
In Section~\ref{subsec:lumped}, the 0D models are illustrated, 
while the 3D--0D coupled model is described and 
briefly analyzed in Section~\ref{subsec:3d-0d-model}.

\subsection{A 0D Model for Closed-Loop Blood Circulation}
\label{subsec:lumped}%
\begin{figure}[t!]
    \centering
    \subfloat[0D model.\label{fig:lumped_nocap}]{
        \includegraphics[height=0.45\linewidth]{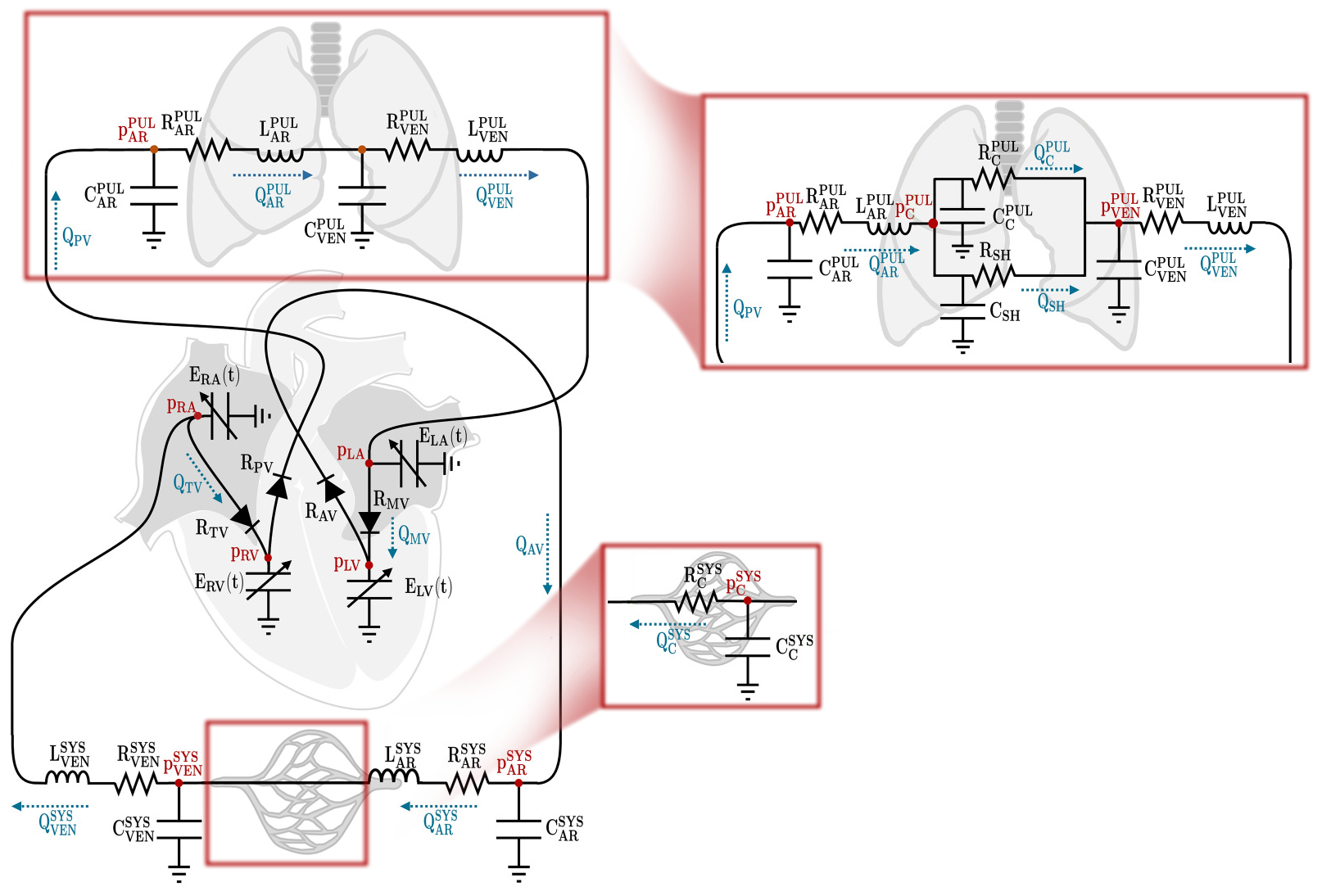}
    }
    \subfloat[3D--0D model.\label{fig:3dmodel}]{
        \includegraphics[height=0.45\linewidth]{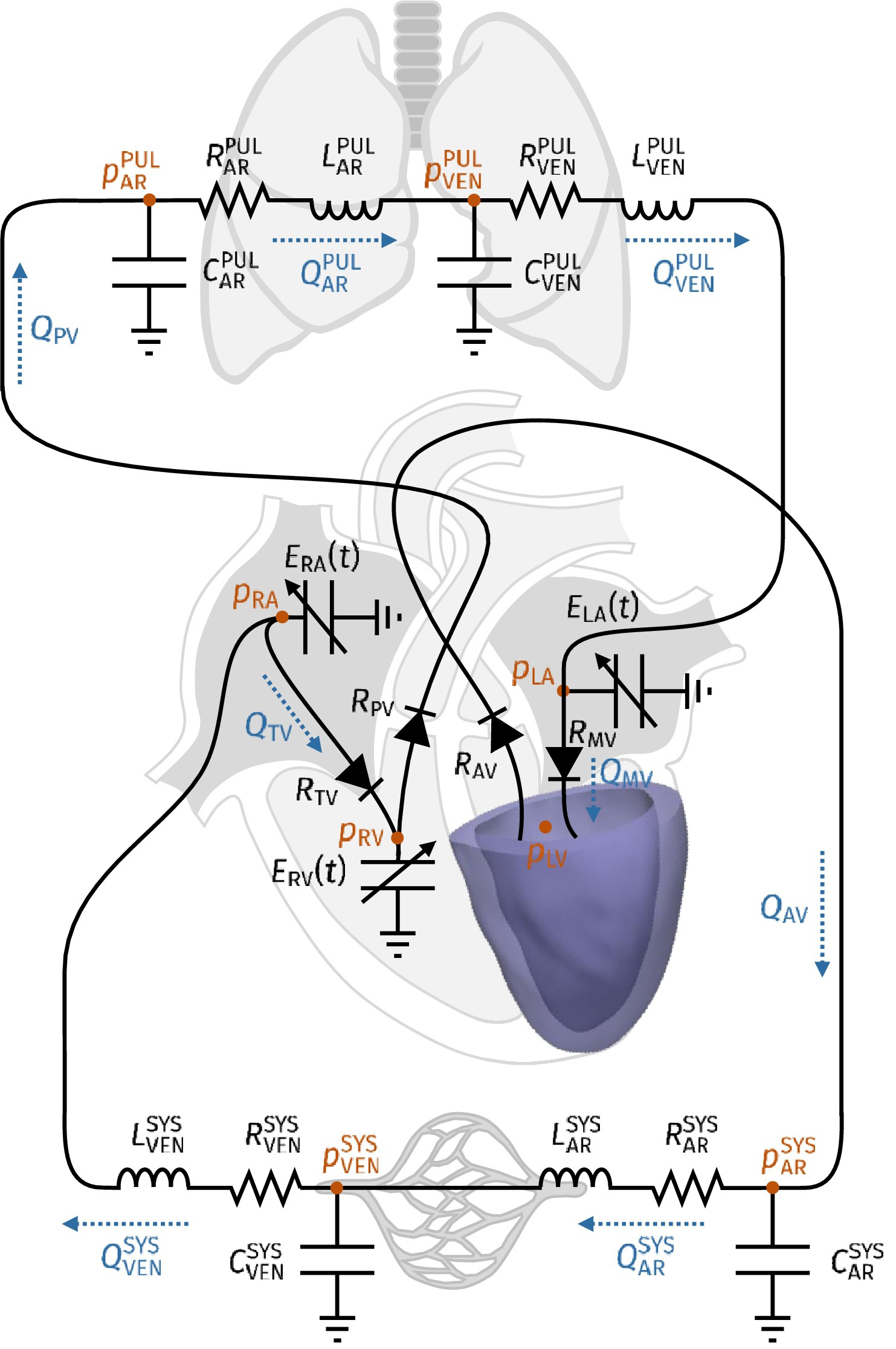}
    }
    \caption{
        0D cardiocirculatory models 
    with (center) and without (left) capillaries. The 3D--0D coupling between the 
    left ventricle 3D electromechanical 
    model and the 0D circulation model without capillaries is shown on the right.
    Source: \cite{3d-0d_ventr,lumped:capillary} and graphical modifications.}
    \label{fig:lumped_model}
\end{figure}

The 0D model $(\mathscr{C})$ \cite{3d-0d_ventr,lumped:capillary}
describes the haemodynamics of the cardiocirculatory system by means of a 
lumped-parameter approach: it represents the entire 
circulatory network as a set of 0D interconnected components, and allows the system to be 
modeled using ordinary differential equations (ODEs). $(\mathscr{C})$ is based on the 
assumption that the cardiocirculatory system can be simplified into discrete elements that 
capture the mean features of blood flow and pressure dynamics. 
This type of models have been implemented in two variants 
(Figure~\ref{fig:lumped_nocap}): $(\mathscr{C}_{\text{C}})$ which includes capillaries
and $(\mathscr{C}_{\text{NC}})$ which does not.
The circulation is 
simulated over a time period T, long enough to reach a periodic limit cycle in the output 
variables, with period that of the cardiac cycle; this periodicity is a property of the 
simulation rather than the model itself. Only results belonging to the last 
heartbeat, lasting $\text{T}_{\text{HB}}$, are then considered.

In $(\mathscr{C})$, both the systemic and pulmonary circulations are represented by a 
resistance-inductance-capacitance (RLC) circuit: the resistance corresponds to the 
opposition to blood flow, the inductance accounts for the inertial effects of the blood, 
and the capacitance represents the elasticity of the vessel walls. 
The heart is modeled as a series of time-varying elastance elements, describing 
dynamic variations of the force developed by the cardiac muscle during the cardiac 
cycle. Such elements characterize each of the four cardiac chambers,
allowing for the simulation of 
contraction and relaxation phases of the heart. 
The four heart valves 
are represented as 
non-ideal diodes:
they capture the unidirectional blood flow and the 
resistance to backflow when the valves are closed. The heart elements are integrated 
into a closed-loop system, where the interactions between the heart and the blood 
vessels are modeled to replicate the complete haemodynamic behaviour of the 
cardiocirculatory system.

The two variations of the model are based on distinct systems of ODEs, 
reflecting the inclusion or exclusion of capillaries; below, the system of ODEs 
for $(\mathscr{C}_{\text{NC}})$ is presented: 
\begin{center}
\rotatebox{90}{
\begin{minipage}{0.8\textheight}
\begin{subequations}\label{eq:system}
    \begin{empheq}[left=\empheqlbrace]{align}
    \dfrac{dV_{\text{RA}}(t)}{dt}&=Q_\text{VEN}^\text{SYS}(t) - Q_\text{TV}(t),  &\dfrac{dV_{\text{LA}}(t)}{dt}&=Q_\text{VEN}^\text{PUL}(t) - Q_\text{MV}(t), \label{eq:system1}\\
    \dfrac{dV_{\text{RV}}(t)}{dt}&=Q_\text{TV}(t) - Q_\text{PV}(t),  & \dfrac{dV_{\text{LV}}(t)}{dt}&= Q_\text{MV}(t) - Q_\text{AV}(t), \label{eq:system2}\\ 
    \dfrac{dp_{\text{VEN}}^{\text{SYS}}(t)}{dt}&=\dfrac{{Q_\text{AR}^\text{SYS}(t) - Q_\text{VEN}^\text{SYS}(t)}}{{\text{C}_\text{VEN}^\text{SYS}}},  & \dfrac{dp_{\text{VEN}}^{\text{PUL}}(t)}{dt}&=\dfrac{{Q_\text{AR}^\text{PUL}(t) - Q_\text{VEN}^\text{PUL}(t)}}{{\text{C}_\text{VEN}^\text{PUL}}},\label{eq:system3}\\ 
    \dfrac{dp_{\text{AR}}^{\text{PUL}}(t)}{dt}&=\dfrac{{Q_\text{PV}(t) - Q_\text{AR}^\text{PUL}(t)}}{{\text{C}_\text{AR}^\text{PUL}}}, & \dfrac{dp_{\text{AR}}^{\text{SYS}}(t)}{dt}&=\dfrac{{Q_\text{AV}(t) - Q_\text{AR}^\text{SYS}(t)}}{{\text{C}_\text{AR}^\text{SYS}}},\label{eq:system4}\\ 
    \dfrac{dQ_{\text{VEN}}^{\text{SYS}}(t)}{dt}&=-\left(Q_\text{VEN}^\text{SYS}(t) + \dfrac{{p_\text{RA}(t) - p_\text{VEN}^\text{SYS}(t)}}{{\text{R}_\text{VEN}^\text{SYS}}}\right) \dfrac{{\text{R}_\text{VEN}^\text{SYS}}}{{\text{L}_\text{VEN}^\text{SYS}}},  &\dfrac{dQ_{\text{VEN}}^{\text{PUL}}(t)}{dt}&=-\left(Q_\text{VEN}^\text{PUL}(t) + \dfrac{{p_\text{LA}(t) - p_\text{VEN}^\text{PUL}(t)}}{{\text{R}_\text{VEN}^\text{PUL}}}\right) \dfrac{{\text{R}_\text{VEN}^\text{PUL}}}{{\text{L}_\text{VEN}^\text{PUL}}},\label{eq:system5}\\
    \dfrac{dQ_{\text{AR}}^{\text{PUL}}(t)}{dt}&=-\left(Q_\text{AR}^\text{PUL}(t) + \dfrac{{p_\text{VEN}^\text{PUL}(t) - p_\text{AR}^\text{PUL}(t)}}{{\text{R}_\text{AR}^\text{PUL}}} \right) \dfrac{{\text{R}_\text{AR}^\text{PUL}}}{{\text{L}_\text{AR}^\text{PUL}}}, & \dfrac{dQ_{\text{AR}}^{\text{SYS}}(t)}{dt}&=-\left(Q_\text{AR}^\text{SYS}(t) + \dfrac{{p_\text{VEN}^\text{SYS}(t) - p_\text{AR}^\text{SYS}(t)}}{{\text{R}_\text{AR}^\text{SYS}}} \right)  \dfrac{{\text{R}_\text{AR}^\text{SYS}}}{{\text{L}_\text{AR}^\text{SYS}}}, \label{eq:system6}
\end{empheq}
\end{subequations}
\end{minipage}
}
\end{center}
for $t\in(0,\text{T}]$, and coupled with suitable initial conditions.

$(\mathscr{C}_{\text{C}})$ introduces modifications to 
Equations~\eqref{eq:system3} and \eqref{eq:system6}, 
respectively:
\begin{subequations}\label{eq:diff}
    \begin{align}
        \dfrac{dp_{\text{VEN}}^{\text{SYS}}(t)}{dt}&=\dfrac{{Q_\text{C}^\text{SYS}(t) - Q_\text{VEN}^\text{SYS}(t)}}{{\text{C}_\text{VEN}^\text{SYS}}}, &\dfrac{dp_{\text{VEN}}^{\text{PUL}}(t)}{dt}&=\dfrac{{Q_\text{SH}(t) + Q_\text{C}^\text{PUL}(t) - Q_\text{VEN}^\text{PUL}(t)}}{{\text{C}_\text{VEN}^\text{PUL}}}, \tag{2c}\label{eq:diff1}\\ 
        \dfrac{dQ_{\text{AR}}^{\text{PUL}}(t)}{dt}&=-\left(Q_\text{AR}^\text{PUL}(t) + \dfrac{{p_\text{C}^\text{PUL}(t) - p_\text{AR}^\text{PUL}(t)}}{{\text{R}_\text{AR}^\text{PUL}}} \right)  \dfrac{{\text{R}_\text{AR}^\text{PUL}}}{{\text{L}_\text{AR}^\text{PUL}}}, & \dfrac{dQ_{\text{AR}}^{\text{SYS}}(t)}{dt}&=-\left(Q_\text{AR}^\text{SYS}(t) + \dfrac{{p_\text{C}^\text{SYS}(t) - p_\text{AR}^\text{SYS}(t)}}{{\text{R}_\text{AR}^\text{SYS}}} \right)  \dfrac{{\text{R}_\text{AR}^\text{SYS}}}{{\text{L}_\text{AR}^\text{SYS}}}, \tag{2f}\label{eq:diff2}
    \end{align}
\end{subequations}\addtocounter{equation}{-1}
for $t\in(0,\text{T}]$, and adds the additional ones, here presented:
\begin{subequations}
\begin{align}
        \dfrac{dp_{\text{C}}^{\text{SYS}}(t)}{dt}&=\dfrac{{Q_\text{AR}^\text{SYS}(t) - Q_\text{C}^\text{SYS}(t)}}{{\text{C}_\text{C}^\text{SYS}}}, & \dfrac{dp_{\text{C}}^{\text{PUL}}(t)}{dt}&=\dfrac{{ Q_\text{AR}^\text{PUL}(t) - Q_\text{SH}(t) - Q_\text{C}^\text{PUL}(t)}}{{\text{C}_\text{C}^\text{PUL}+ \text{C}_\text{SH}}}, \tag{2g}\label{eq:add_cap}
    \end{align}
\end{subequations}
for $t\in(0,\text{T}]$, 
to account for the haemodynamics of the capillary network, 
while the others remain unchanged between the two models.

In both the models, the pressures in the four heart chambers are computed as:
\begin{equation}
    p_i(t)=p_{\text{EX}}(t)+E_i(t)\,(V_i(t)-\text{V}_{0,i}),\quad i\in\{\text{LA},\text{LV},\text{RA},\text{RV}\}, \label{chamber:pressure}
\end{equation}
while the fluxes through the four heart valves are defined as:
\begin{subequations}\label{eq:fluxes}
\begin{align}
    Q_{\text{TV}}(t)&=\dfrac{p_{\text{RA}}(t)-p_{\text{RV}}(t)}{R_{\text{TV}}\bigl(p_{\text{RA}}(t),p_{\text{RV}}(t)\bigr)}, & Q_{\text{PV}}(t)&=\dfrac{p_{\text{RV}}(t)-p_{\text{AR}}^{\text{PUL}}(t)}{R_{\text{MV}}\bigl(p_{\text{RV}}(t),p_{\text{AR}}^{\text{PUL}}(t)\bigr)},\label{QPV}\\
    Q_{\text{MV}}(t)&=\dfrac{p_{\text{LA}}(t)-p_{\text{LV}}(t)}{R_{\text{MV}}\bigl(p_{\text{LA}}(t),p_{\text{LV}}(t)\bigr)}, & Q_{\text{AV}}(t)&=\dfrac{p_{\text{LV}}(t)-p_{\text{AR}}^{\text{SYS}}(t)}{R_{\text{AV}}\bigl(p_{\text{LV}}(t),p_{\text{AR}}^{\text{SYS}}(t)\bigr)},\label{QAV}
    \end{align}
\end{subequations}
where $\text{V}_{0,i}$ is the resting volume of heart chamber $i$; $p_{\text{EX}}(t)$ 
is the pressure exerted outside the heart by the surrounding organs and respiration, and 
in this work will always be set to 0.

Equations~\eqref{eq:diff1}, \eqref{eq:diff2} and \eqref{eq:add_cap}
introduce the following new flux variables:
\begin{align}
    Q_\text{C}^\text{SYS}&=\dfrac{p_\text{C}^\text{SYS}(t)-p_\text{VEN}^\text{SYS}(t)}{\text{R}_\text{C}^\text{SYS}}, & Q_\text{C}^\text{PUL}&=\dfrac{p_\text{C}^\text{PUL}(t)-p_\text{VEN}^\text{PUL}(t)}{\text{R}_\text{C}^\text{PUL}}, & Q_\text{SH}&=\dfrac{p_\text{C}^\text{PUL}(t)-p_\text{VEN}^\text{PUL}(t)}{\text{R}_\text{SH}}.\label{eq:fluxes-cap}
\end{align}
A compact representation suitable for both the 0D models is provided here:
\begin{align}\label{0d-eq}
    (\mathscr{C})\,\,\,
    \begin{cases}
        \dfrac{d\boldsymbol{c}_1(t)}{d t}=\mathbf{D}\bigl(t,\boldsymbol{c}_1(t),\boldsymbol{c}_2(t)\bigr) &\text{for }t\in(0,\text{T}],\\
        \boldsymbol{c}_2(t)=\mathbf{W}\bigl(t,\boldsymbol{c}_1(t)\bigr) &\text{for }t\in[0,\text{T}],\\
        \boldsymbol{c}_1(0)=\boldsymbol{c}_{1,0}.
    \end{cases}
\end{align}
The same notation of Equation~\eqref{0d-eq} can be used to represent both 
$(\mathscr{C}_\text{C})$ and $(\mathscr{C}_\text{NC})$.

\subsection{The 3D Electromechanical Model coupled with the Lumped-Parameter Model}
\label{subsec:3d-0d-model}%
The 3D--0D model employed in this paper is the one 
presented in \cite{3d-0d_ventr}, which couples a detailed 3D electromechanical model of the left 
ventricle with a 0D model of the systemic and pulmonary circulations. 
The 3D left ventricle is represented by the computational domain 
$\Omega_0 \subset \mathbb{R}^3$, with its boundary $\partial\Omega_0$ being split 
into three distinct sections:
the endocardium ($\Gamma^\text{endo}_0$), 
the epicardium ($\Gamma^\text{epi}_0$)
and the ventricular base ($\Gamma^\text{base}_0$); the latter represents the artificial 
boundary where the left ventricle geometry is cut.

The multiphysics and multiscale model of cardiac electromechanics 
comprises five separate core models:
\begin{enumerate}
    \item transmission of cardiac electrical potential $(\mathscr{E})$;
    \item ion dynamics $(\mathscr{I})$;
    \item contraction of cardiac muscle cells $(\mathscr{A})$;
    \item mechanical behaviour of tissue $(\mathscr{M})$;
    \item blood circulation $(\mathscr{C})$.
\end{enumerate}
These core models capture the various physical processes involved in cardiac function, 
at different spatial and temporal dimensions. 
Furthermore, a volume conservation constraint $(\mathscr{V})$ ensures the proper 
interaction between $(\mathscr{M})$ and $(\mathscr{C})$.
In Figure~\ref{fig:3dmodel}, the complete 3D--0D model is presented.
Such a model includes multiple unknowns, depending on both time and space, that 
represent the variables related to each core model and describe the interactions 
between them; $u$ denotes the transmembrane potential, 
$\boldsymbol{w}$ and $\boldsymbol{z}$ represent the ionic variables, 
$\boldsymbol{s}$ indicates the state variables of the force generation model, 
$\boldsymbol{d}$ refers to the mechanical displacement of the tissue, 
$\boldsymbol{c}_1$ is the state vector of the circulation model, 
and $p_{\text{LV}}$ represents the pressure of left ventricle.

The 3D--0D complete model is defined by sets of partial differential equations 
(PDEs); the detailed formulation is reported in \ref{appendix:3d-0d_eq}.
In $(\mathscr{E})$, the dynamics of the transmembrane potential $u$ are captured by a diffusion-reaction 
PDE, which 
governs the electrical behaviour of cardiac muscle cells 
\cite{3d-0d:model-e:behaviour1,3d-0d:model-e:behaviour2,3d-0d:model-e:behaviour3}.
For the purpose of 
electromechanic coupling, a proper modeling of the ionic fluxes across the 
cell membrane is mandatory \cite{3d-0d:model-e:book}; the latter involves the 
inclusion of two variables: the recovery variables $\boldsymbol{w}$, 
which characterize the proportion of open ionic channels, and the 
concentration variables $\boldsymbol{z}$, describing the concentrations of critical 
ionic species, 
such as intracellular calcium ions $\left[\text{Ca}^{2+}\right]_i$, which play a pivotal 
role in triggering mechanical contraction.

In the electrophysiological framework defined by $(\mathscr{E})-(\mathscr{I})$, 
the transmembrane potential $u$ propagates through gap junctions; this phenomenon is 
governed by the diffusion term 
$\nabla\cdot\left(J\mathbf{F}^{-1}\mathbf{D}_\text{M}\mathbf{F}^{-T}\nabla u\right)$, 
which incorporates the influence of tissue deformation and 
introduces a feedback mechanism, linking electrical behaviour to mechanical stretch.
Here, $\mathbf{F}=\mathbf{I}+\nabla \boldsymbol{d}$ is the deformation gradient tensor, while $\mathbf{D}_\text{M}$ 
represents the diffusion tensor, which captures the anisotropic characteristics of cardiac tissue.
The applied current, $\mathcal{I}_{\text{app}}(t)$, simulates the role of the Purkinje 
network: these fibers form the 
outermost component of the cardiac conduction system and their primary function is to 
ensure the swift and synchronized activation of the ventricular myocardium \cite{purkinje}.
Additionally, the ionic current 
$\mathcal{I}_{\text{ion}}(u,\boldsymbol{w},\boldsymbol{z})$ encapsulates multiscale 
interactions at both cellular and tissue levels. Electrically isolated interfaces 
are modeled using homogeneous Neumann boundary conditions.

Cardiac contraction arises from 
interactions between actin and myosin proteins within sarcomeres, the fundamental 
contractile units of heart muscle \cite{3d-0d:model-a}. 
To simulate the mechanism in $(\mathscr{A})$, the RDQ20-MF 
model has been employed \cite{RDQ20} instead of RDQ18 in \cite{3d-0d_ventr}: 
it provides a complete biophysical model that characterizes the behaviour of sarcomeric 
proteins, along with the causes that support the responsiveness of
the heart towards variations in calcium levels.
RDQ20-MF is a system 
of differential equations, where the state variables are represented 
by the vector $\boldsymbol{s}$; it outputs the permissivity $P$, 
representing the fraction 
of contractile units in the force-generating state, from which the effective active 
tension $T_a$ is derived. 
Inputs to the model include calcium concentration $\left[\text{Ca}^{2+}\right]_i$ (properly 
described using a ionic model), sarcomere length (SL) and its time derivative from the mechanical model.

Tissue displacement $\boldsymbol{d}$ is described by the momentum conservation 
equation $(\mathscr{M})$. The Piola-Kirchhoff stress tensor 
$\mathbf{P}$ represents both the passive and active 
mechanical properties of the tissue. Under the assumption of hyperelasticity, the 
passive component of the stress tensor is derived from the strain energy density 
function $\mathcal{W}$, while the active component depends on tissue stretch along the 
fiber direction and on the active tension $T_a$.
To capture the anisotropic properties of cardiac muscle, the Guccione strain energy 
density function as $\mathcal{W}(\mathbf{F})$ \cite{guccione1,guccione2} is used, which 
depends on the Green-Lagrange strain energy tensor.
Additionally, to penalize large volume variations and enforce weak incompressibility, 
a term involving $J = \det \mathbf{F}$ is added to $\mathcal{W}(\mathbf{F})$ \cite{3d-0d:model-m}.

To model the interaction between the left ventricle and the pericardium 
\cite{3d-0d:model-m2}, a generalized 
Robin boundary condition is applied on the epicardium $\Gamma_0^\text{epi}$: 
this condition is designed to simulate the effect of 
the right ventricle, major veins and arteries, i.e. to restrict rigid rotations of 
the ventricle around the 
apico-basal axis (the direction from the base of the heart to its apex), while permitting torsion.
At the base $\Gamma^\text{base}_0$, an energy-consistent boundary condition is enforced: 
this ensures the correct stress distribution at the base boundary.
For the endocardium $\Gamma^\text{endo}_0$, the applied boundary condition 
accounts for the pressure $p_{\text{LV}}(t)$, which reflects the pressure exerted by the 
blood within the ventricle. 
Lastly, the mechanical model $(\mathscr{M})$ influences 
$(\mathscr{A})$ by determining the local sarcomere length (SL); the sarcomeres are 
aligned with the muscle fibers, so the local sarcomere length directly correlates with 
the tissue stretch along the fiber direction.

In the 3D--0D coupled framework, the 0D cardiocirculatory model used does not include 
the capillary network. The 0D model, as introduced in Section~\ref{subsec:lumped}, 
represents each cardiac chamber as a time-varying elastance element and yields a 
simplified 0D representation of the circulatory system. However, the coupled model 
integrates this 0D circulatory model with the 3D left ventricular model, which is 
governed by equations $(\mathscr{E})$, $(\mathscr{I})$, $(\mathscr{A})$ and 
$(\mathscr{M})$.
To achieve this integration, the elastance element of the left ventricle is removed from 
the 0D model and replaced by the 3D electromechanical model. Additionally, 
the 3D--0D coupled model must satisfy a 
volume consistency condition $\forall t \in (0, \text{T}]$, as described by 
Equation~\eqref{3d0d-eq:coupling}; 
$V^{\text{0D}}_{\text{LV}}(\boldsymbol{c}(t))$ 
and $V^{\text{3D}}_{\text{LV}}(\boldsymbol{d}(t))$ 
stand for the left ventricular volumes in the 0D and 3D models, respectively.
The introduction of this volume-consistency condition leads to an additional unknown 
variable, $p_{\text{LV}}$, which is no longer determined by 
Equation~\eqref{chamber:pressure}, but serves as a Lagrange multiplier enforcing 
the constraint $(\mathscr{V})$. Consequently, a ``reduced'' vector 
$\tilde{\boldsymbol{c}}_2$ is defined, where 
$\boldsymbol{c}_2^T=\left(p_\text{LV}, \tilde{\boldsymbol{c}}_2^T\right)$: this allows 
for the reformulation of Equations~\eqref{chamber:pressure} and \eqref{eq:fluxes} as:
\begin{equation*}
    \tilde{\boldsymbol{c}}_2(t)=\overset{\sim}{\mathbf{W}}\left(t,\boldsymbol{c}_1(t),p_\text{LV}(t)\right).
\end{equation*}
As a result, the reduced version of Equation~\eqref{0d-eq} can be expressed as 
Equation~\eqref{3d0d-eq:circ}, where the following expression for 
$\overset{\sim}{\mathbf{D}}$ is employed:
\begin{equation*}
    \overset{\sim}{\mathbf{D}}\left(t,\boldsymbol{c}_1,p_\text{LV}\right)=\mathbf{D}\left(t,\boldsymbol{c}_1,\binom{p_\text{LV}}{\overset{\sim}{\mathbf{W}}\left(t,\boldsymbol{c}_1,p_\text{LV}\right)}\right).
\end{equation*}\\
Finally, the coupled 3D--0D model is obtained, as reported in Equation~\eqref{3d0d-eq}, with the number 
of equations in the model matching the number of unknowns.

In this paper, a left ventricle geometry
is utilized, 
derived from the Zygote Solid 3D heart model \cite{mesh}; this model represents the 
50$^\text{th}$ percentile of a healthy Caucasian male in the U.S., reconstructed from 
a high-resolution computed tomography scan. The model is acquired at 70\% of diastole and
captures the late phase of diastasis when passive filling of the ventricles has slowed, 
but before atrial contraction completes active filling. The mesh used, denoted as 
$\mathcal{T}_{h}$, consists of 6079 vertices, 6682 triangles and 27492 tetrahedra, 
with an average element size $h_\text{mean}=3.7$ mm. The numerical methods, as 
briefly introduced above, are implemented within the \lifex \cite{lifex,lifex-home} framework in a parallel 
environment, that utilizes high-performance computing (HPC) resources (total available: 48 Intel Xeon ES-2640 CPUs) at 
MOX, Politecnico di Milano.
For more details about \lifex, its documentation can be retrieved 
\href{https://lifex.gitlab.io/lifex/}{here}, while refer to \cite{lifex-2024} 
for the last updates.

\section{Computational Models in Hypertension}
\label{hyper-model}
In this section, the computational models used to investigate hypertensive conditions 
are analyzed, focusing on their applicability in 
reproducing clinically relevant haemodynamic behaviours.

\subsection{Sensitivity Analysis in 0D Models}
\label{subsec:sensitivity}%
The sensitivity analysis of the 0D models is performed to identify the most influential 
parameters on the 
outputs of interest. This analysis is conducted using single- and multi-parameter
approaches, respectively in Sections~\ref{subsubsec:sobol} and \ref{subsec:multi-sensitivity}.

\subsubsection{Single-Parameter Sensitivity Analysis by means of Sobol Indices}
\label{subsubsec:sobol}%
Sobol indices \cite{sensitivity-analysis08} are a form of global sensitivity analysis: 
they rely on the decomposition of
the variance of the model's output into components, which correspond to different input
variables and their interactions.
In this paper, 
total Sobol indices $\mathcal{S}^{j,\mathcal{T}}_k$ have been employed: they take into consideration 
both the direct
effect of a parameter $p_k$ and its interactions with the others.
The methodology and analysis of this section are directly derived from \cite{sobol-tonini}.

\begin{figure}[t!]
    \centering
    \includegraphics[width=165 mm]{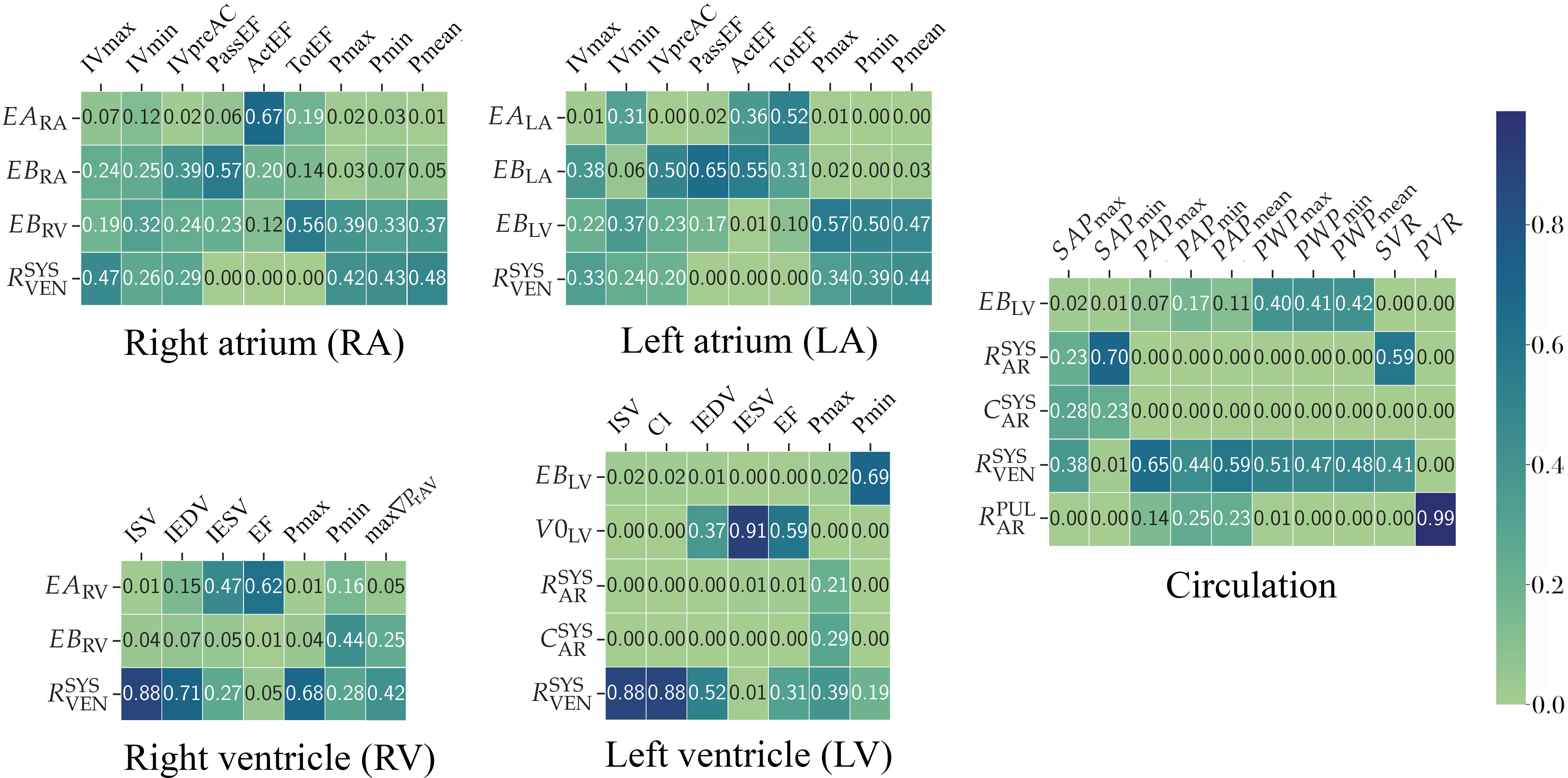}
    \caption{Most relevant total 
    Sobol indices, computed for $(\mathscr{C}_\text{NC})$. The index 
    $\mathcal{S}_k^{j,\mathcal{T}}$ in position $(k,j)$ quantify the contribution of 
    parameter $p_k$ to the output $y_j$. Only $\mathcal{S}_k^{j,\mathcal{T}}\ge 0.2$ are shown.}
    \label{fig:sobol-nc}
\end{figure}

The 0D cardiocirculatory model includes a large number of parameters, 
however, for practical and methodological reasons, some parameters are 
not included in the sensitivity analysis: 
HR is excluded, as can be 
directly and clinically measured; 
parameters related to the
timings of the cardiac cycle are excluded due to their significant influence on the
shape of pressure-volume (PV) loops, as they can 
lead to
non-physical results, and 
complicate the analysis without providing further
insights.
Moreover, the total blood volume was not included among the
parameters, since in this framework it is implicitly determined by the choice of initial
conditions rather than treated as an independent variable.
For both $(\mathscr{C}_\text{NC})$ and $(\mathscr{C}_\text{C})$, 
quantification of the influence of parameters and interpretation of their impact on 
model outputs follow the same principles and hence makes both configurations compatible 
and comparable.

The $\text{N}_p$ starting parameters, denoted as $\mathbf{p}^\text{R}\in\mathbb{R}^{\text{N}_p}$ (\ref{appendix:param-healthy}), are selected to be those which generate a 
simulation for a healthy individual.
To simplify the analysis, the parameters are assumed to be independent; furthermore, 
to reflect the diversity across different healthy and hypertensive individuals, the 
parameters are allowed to vary in a
multidimensional space, a hypercube with its center in
$\mathbf{p}^\text{R}$. The amount of variation along each axis is
$\frac{2}{3}$ of the reference parameter value, therefore each parameter $p_k$
is sampled within the interval:
\begin{equation*}
    \text{I}_k=\left[\left(1-\dfrac{2}{3}\right)p^\text{R}_k, \left(1+\dfrac{2}{3}\right)p^\text{R}_k\right], \quad k\in\mathbb{N} \text{ and }1 \leq k \leq \text{N}_p,
\end{equation*}
ensuring the representation of various physiological conditions.

In the parameter sampling process, Sobol sequences are used 
to ensure a more uniform and space-filling coverage of the hypercube, 
by minimizing clustering and gaps. This quasi-random structure allows for faster 
convergence of the sample mean to the expected value compared to standard Monte Carlo 
sampling, thereby reducing the number of simulations needed to explore the parameter 
space effectively \cite{sobol-sequence:discrepancy}.
Sobol indices are computed using Saltelli's method \cite{saltelli-method}:
this technique guarantees a high-resolution variance decomposition, providing an 
accurate information
on the contribution of each parameter and its interactions.
Though computational effort of Saltelli's method
scales linearly with the number of parameters, it is still feasible, since
$(\mathscr{C}_\text{NC})$ and $(\mathscr{C}_\text{C})$ are computationally cheap, and 
parallelization limits the time taken for processing.
In this paper, 
$\text{N}_p = 26$ parameters 
are analyzed for $(\mathscr{C}_\text{NC})$, while, in the case of $(\mathscr{C}_\text{C})$, 
the number of parameters taken in account rises up to 32; the method 
generates $2\text{N}\left(\text{N}_p + 1\right)$ 
samples ($\text{N}= 2^{12}$ is a user-defined value), resulting in 221,184 samples 
for $(\mathscr{C}_\text{NC})$ and 270,336 samples for $(\mathscr{C}_\text{C})$. 
The most relevant total Sobol indices $(\ge 0.2)$ for $(\mathscr{C}_\text{NC})$, 
obtained through this procedure and subsequently utilized, are shown 
in Figure~\ref{fig:sobol-nc}.

\subsubsection{Multi-Parameter Sensitivity Analysis}
\label{subsec:multi-sensitivity}%
Multi-parameter sensitivity analysis extends the traditional sensitivity 
analysis in order to understand the combined effects that related parameters have on 
model outputs: 
this approach 
aggregates them into meaningful groups, often based on physiological, structural, or 
functional similarities. This approach is particularly beneficial for complex models 
characterized by a large number of parameters, whose collective influence
can not generally be interpreted in an elementary way 
with respect to their individual influence on the model. 

In this section, groups of parameters (specifically, resistances, inductances and capacitances) related to specific parts of the cardiocirculatory 
system are considered in order to understand how changes in those groups produce 
variations in the model's behaviour, in both the systemic and pulmonary circulation 
cases. Namely, the impact on the simulation of the variation of arterial, venous and 
capillary parameters is studied; 
the parameters related to the left ventricle 
($\text{E}^a_\text{LV}$, $\text{E}^p_\text{LV}$ and $\text{V}_{0,\text{LV}}$) are 
analyzed to assess how keeping them fixed in the 3D--0D model may limit its ability to 
accurately reproduce certain hemodynamic conditions 
(Sections~\ref{subsec:3d-0d-model}). The goal is to understand whether 
fixing these parameters hinders the model from capturing relevant physiological changes, 
and to what extent this choice influences the simulation outcomes.
A noteworthy aspect is the arbitrary choice of groups of 
parameters; another drawback of the current approach is the fact that other groups of 
possibly influential parameters may not be considered at all.
The sensitivity analysis of these groups of parameters is performed for both 
$(\mathscr{C}_\text{NC})$
and $(\mathscr{C}_\text{C})$; 
for $(\mathscr{C}_\text{C})$, the systemic
and pulmonary capillary circulation is considered separately.

All the mentioned parameters undergo modification by the multiplication with 
coefficients in such a way that they are constrained within realistic intervals, while 
all the others remain unchanged. 
The modeling of hypertension by means of lumped-parameter models is a topic that 
has not been systematically exploited yet; different 
studies employ different models, each with different structures and peculiarities, such 
as in \cite{coeff-sens-analys1,coeff-sens-analys2,coeff-sens-analys3}. 
Since our sensitivity analysis focuses on the hypertensive setting, 
constructing realistic parameter intervals is not a trivial task; in this 
paper, inspiration has been taken from 
\cite{coeff-sens-analys1,coeff-sens-analys2,coeff-sens-analys3} and 
Section~\ref{hyper}, in order to obtain coefficients 
that would be acceptable and useful for carrying out the sensitivity analysis.
These coefficients have been chosen to reflect the typical parameter variations 
observed during the progression of hypertension: 
resistances and 
inductances are increased while conductances are decreased, and vice versa; 
active and passive elastance of the left ventricle are enhanced while the resting volume 
is reduced, and vice versa.
Such modifications try to reproduce the haemodynamic alterations reported in clinical and 
modeling studies, while keeping the parameters within physiologically plausible ranges. 
The coefficients are defined as follows:
The coefficients are defined as follows:
\begin{subequations}\label{sobol-analysis:coeff-to-change}
\begin{align}
        \eta_\text{R}&=1 + 0.15 \rho,  & \eta_\text{L}&=1+0.075\rho, & \eta_\text{C}&=1-0.06\rho,\label{sobol-analysis1}\\
        \eta_a&=1 + 0.12 \rho, &\eta_p&=1 + 0.14 \rho, &\eta_0&=1 - 0.05 \rho,\label{sobol-analysis2}
\end{align} 
\end{subequations}
where $\eta_\text{R}$, $\eta_\text{L}$ and $\eta_\text{C}$ are the coefficients 
applied to resistances, inductances and 
conductances, respectively; $\eta_a$, $\eta_p$ and $\eta_0$ represent the coefficients 
applied to active 
elastance, passive elastance and resting volume of left ventricle, respectively; 
the auxiliary parameter $\rho$ is defined as 
$\rho\in\mathbb{Z}\backslash\{0\}$, with the constraint 
$\left|\rho\right|\le 5$.

In this section, to simplify the notation, given a generic parameter or output $\alpha$ 
(time dependent or independent), the symbol $\bar{\alpha}$ will denote the reference 
value of the parameter 
or the output computed using the reference parameters, while $\hat{\alpha}$ will 
indicate the value of the parameter after being modified according to the guidelines 
outlined earlier, or the output computed with the modified parameters.
Therefore, the following relation holds: $\hat{\alpha}=\eta_\alpha \bar{\alpha}$, 
for some parameter $\alpha$ and the corresponding coefficient $\eta_\alpha$.

The time-independent outputs monitored in this analysis are listed in 
\ref{appendix:ranges-healthy}.
The variation of an output $\zeta$ is 
quantified using the following formula:
\begin{equation}\label{indicator-out}
    \Gamma(\zeta)=\left(\dfrac{\hat{\zeta}}{\raisebox{-0.2ex}[0pt][0pt]{$\bar{\zeta}$}}-1\right)\cdot 100.
\end{equation}
The variables that represent the solutions to the ODE system (Section~\ref{subsec:lumped}) within the 
simulation are also monitored, along with the pressures in the heart chambers, the 
fluxes through the cardiac valves and, only for $(\mathscr{C}_\text{C})$, the non 
oxygenated pulmonary capillary flux. 
Since these variables are functions of time (represented as vectors in the numerical 
solution), four indicators are introduced to 
simplify the monitoring process and provide scalar values for analysis. They are 
computed for each variable $v$ as follows:
\begin{equation}\label{indicator-variable}
    \begin{aligned}
    &\Gamma_\text{T}(v)=\left(\dfrac{\int_{\text{T}-\text{T}_\text{HB}}^{\text{T}} (\hat{v}-\bar{v})\,dt}{\int_{\text{T}-\text{T}_\text{HB}}^{\text{T}} \bar{v}\,dt}\right)\cdot 100,
    &\bar{\Gamma}(v)=\left(\dfrac{\text{mean}_t\,{\hat{v}}}{\text{mean}_t\,{\bar{v}}}-1\right)\cdot 100,\\
    &\Gamma_\text{M}(v)=\left(\dfrac{\max_t{\hat{v}}}{\max_t{\bar{v}}}-1\right)\cdot 100,
    &\Gamma_\text{m}(v)=\left(\dfrac{\min_t{\hat{v}}}{\min_t{\bar{v}}}-1\right)\cdot 100.
\end{aligned}
\end{equation}
Additionally,
the variation in the PV loops $\Gamma_i$, for 
$i\in\left\{\text{LA},\text{LV},\text{RA},\text{RV}\right\}$, is measured as the ratio 
between the area of the cycle obtained with modified parameters and the area of the 
cycle obtained with unmodified parameters.

As the analysis begins, the parameters of a group
are modified by means of the coefficients in 
Equation~\ref{sobol-analysis:coeff-to-change}, while all the other parameters remain 
fixed at reference values. When a simulation is performed, a certain number of 
time-independent outputs and time-dependent variables are computed: in particular, 
45 outputs and 20 variables are computed if 
$(\mathscr{C}_\text{NC})$
is employed, while 46 outputs and 25 variables are 
calculated if $(\mathscr{C}_\text{C})$ is utilized (refer 
to \ref{appendix:ranges-healthy} for the outputs, to 
Section~\ref{subsec:lumped} for the variables).
Additionally, in the case of $(\mathscr{C}_\text{NC})$, five groups are 
considered, while seven groups are taken for
$(\mathscr{C}_\text{C})$.
By using the coefficients in Equation~\eqref{sobol-analysis:coeff-to-change}, 
both models' parameters are changed over ten nuances of variation for each 
considered group: whence, 50 simulations are run for $(\mathscr{C}_\text{NC})$ 
and 70 simulations for $(\mathscr{C}_\text{C})$.
For each output, the indicator in Equation~\eqref{indicator-out}
is calculated, 
while, for each variable, the four indicators in Equation~\eqref{indicator-variable}
are computed; four additional indicators ($\Gamma_\text{LA}$, $\Gamma_\text{LV}$, 
$\Gamma_\text{RA}$, $\Gamma_\text{RV}$) are also considered, once for
each nuance 
of variation of each group.
All in all, a simulation employing $(\mathscr{C}_\text{NC})$ produces 129 
indicators and a simulation adopting $(\mathscr{C}_\text{C})$ generates
150 indicators, which 
gives an overall total of
16,950 indicators across all simulations, taking in account every group, every 
nuance and every model.

At the end of the analysis, the parameters to calibrate comprehend: the active and passive elastances, and the 
unstressed volume of the four cardiac chambers, the valve maximal and minimal resistances 
and all the RLC components of the vascular circulation (Section~\ref{subsec:lumped}).

\subsection{Implementation of Hypertensive Scenarios}
\label{subsec:implem-hyper}%
In order to simulate the impact of hypertension in the 0D cardiocirculatory model, 
it is necessary to introduce systematic modifications of model parameters. Changes in 
parameters are implemented in such a way that simulated alterations in haemodynamics 
align with clinical and physiological expectations: this would present a 
realistic reproduction of hypertensive pathophysiology within the framework of 0D 
modeling. The modifications are done in the form of percentage changes relative to 
baseline parameters of a healthy individual, 
collected in \ref{appendix:param-healthy}: all those adjustments try to represent
the physiological processes explained in Section~\ref{hyper}, maintaining a 
consistent representation of hypertensive condition.
The result obtained from sensitivity analysis
(Section~\ref{subsec:sensitivity}) are supporting guidelines to enhance the modifications, 
helping to achieve the desired haemodynamic effects within the model.

In this study, three distinct types of hypertension are considered 
(Section~\ref{sec:cause_sympt}), 
each characterized by different haemodynamic implications and clinical relevance:
\begin{enumerate}
\item systemic hypertension, mostly affecting the systemic arterial
circulation and commonly associated with increased peripheral resistance;
\item pulmonary hypertension, involving elevated pressure in the pulmonary
circulation;
\item renovascular hypertension, due to abnormalities in kidneys, which gives rise to 
compensatory mechanisms, further increasing systemic blood pressure; the
examined condition is also aggravated by secondary pulmonary hypertension.
\end{enumerate}
For each of these conditions, three different levels of severity (mild, moderate and
severe) are investigated, in an effort to capture the progressive nature of the disease 
and its effects on the circulatory dynamics; the choice of these three levels of 
severity was made in such a way as to represent the range of hypertensive conditions 
that might be observed in a realistic patient population.
The parameter changes are listed in Table~\ref{table:hyper-changes}, 
indicating the parameter adjustments for each condition and severity level.

\begin{table}[t!]
    \centering
    \begin{tabular}{llllllllll}
    \hline
     & \multicolumn{3}{c}{Systemic} &\multicolumn{3}{c}{Pulmonary} &\multicolumn{3}{c}{Renovascular} \\
    \cmidrule(lr){2-4} \cmidrule(lr){5-7} \cmidrule(lr){8-10}
                           Parameter & Mild & Moderate & Severe
                           & Mild & Moderate & Severe
                           & Mild & Moderate & Severe \\
    \hline 
    $\text{HR}$ & $-$ & $+5\%$ & $+10\%$                            & $-$ & $+5\%$ & $+10\%$                & $-$ & $-$ & $+10\%$\\[0.2em]
    $\text{E}^a_\text{LV}$ & $+10\%$ & $+20\%$ & $+40\%$            & $+5\%$ & $+10\%$ & $+15\%$            & $-$ & $-$ & $-$\\[0.2em]
    $\text{V}_{0,\text{LV}}$ & $-$ & $-$ & $-10\%$                  & $-$ & $-$ & $-10\%$                   & $-$ & $-$ & $-$\\[0.2em]
    $\text{R}_\text{AR}^\text{SYS}$ & $+10\%$ & $+30\%$ & $+50\%$   & $-$ & $-$ & $+5\%$                    & $+15\%$ & $+40\%$ & $+85\%$\\[0.2em]
    $\text{C}_\text{AR}^\text{SYS}$ & $-10\%$ & $-20\%$ & $-40\%$   & $-$ & $-$ & $-10\%$                   & $-10\%$ & $-15\%$ & $-25\%$\\[0.2em]
    $\text{R}_\text{VEN}^\text{SYS}$ & $-$ & $-$ & $+5\%$           & $-$ & $-$ & $-$                       & $-$ & $-$ & $-$\\[0.2em]
    $\text{C}_\text{VEN}^\text{SYS}$ & $-$ & $-$ & $-5\%$           & $-$ & $-$ & $-$                       & $-$ & $-$ & $-$\\[0.2em]
    $\text{R}_\text{AR}^\text{PUL}$ & $-$ & $-$ & $+10\%$           & $+20\%$ & $+50\%$ & $+100\%$          & $+5\%$ & $+15\%$ & $+40\%$\\[0.2em]
    $\text{C}_\text{AR}^\text{PUL}$ & $-$ & $-$ & $-10\%$           & $-10\%$ & $-25\%$ & $-50\%$           & $-$ & $-$ & $-$\\[0.2em]
    $\text{R}_\text{VEN}^\text{PUL}$ & $-$ & $-$ & $-$              & $+10\%$ & $+25\%$ & $+50\%$           & $-$ & $-$ & $-$\\[0.2em]
    $\text{C}_\text{VEN}^\text{PUL}$ & $-$ & $-$ & $-$              & $-10\%$ & $-25\%$ & $-50\%$           & $-$ & $-$ & $-$\\[0.2em]
    $\text{R}_\text{C}^\text{SYS}$ & $-$ & $-$ & $+10\%$            & $-$ & $-$ & $-$                       & $-$ & $-$ & $-$\\[0.2em]
    $\text{C}_\text{C}^\text{SYS}$ & $-$ & $-$ & $-20\%$            & $-$ & $-$ & $-$                       & $-$ & $-$ & $-$\\[0.2em]
    $\text{R}_\text{C}^\text{PUL}$ & $-$ & $-$ & $-$                & $+10\%$ & $+30\%$ & $+60\%$           & $+5\%$ & $+10\%$ & $+25\%$\\[0.2em]
    $\text{C}_\text{C}^\text{PUL}$ & $-$ & $-$ & $-$                & $-$ & $-$ & $-20\%$                   & $-$ & $-$ & $-$\\[0.2em]
    $R_\text{AV}$ & $+25\%$ & $+50\%$ & $+100\%$                    & $-$ & $-$ & $-$                       & $+25\%$ & $+50\%$ & $+100\%$\\[0.2em]
    $R_\text{MV}$ & $-$ & $+10\%$ & $+20\%$                         & $-$ & $-$ & $-$                       & $+10\%$ & $+20\%$ & $+30\%$\\[0.2em]
    $R_\text{PV}$ & $-$ & $-$ & $-$                                 & $+25\%$ & $+50\%$ & $+100\%$          & $+25\%$ & $+50\%$ & $+100\%$\\[0.2em]
    $R_\text{TV}$ & $-$ & $-$ & $-$                                 & $+10\%$ & $+25\%$ & $+50\%$           & $+10\%$ & $+25\%$ & $+50\%$\\[0.2em]
    \hline
    \end{tabular}
    \caption{Modifications 
        applied to the model parameters to simulate mild, moderate and severe 
        hypertension. The changes are expressed as percentage variations relative to the 
        baseline parameters of a healthy individual, which are reported in 
        \ref{appendix:param-healthy}.}
    \label{table:hyper-changes}
\end{table}

\subsection{Models Calibration}
\label{subsec:calibration}%
In this section, two different calibration procedures for the 0D and the 3D--0D models 
are presented, respectively in Sections~\ref{subsubsec:calibration:0d} 
and \ref{subsubsec:calibration:3d}; since the models are profoundly different also
in terms of the computational cost required to solve them, it was deemed appropriate 
to adopt two different calibration strategies. 

\subsubsection{Calibration of the 0D model}
\label{subsubsec:calibration:0d}%
The calibration procedure that is presented in this section follows the framework 
described in \cite{lumped:capillary}, which, in turn, relies on the method presented in 
\cite{calibration-model}. The process begins by defining a loss function 
that quantifies 
the discrepancy between simulation results and target values. 
The aim of calibration is identifying the parameter configuration 
$\bar{\mathbf{p}}^j\in \Theta=\mathbb{R}^{\text{N}_p}$, with $\text{N}_p$ being the number 
of parameters to calibrate, 
that minimizes the loss function; in order to find the configuration $\bar{\mathbf{p}}^j$, 
parameters are iteratively modified within predefined bounds, untill the loss function 
is sufficiently small.
In \cite{lumped:capillary}, the loss function is given by the 
sum of squared relative errors between 
model outputs and target data, for the specific individual $j$, as follows:
\begin{equation}\label{loss-p-opt}
    \bar{\mathbf{p}}^j=\underset{\mathbf{p}\in\Theta}{\argmin} \,\,\mathcal{L}^j(\mathbf{p}),
\end{equation}
with the relative error $\delta^j_{k}$ and the loss function $\mathcal{L}^j$, defined as:
    \begin{align}\label{loss-function-article}
    \delta^j_{k}(\mathbf{p})&=\left|\dfrac{q_{\text{m}_j(k)}(\mathbf{p}) - d^j_k}{d_k^j}\right|,    & \mathcal{L}^j(\mathbf{p}) &= \sum_{k=1}^{\text{N}_j} \delta^j_{k}(\mathbf{p})^2,
    \end{align}
where $\text{N}_j$ is the number of available data for individual $j$, 
$d_k^j$ represents the $k$-th target data value of individual $j$, and
$q^j_{\text{m}_j(k)}(\mathbf{p})$ is the corresponding model output. 
The index m for $q_\text{m}^j$ is $\text{m}\in\{1, \dots, 
\text{N}_q\}$, where $\text{N}_q$ is the total number of distinct model outputs. 
For more details about parameters $\mathbf{p}$ to be calibrated, refer to Section~\ref{subsec:sensitivity}.

To minimize $\mathcal{L}^j$, parameters are optimized within predefined intervals 
$\text{I}_i$, where $i\in\{ 1, \dots, \text{N}_p\}$. These intervals are set around some initial reference values 
in such a way that, after the optimization, the changes in parameters remain physically 
and clinically relevant.
The value of the loss function falling below some given small threshold, set to $10^{-3}$, 
is assumed to indicate a successful calibration.
Nevertheless, there are cases where the procedure might be unsuccessful, such as when 
the loss function's 
minimum is located above the threshold, meaning the parameter set is not able to reproduce the 
reference data within the given constraints.
The minimization of $\mathcal{L}^j$ is performed in Python using the Quasi-Newton 
optimization algorithm L-BFGS-B \cite{bfgs}, as implemented in the \texttt{SciPy} library 
\cite{scipy}.

\begin{figure}[t!]
    \centering
    \subfloat[Loss functions in linear scale.\label{fig:loss-functions}]{
        \includegraphics[width=0.475\linewidth]{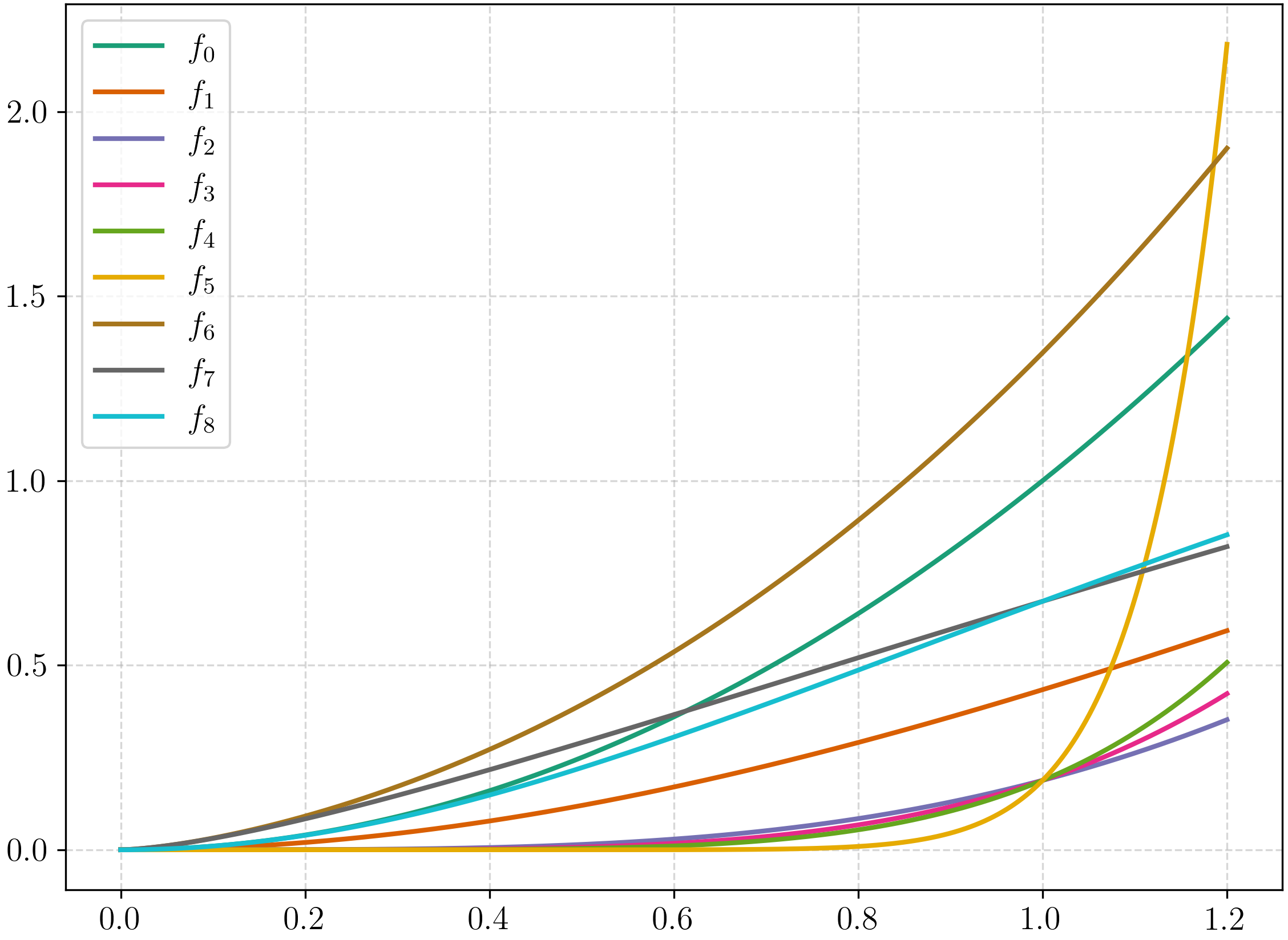}
    }
    \hfill
    \subfloat[Loss functions in semi-logarithmic scale.\label{fig:loss-functions-log}]{
        \includegraphics[width=0.475\linewidth]{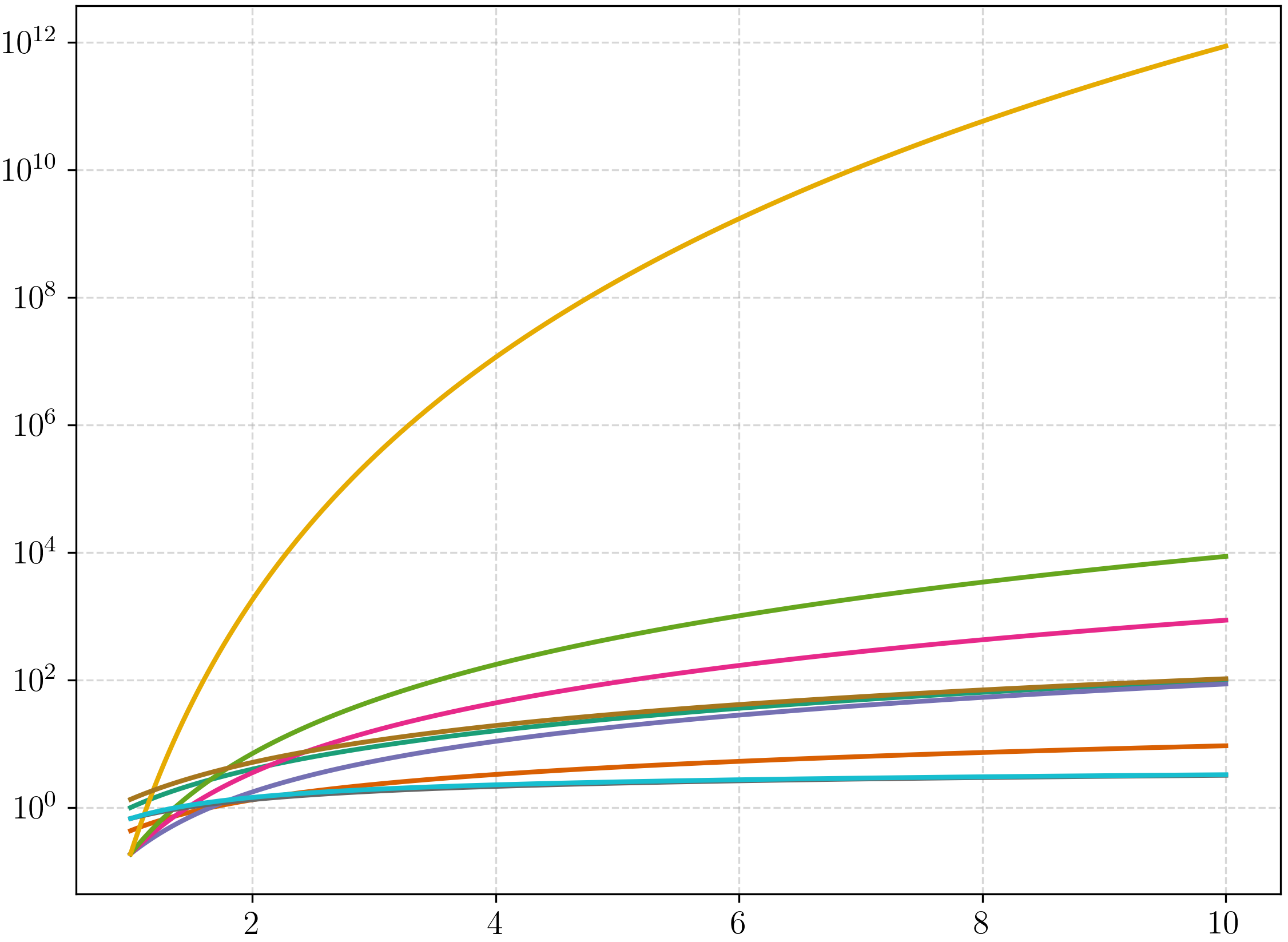}
    }
    \caption{Comparison of the eight loss functions used in this study: 
    on the left, they are represented in linear scale for $\delta\in[0, 1.2]$,
    while on the right in semi-logarithmic scale on $y$ axis for $\delta\in[1, 10]$.}
    \label{fig:loss}
\end{figure}

\paragraph{An Alternative Choice for the Loss Function.}
\label{paragraph:loss-chioce}%
Apart from the choice of an optimization algorithm, another important issue in model 
calibration is related to the selection of a suitable loss function, which drives the 
optimization procedure. 
The selection or design of an 
appropriate loss function is particularly critical in the context of cardiocirculatory 
models, where complex physiological dynamics must be captured with high fidelity during 
the calibration process.
The method is based on constructing alternative loss functions that apply distinct 
penalizations to different magnitudes of parameter values;
subsequently, these various loss functions are 
applied for the calibration process so that the results of the estimated parameters and 
the associated calibration performance could be compared.
Since the following analysis does not focus on the individual from whom the data is 
provided but rather the functional form of the loss function, the index $j$ in 
Equations~\eqref{loss-p-opt} and \eqref{loss-function-article} will
be omitted to simplify the notation; additionally, the index $r$ and the total 
parameter deviation $\Delta_{r}$ are introduced as follows:
\begin{subequations}\label{loss-f_r}
    \begin{align}
        \mathcal{L}_r(\mathbf{p}) &= \sum_{k=1}^{\text{N}} f_r\left(\delta_{k}(\mathbf{p})\right),
        &\Delta_{r}(\mathbf{p}) &= \sum_{k=1}^{\text{N}} \delta_{k}(\mathbf{p}),
    \end{align}
\end{subequations}
such that $f_r$ denotes a different function as $r$ changes in the subsequent analysis.
By defining the following functions:
\begin{align*}
    f_0(\delta)&=\delta^2, & g(\delta)&=\log{(\cosh{\delta})}, & h(\delta)&=1+\dfrac{1}{2}\log{\left(1+\delta^2\right)},
\end{align*}
it is possible to express the different functions $f_r$ besides $f_0$ (Figure~\ref{fig:loss}), chosen to be subsequently analyzed,
as follows:

\begin{subequations}\label{loss-allfunc}
    \begin{align}
        f_1 (\delta)&=  g(\delta), & f_2 (\delta)&= g(\delta)^2,  & f_3 (\delta)&= \left|\delta\right|\,g(\delta)^2,  & f_4 (\delta)&= \delta^2 \,g(\delta)^2, & f_5 (\delta)&= \delta^{10} \,g(\delta)^2,\label{loss-log2-x10}\\
        f_6 (\delta)&=  |\delta|^{\frac{3}{2}}\,h(\delta), & f_7 (\delta)&= \dfrac{|\delta|^{\frac{3}{2}}}{1+|\delta|^{\frac{3}{2}}}\,h(\delta), & f_8 (\delta)&= \dfrac{\delta^{2}}{1+\delta^{2}}\,h(\delta). \label{loss-frac-x2}
    \end{align}
\end{subequations}

\begin{figure}[t!]
    \centering
    \subfloat[Systemic hypertension.\label{fig:cal:time:sys}]{
        \includegraphics[width=0.47\linewidth]{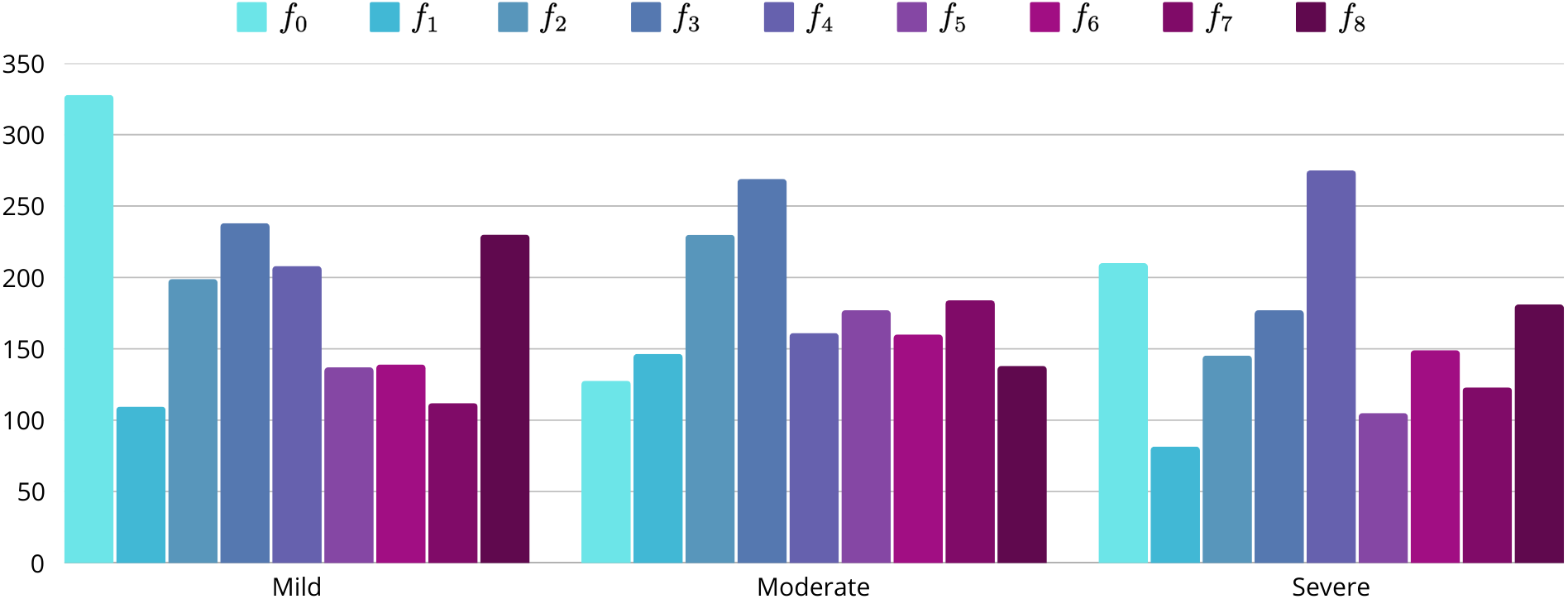}
    }\hfill
    \subfloat[Systemic hypertension.\label{fig:cal:loss:sys}]{
        \includegraphics[width=0.47\linewidth]{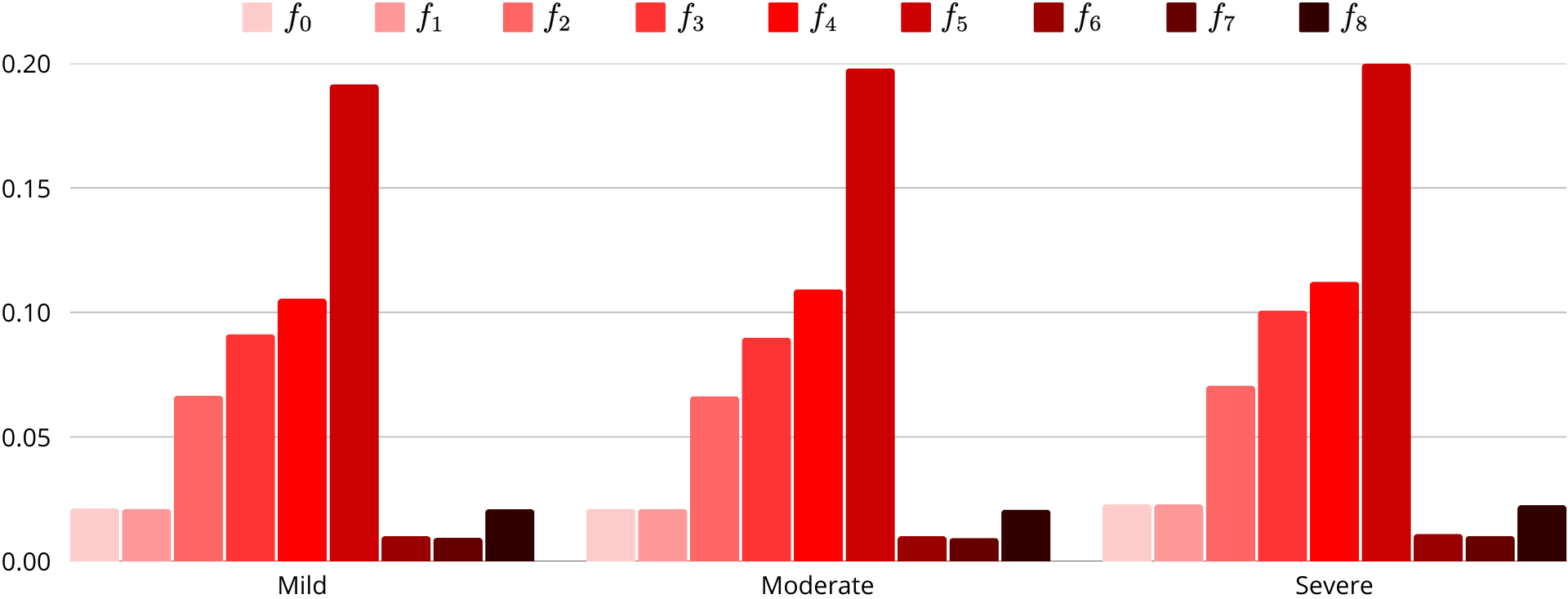}
    }\\[0.4cm]
    \subfloat[Pulmonary hypertension.\label{fig:cal:time:pul}]{
        \includegraphics[width=0.47\linewidth]{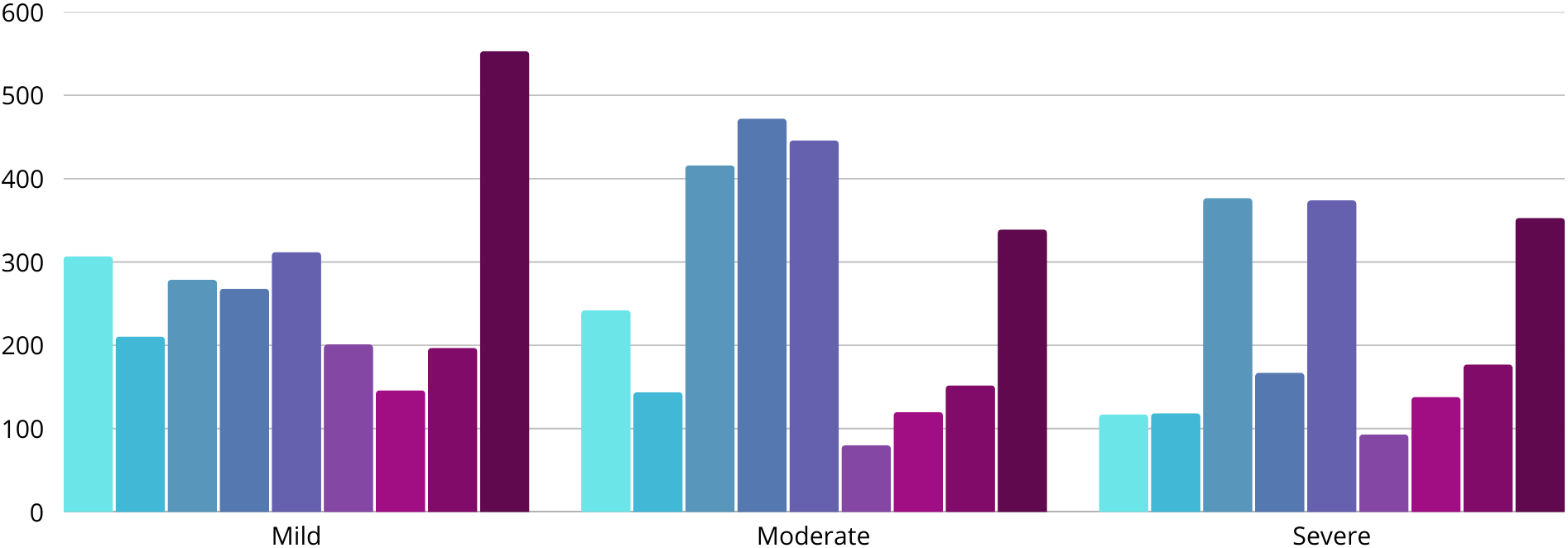}
    }\hfill
    \subfloat[Pulmonary hypertension.\label{fig:cal:loss:pul}]{
        \includegraphics[width=0.47\linewidth]{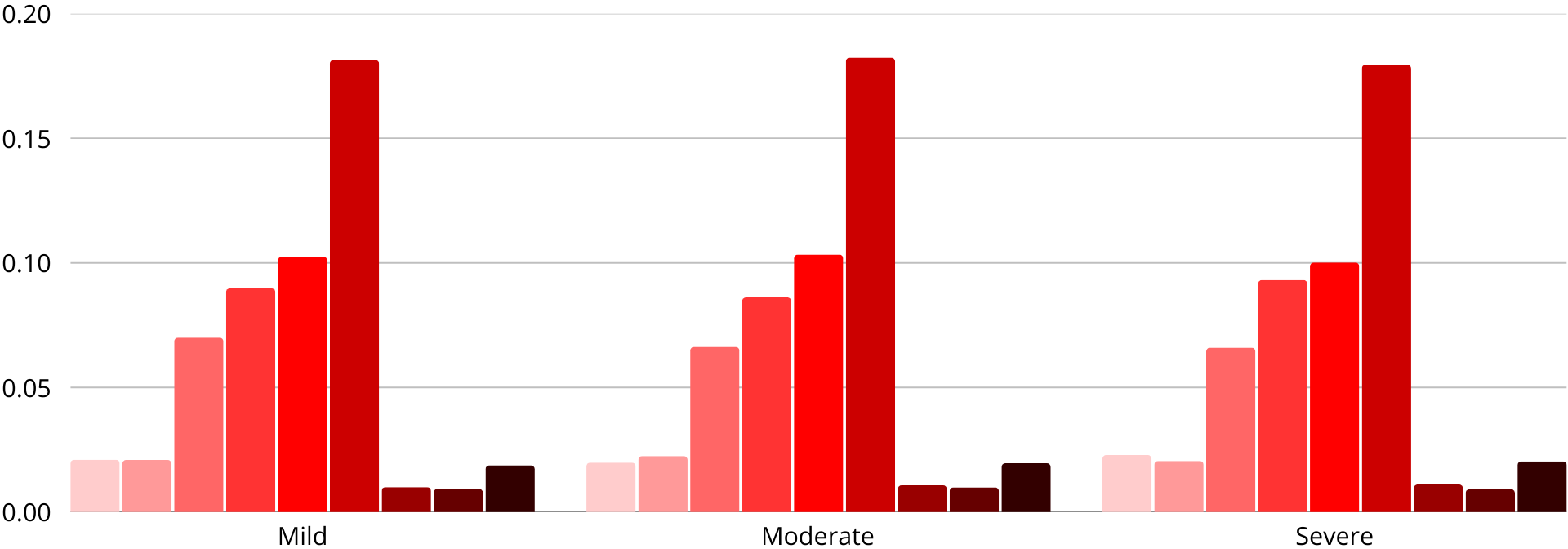}
    }\\[0.4cm]
    \subfloat[Renovascular hypertension with secondary pulmonary hypertension.\label{fig:cal:time:reno}]{
        \includegraphics[width=0.47\linewidth]{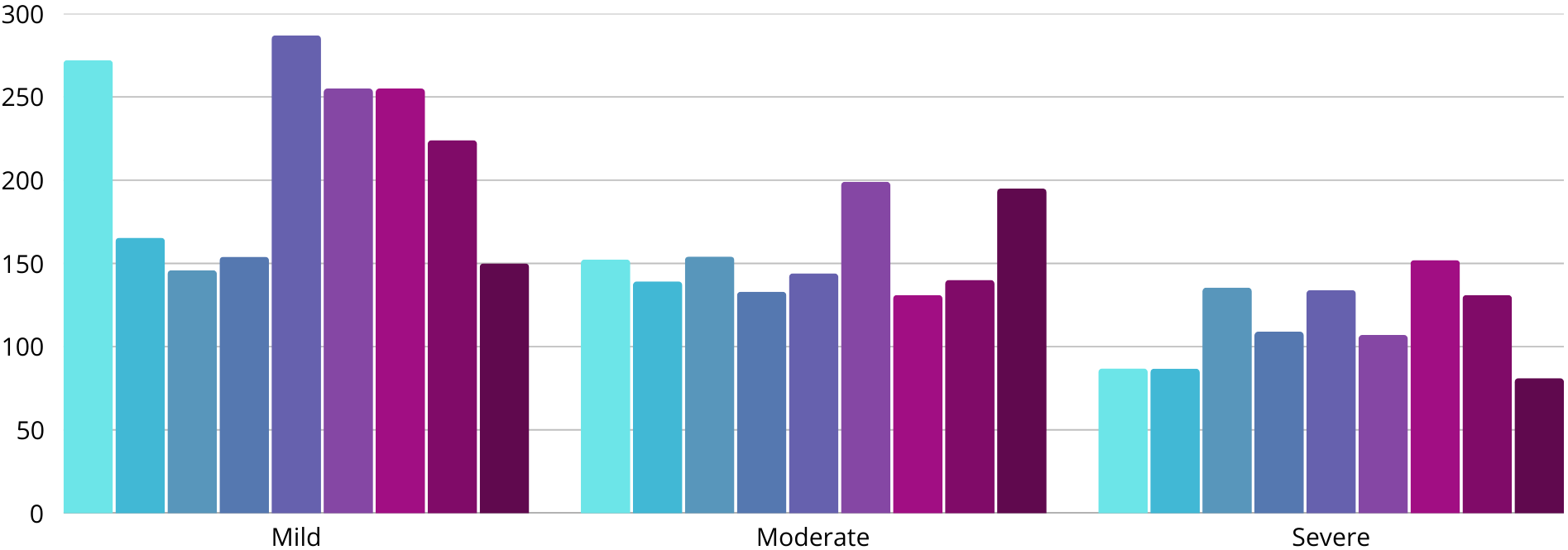}
    }\hfill
    \subfloat[Renovascular hypertension with secondary pulmonary hypertension.\label{fig:cal:loss:reno}]{
        \includegraphics[width=0.47\linewidth]{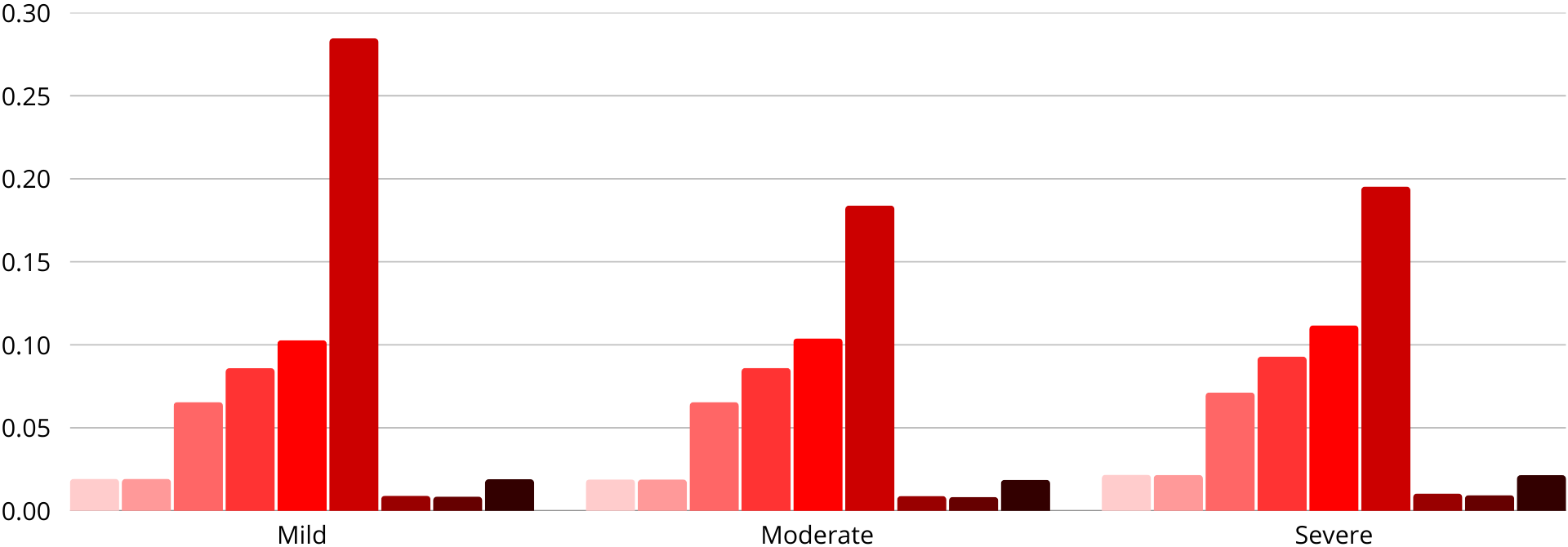}
    }
    \caption{Calibration 
    times (in seconds, on the left) and total parameter deviations $\Delta_{r}$ (on the right) at the end of calibration 
    for $(\mathscr{C}_\text{NC})$. The following conditions are 
    studied: systemic, pulmonary and renovascular hypertension 
    with secondary pulmonary hypertension. Three severity levels (mild, moderate and 
    severe) are considered for each condition, and the calibration is performed using 
    different choices for the loss function.}
    \label{fig:calibration}
\end{figure}

To investigate the impact of the different choices of $f_r$ in the loss function 
definition in Equation~\eqref{loss-f_r}, a systematic 
study has been conducted. In the absence of real clinical data, a synthetic approach is 
adopted: data from a healthy individual (\ref{appendix:param-healthy}) are used 
as a baseline, and modifications 
are introduced to replicate the causes and effects of hypertension. These modified 
datasets are employed as input for simulations performed using the $(\mathscr{C}_\text{C})$ 
model, 
that has been discussed in Section~\ref{subsec:lumped}; 
time-independent outputs from these simulations are then collected to serve as reference data.
The next step involves calibrating $(\mathscr{C}_\text{NC})$, presented 
in Section~\ref{subsec:lumped}, using these time-independent 
outputs as targets: 
for the selection of the outputs utilized in the calibration process, all the outputs 
included in \ref{appendix:ranges-healthy}
are considered. 
As for the parameters to be calibrated, 
reference is made to the total Sobol indices $\mathcal{S}^\mathcal{T}$ presented in 
Figure~\ref{fig:sobol-nc}: 
specifically, only those parameters $p_k$ such that 
$\exists j: \mathcal{S}^{j,\mathcal{T}}_k\ge 0.2$
are selected for calibration.

This setup is intentionally designed to enable a focused comparison of different loss 
functions, without introducing additional variability related to capillary dynamics in 
the model being calibrated. Using $(\mathscr{C}_\text{C})$ to generate the reference 
data and calibrating $(\mathscr{C}_\text{NC})$ ensures that the evaluation of the loss 
functions is not influenced by the added complexity of capillary parameters.
Moreover, $(\mathscr{C}_\text{C})$ is not calibrated anywhere in this study, both 
because the calibration procedure is more challenging and because it tends to produce 
parameter values that result in less physiologically realistic cardiac and circulatory 
dynamics compared to those obtained with $(\mathscr{C}_\text{NC})$.

For each candidate loss function in Equation~\eqref{loss-allfunc}, the calibration 
process is carried out to determine the optimal parameter configuration, and the time needed 
is recorded; then, the vectorial outputs 
from the two models as in \ref{appendix:capillary} are compared. This procedure 
is repeated for the cases presented in Section~\ref{subsec:implem-hyper}, encompassing 
systemic hypertension, pulmonary 
hypertension and renovascular hypertension with secondary pulmonary hypertension. 
Each of these conditions is further studied at three different levels of gravity (mild, 
moderate and severe hypertension), allowing a thorough check of the loss functions in 
a wide variety of situations. The whole analysis has been performed on a MacBook Pro 
(Apple M1 Pro, 10 cores, 3.2 GHz, 16 GB RAM). Times necessary to complete the 
calibration for each analyzed case with a different choice of the loss function are 
summarized in Figure~\ref{fig:calibration}.
Among all the candidate functions, 
$f_1$ stands out as the 
optimal selection: even if the total parameter deviations in Figure~\ref{fig:calibration} 
are generally smaller for $f_6$ and $f_7$, $f_1$ requires the shortest 
calibration time, being the best compromise between quality of the parameters and 
calibration time.
As a result of this analysis, all subsequent calibrations are performed using 
$f_1$ as the function defining the loss; for simplicity, from this 
point onward, the loss function will be denoted as $\mathcal{L}$ in place of $\mathcal{L}_1$.

\subsubsection{Calibration of the 3D--0D model}
\label{subsubsec:calibration:3d}%
\begin{figure}[t]
    \centering
    \includegraphics[width=\linewidth]{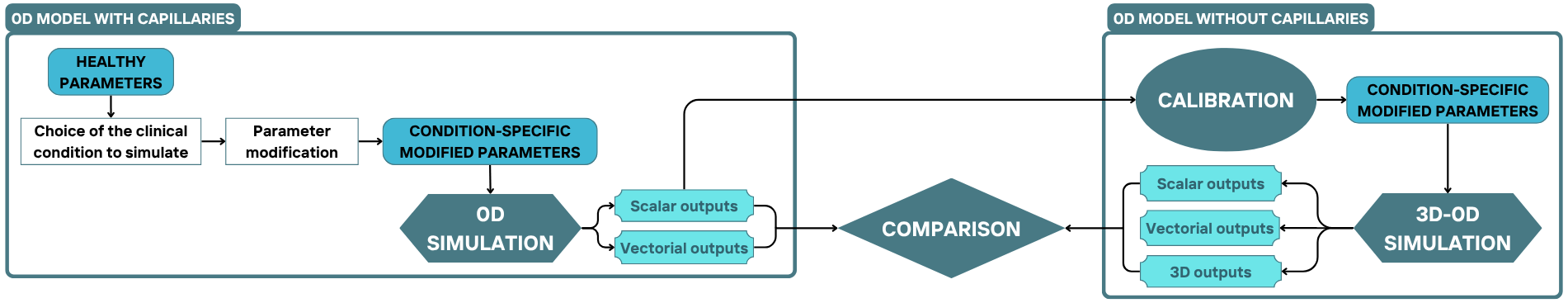}
    \caption{Workflow illustrating the calibration process of 3D--0D models for 
    a generic clinical scenario.}
    \label{fig:calibration:workflow}
\end{figure}

In general, the calibration of a 3D--0D model is computationally expensive, due 
to the high complexity of the three-dimensional representation of the left ventricle, 
making both single-parameter and multi-parameter sensitivity analysis unfeasible.
Even a direct calibration would be impractical, as it would require, 
among other things, solving the 3D--0D model many times, an operation that is computationally expensive.
In addition, a whole calibration process would be required for each case of 
hypertension that is to be analysed, nine in total 
(Section~\ref{subsec:implem-hyper}), making the approach even less feasible;
lastly, all the adjustments performed to simulate hypertension involve only parameters 
of the 0D circulation model (as detailed in Section~\ref{subsec:implem-hyper}). 
In light of all these considerations, an alternative and more
computationally efficient approach, depicted in Figure~\ref{fig:calibration:workflow}, 
is preferred.

The workflow starts by exploiting the 0D model with capillary circulation, 
$(\mathscr{C}_\text{C})$, which allows a direct manipulation on parameters related to 
capillaries: through modification of $(\mathscr{C}_\text{C})$, a direct action is 
exerted on the capillary circulation in order to ensure that, after the calibration of the 
electromechanical model, the 3D--0D representation also reflects the effects of 
hypertension at the capillary level, despite capillaries are not physically included in 
the 3D--0D model. 
Taking the parameters of a healthy individual 
(\ref{appendix:param-healthy}), the adjustments to simulate hypertension, 
as outlined in Section~\ref{subsec:implem-hyper}, are then applied to this model. 
A simulation is performed, and its time-independent outputs are used to calibrate, 
following the procedure illustrated in Section~\ref{subsubsec:calibration:0d}, a 0D model 
without capillary circulation, $(\mathscr{C}_\text{NC})$: such a procedure allows for 
direct manipulations of the capillary network and ensures that $(\mathscr{C}_\text{NC})$
properly captures the effects produced by these manipulations.

It is important to remark here that all the parameters related to the left ventricle, 
which are active and passive elastance and the unstressed volume, are not included in 
this calibration procedure: such a choice is indeed motivated by the fact that, in the 
3D--0D framework, the left ventricle is represented by means of a 3D electromechanical 
structure, as described in Section~\ref{subsec:3d-0d-model}, which is inherently not 
including these latter parameters; 
the mesh employed in the 3D--0D model is generated from the left ventricle 
of a healthy individual, as explained in Section~\ref{subsec:3d-0d-model}.
It would hence 
be inconsistent 
to modify the 0D left ventricle parameters because changes in such parameters can not be 
replicated in the 3D electromechanical model. For additional information, refer to the 
details about the left ventricle, in Chapter~\ref{hyper}, and 
Sections~\ref{subsec:sensitivity} and \ref{subsec:implem-hyper}.

As illustrated in Figure~\ref{fig:calibration:workflow}, the parameters of 
$(\mathscr{C}_\text{NC})$ are calibrated, using the time-independent 
outputs of $(\mathscr{C}_\text{C})$, and then included
in the 3D--0D model. 
It is noteworthy that the parameters of the 3D--0D model related to the 
electromechanical component of the left ventricle are never modified: 
as they 
are very delicate quantities and closely linked to the specific model itself, it would 
be unfeasible to try to adapt them to simulate hypertension. For this reason, the same 
parameters defined in \cite{3d-0d_ventr} are employed.
Finally, a simulation is 
performed under this comprehensive framework, allowing for the analysis of the model's 
behaviour under hypertensive conditions. The results obtained from 
$(\mathscr{C}_\text{C})$ and from the 3D--0D model are then compared to check if the 
calibration is performed successfully: by effectively aligning the outputs of 
$(\mathscr{C}_\text{C})$ and the 3D--0D model, the latter is now able to give real 
insight into the complex dynamics of hypertension, achieving this goal with minimal 
computational effort.

This calibration workflow also lays the foundation for cardiovascular DTs. By aligning 
the outputs of the 0D and 3D--0D models, it becomes possible to generate patient-specific 
representations of cardiac dynamics under hypertensive conditions. Such personalized 
models could be further exploited for predictive simulations, risk assessment and 
clinical monitoring, highlighting the potential of DTs in supporting 
individualized diagnosis and therapy planning.

\section{Numerical Results and Discussion}
\label{results}
Each simulation generates a total of 25 time-dependent variables 
and 46 outputs, if carried out with $(\mathscr{C}_\text{C})$, or 20 time-dependent 
variables and 45 outputs, if the 3D--0D model is employed (Section~\ref{subsec:lumped} 
for further details); in addition, the 3D--0D model 
calculates key quantities from the electromechanical model 
(Section~\ref{subsec:3d-0d-model}), of the left ventricle mesh.
For both 0D and 3D--0D models, right and left ventricular stroke work (respectively, $\text{RV}_\text{SW}$ and
$\text{LV}_\text{SW}$) is computed as the area enclosed within the corresponding PV loop.
Stroke work was not used for parameter calibration because it showed low sensitivity to 
model parameters, making it less informative for optimization purposes. Moreover, the 
information it provides is largely redundant with that already captured by the pressure 
and volume curves, which are directly included in the loss function.
Although these quantities are not directly used for parameter calibration, they provide 
valuable insight into ventricular function.

Each condition is simulated through parameter adjustment, as outlined in 
Section~\ref{subsec:implem-hyper}; the timestep employed in 0D simulation is $10^{-3}\,\text{s}$, 
whereas a more restrictive value of $10^{-4}\,\text{s}$ is preferred for the 3D--0D 
model to account for the fast time scale charachtering electrophysiology. 
The calibration process for $(\mathscr{C}_\text{C})$ 
is performed as described in Section~\ref{subsubsec:calibration:0d}, 
while the calibration for the 3D--0D model is carried out as outlined in 
Section~\ref{subsubsec:calibration:3d}. 
As mentioned in Section~\ref{subsubsec:calibration:3d}, in order to ensure consistency in 
the calibration process, the left ventricle parameters remain unchanged in the 3D--0D 
model, since the left ventricle mesh is constructed on the basis of a healthy 
individual: from a clinical point of view, this fact may be interpreted as a 
representation of an individual with a hypertensive 
condition, which might also be very severe, but not of long-standing duration, i.e. 
the high pressure has not yet seriously compromised the left ventricular 
structure.
It is worth noting that the results presented in the following sections also illustrate 
the potential of the 
calibrated models to serve as the foundation for patient-specific cardiovascular DTs, 
enabling personalized insights and predictive analyses.

\subsection{Systemic Hypertension}
\label{subsec:sys-hyper}%

\begin{table}[t]
    \centering 
    \begin{tabular}{l l l l l l l l l l}
    \hline
    \multicolumn{2}{c}{} & \multicolumn{4}{c}{0D model} & \multicolumn{4}{c}{3D--0D model}\\
    \cmidrule(lr){3-6} \cmidrule(lr){7-10} 
    Output & Unit & Healthy & Mild & Moderate & Severe & Healthy & Mild & Moderate & Severe \\
    \hline
    $\text{LV}_\text{SW}$    & $\text{mmHg}\cdot \text{L}$   & 5.9 & 6.3 &   6.9   &  7.5   & 6.3 & 6.5 &   7.1   &  7.2 \\ 
    CI                         &$\text{m}^2\cdot \text{L}/\text{min}$ & 2.9 & 2.9     &3.0      &3.0   & 3.0  &  3.0  &  3.0    &  3.1      \\
    $\text{SAP}_\text{max}$    & mmHg   & 109.6 &  117.0  &   128.2  &  144.4  & 106.9 &  111.4  &   120.9  &  127.4\\
    $\text{SAP}_\text{min}$    & mmHg   & 71.2  & 74.5  &   80.2   &  83.8    & 72.8  & 76.5  &   81.5   &  80.8 \\
        \hline
    \end{tabular}
    \caption{Variations of selected time-independent 
    outputs in systemic hypertension, from both $(\mathscr{C}_\text{C})$ and 3D--0D model.
    }
    \label{table:sys-out}
\end{table} 

\begin{figure}[p]
    \centering
    \subfloat[LV PV loop.\label{fig:sys0d:lv-loop}]{
        \includegraphics[height=0.25\linewidth]{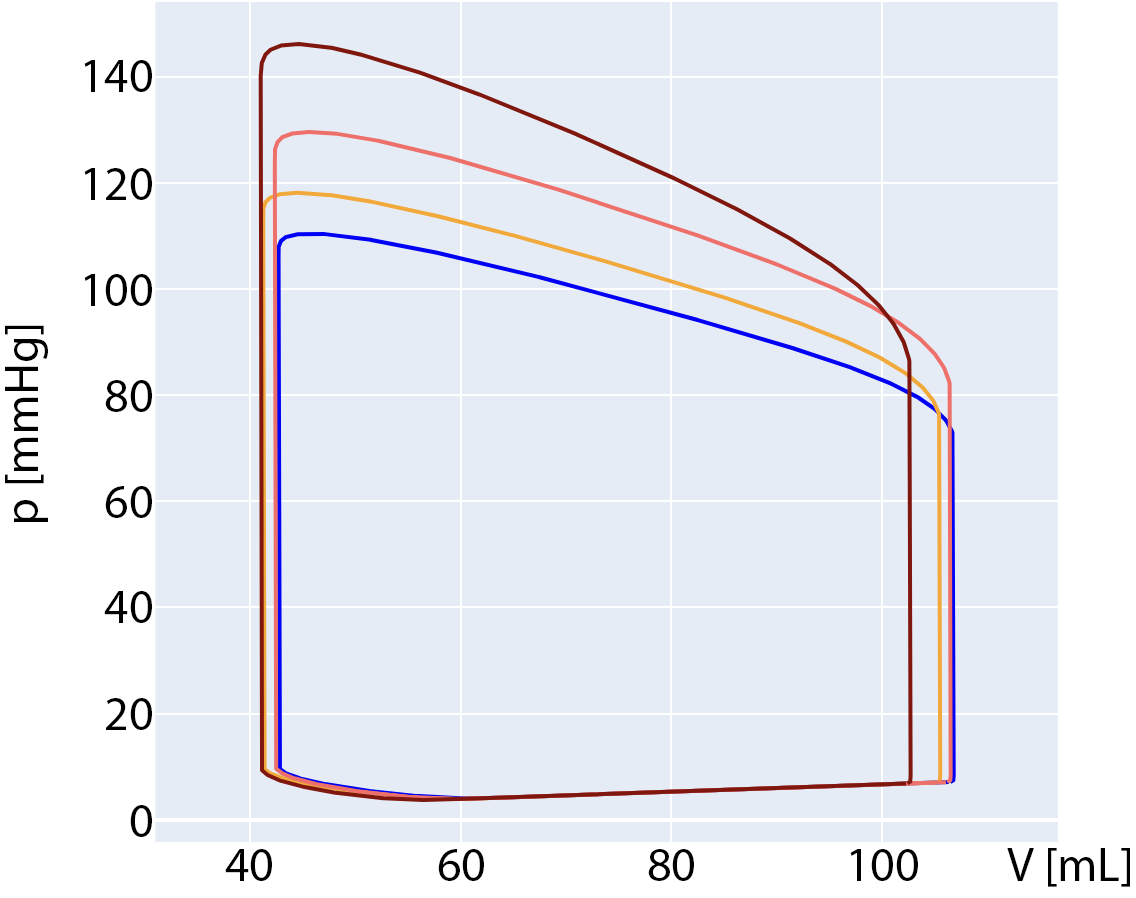}
    }
    \subfloat[LV PV loop.\label{fig:sys3d:lv-loop}]{
        \includegraphics[height=0.25\linewidth]{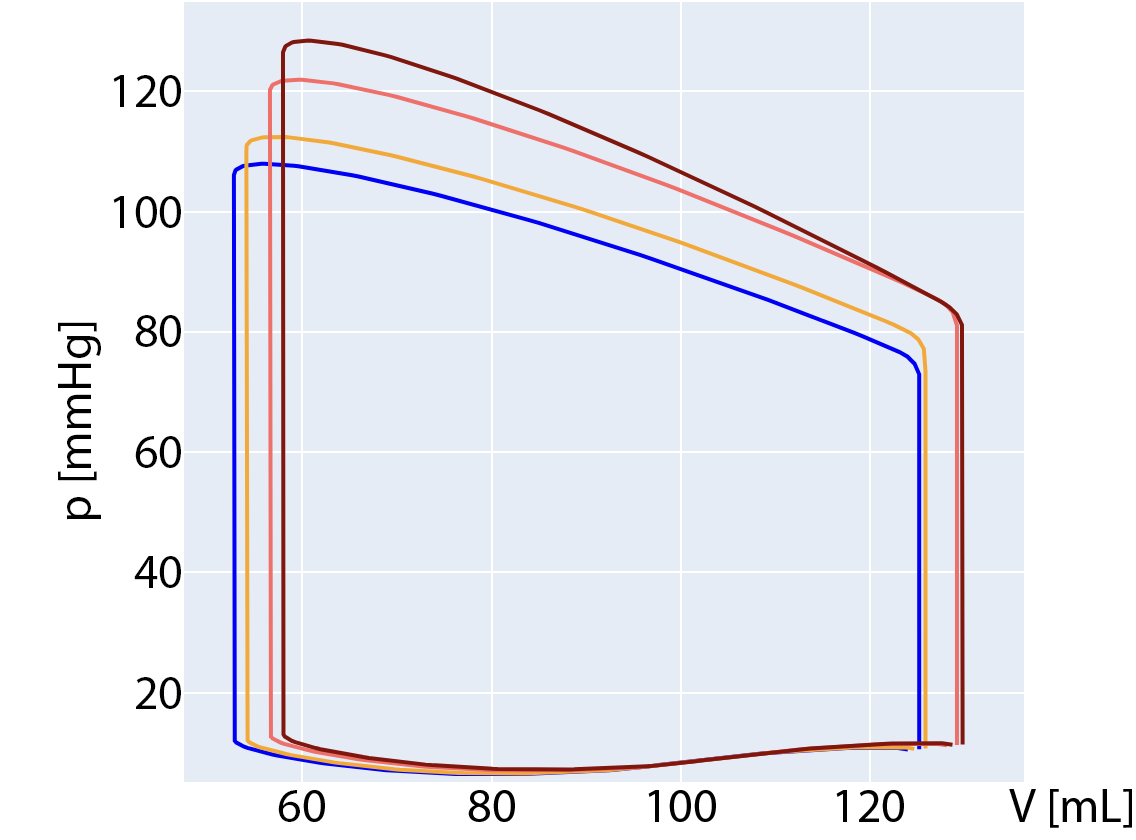}
    }\\
    \subfloat[$p_\text{LV}$.\label{fig:sys0d:plv}]{
        \includegraphics[width=0.32\linewidth]{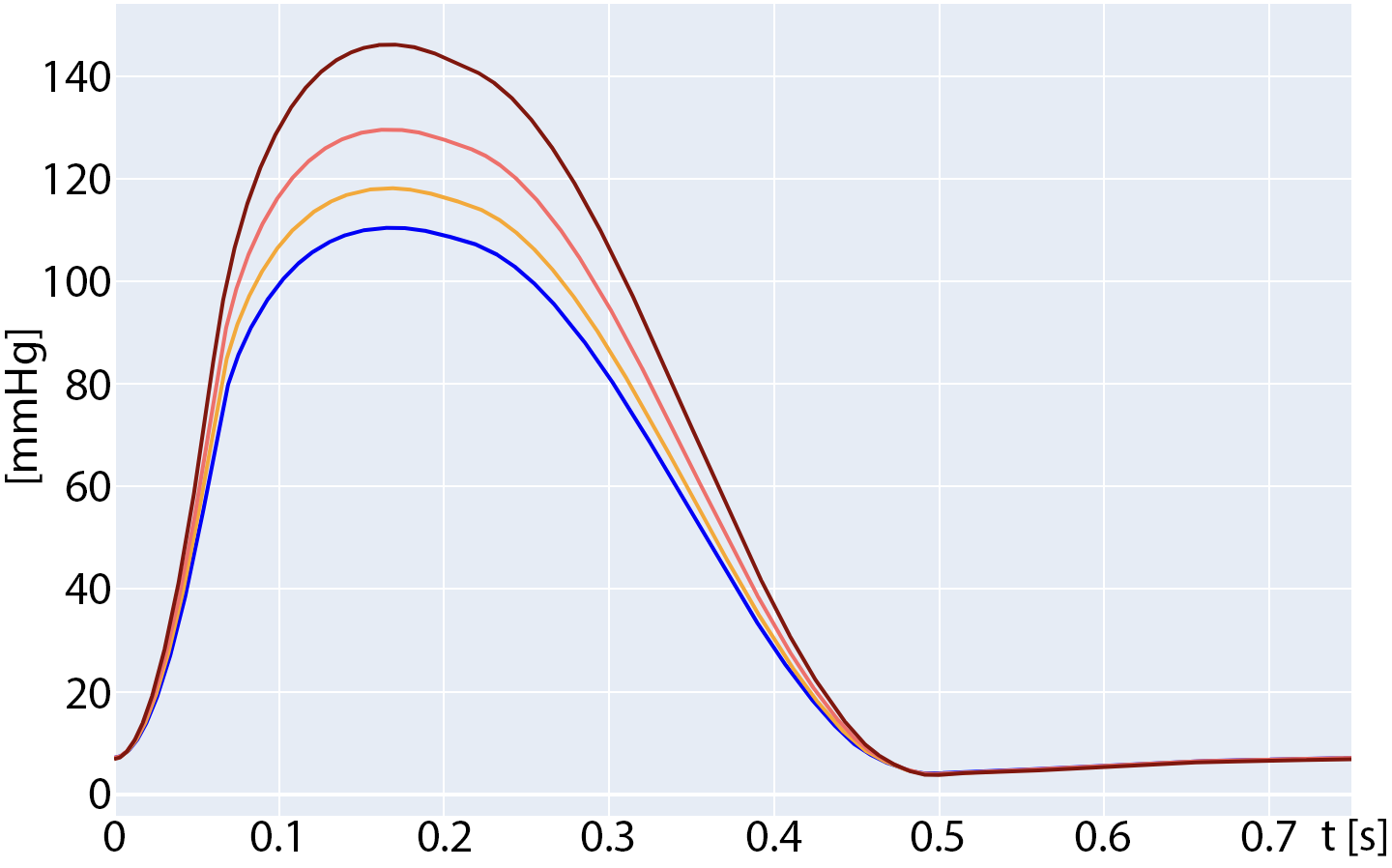}
    }
    \subfloat[$p_\text{LV}$.\label{fig:sys3d:plv}]{
        \includegraphics[width=0.32\linewidth]{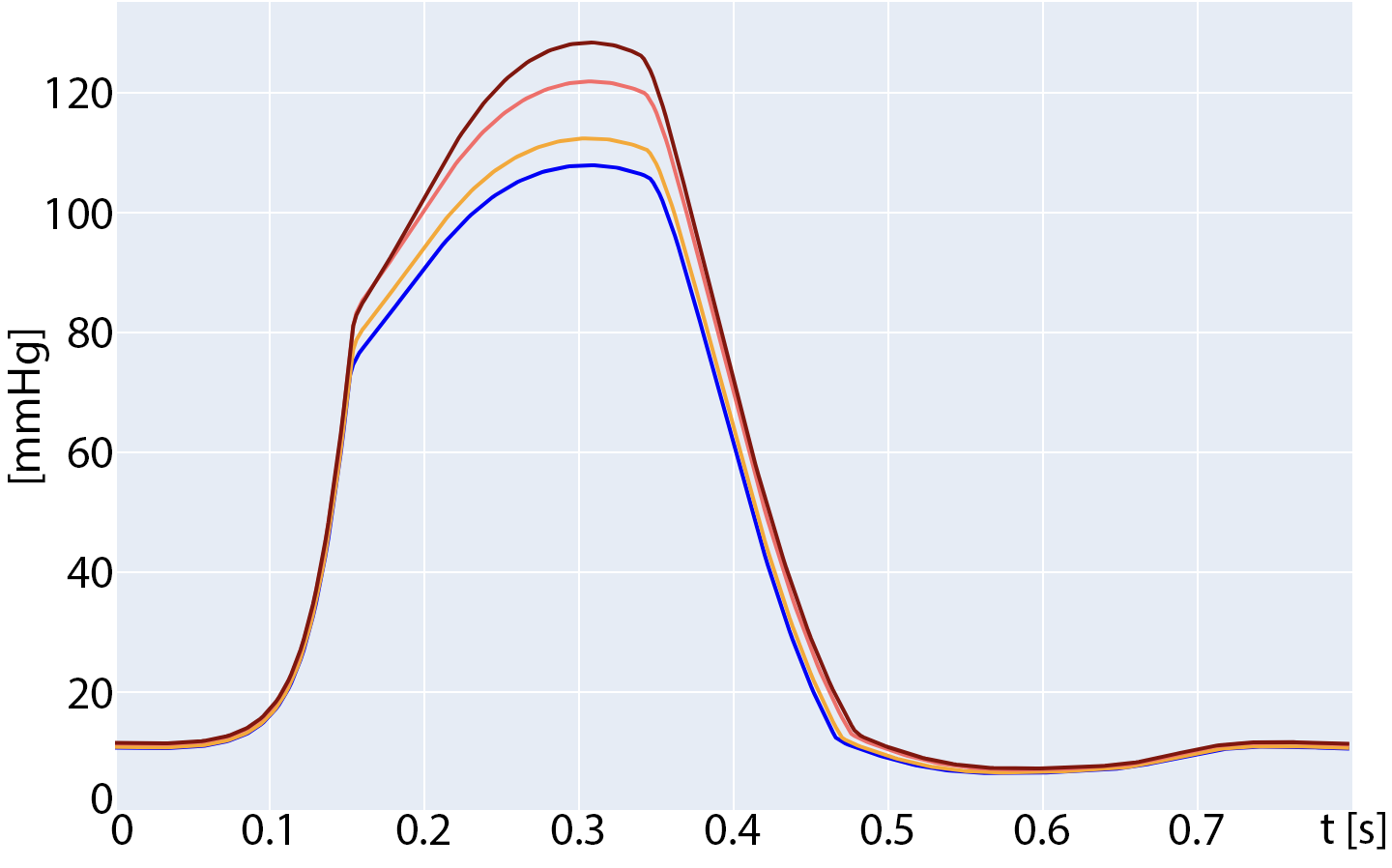}
    }\\\hspace{0.112\linewidth}
    \subfloat[$Q_\text{AV}$.\label{fig:sys0d:qav}]{
        \includegraphics[width=0.32\linewidth]{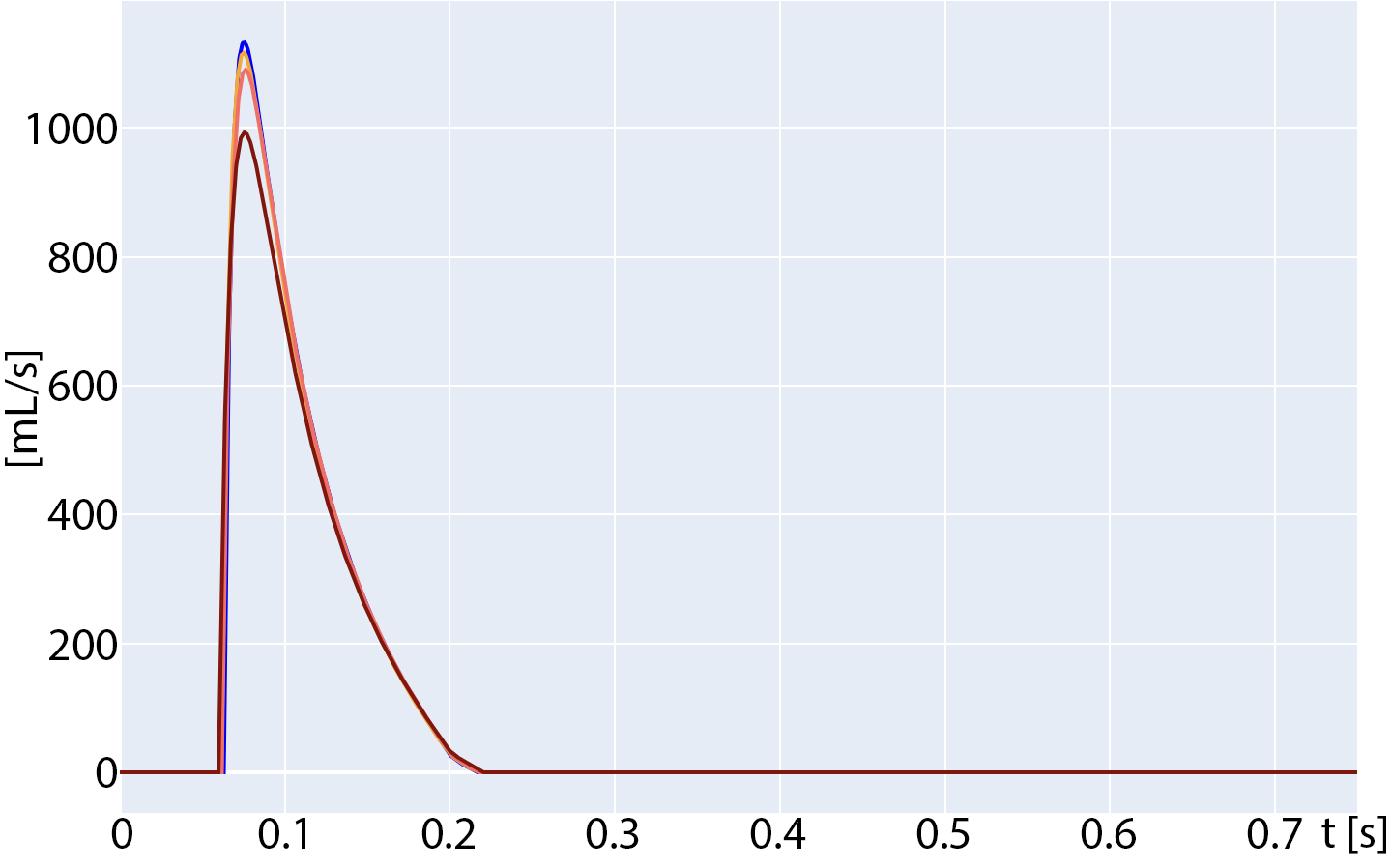}
    }
    \subfloat[$Q_\text{AV}$.\label{fig:sys3d:qav}]{
        \includegraphics[width=0.32\linewidth]{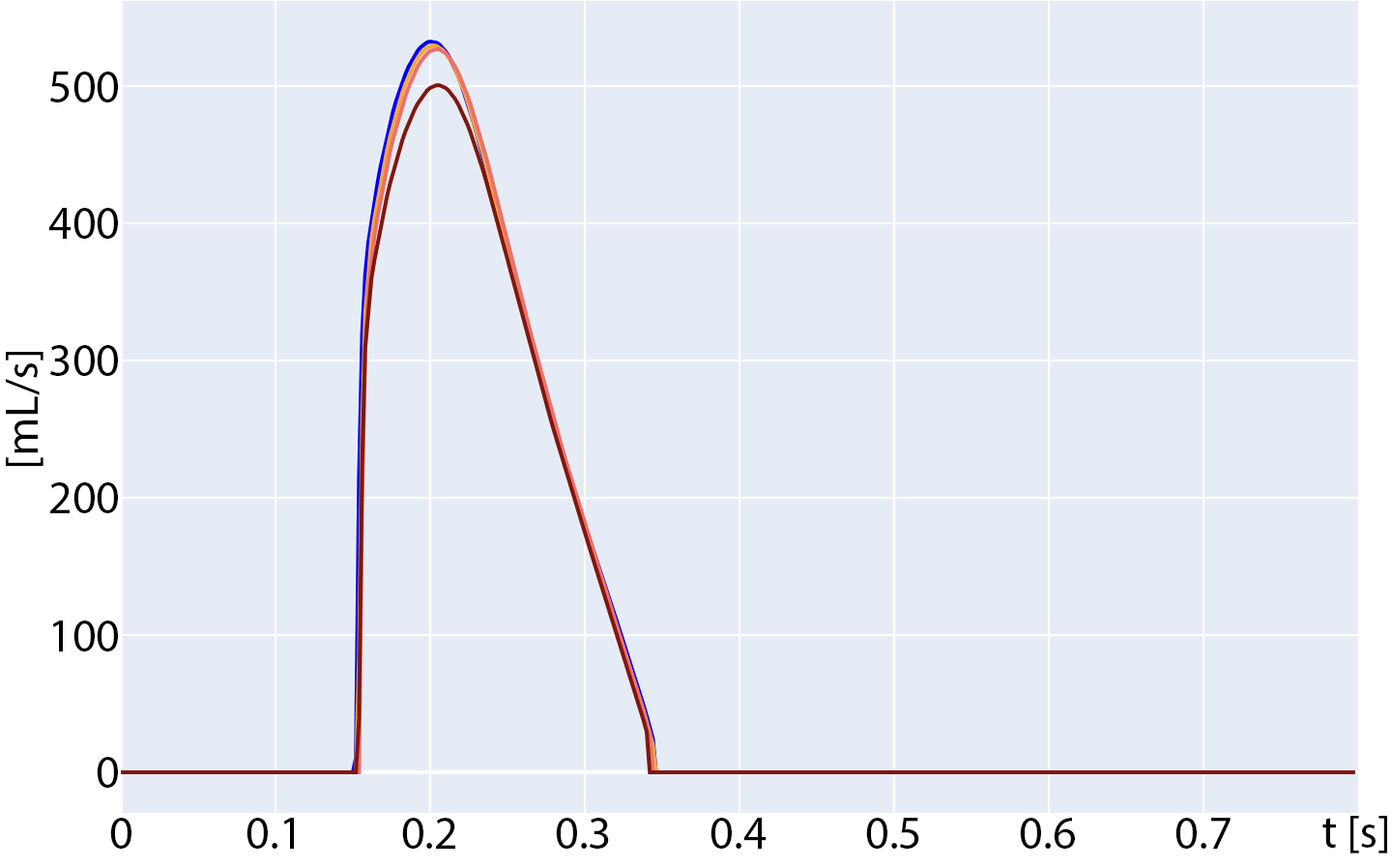}
    }
    \begin{subfigure}[b]{0.12\textwidth}
        \raisebox{0.6\height}{\includegraphics[width=\textwidth]{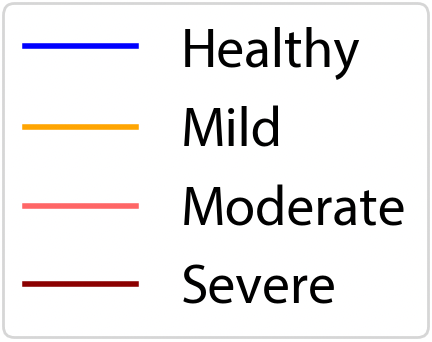}}
    \end{subfigure}\\
    \subfloat[$p_\text{AR}^\text{SYS}$.\label{fig:sys0d:parsys}]{
        \includegraphics[width=0.32\linewidth]{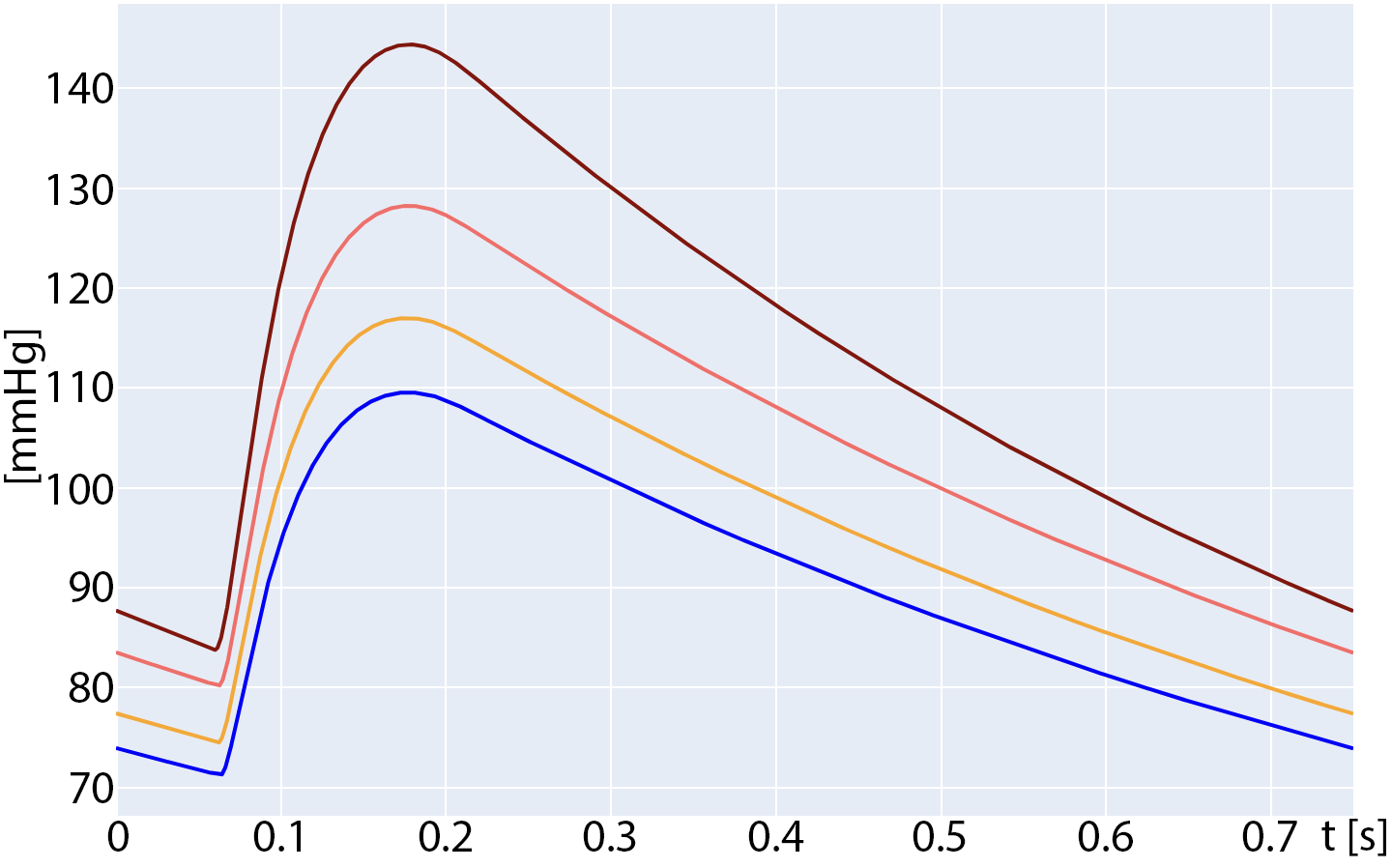}
    }
    \subfloat[$p_\text{AR}^\text{SYS}$.\label{fig:sys3d:parsys}]{
        \includegraphics[width=0.32\linewidth]{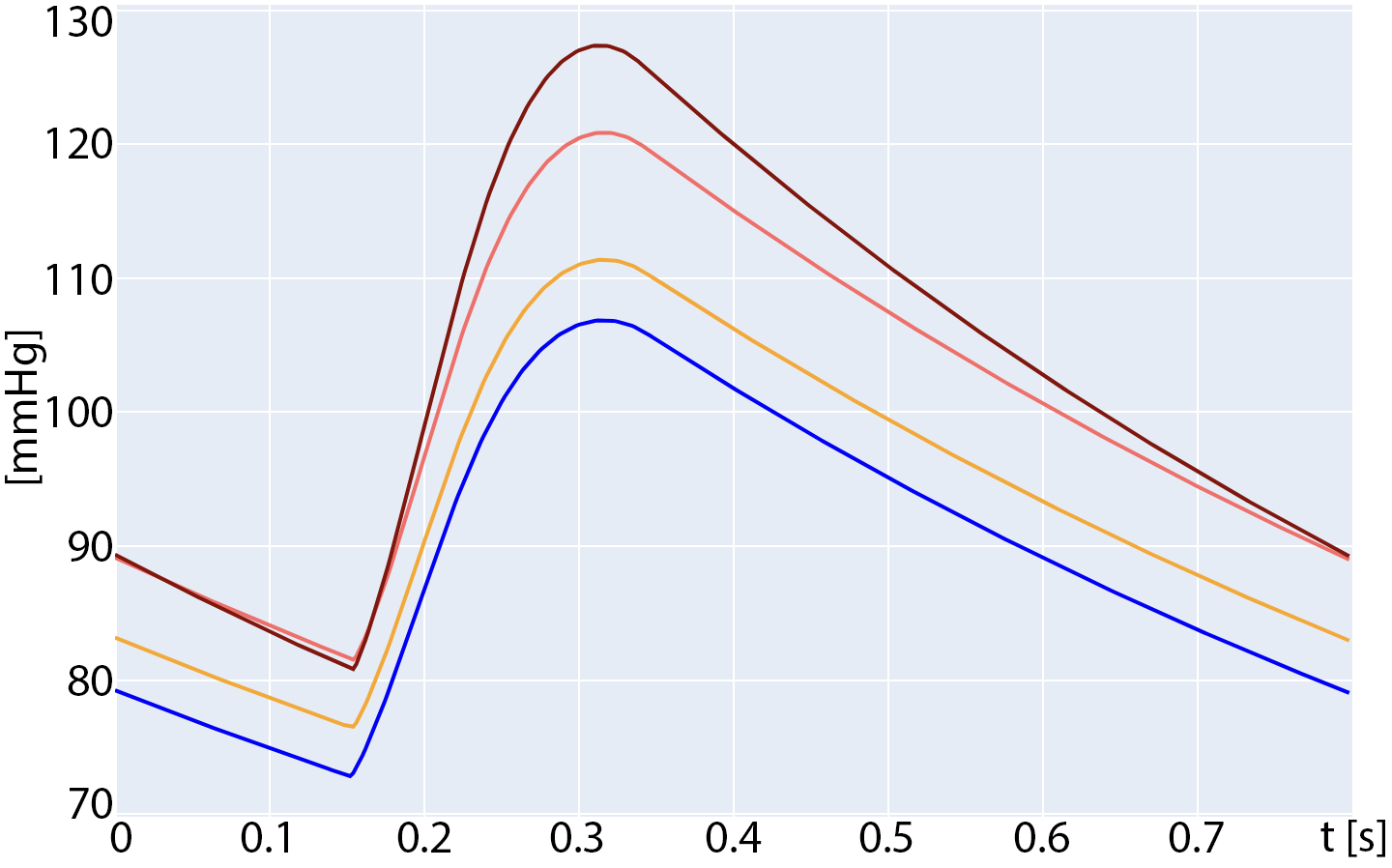}
    }\\
    \subfloat[$Q_\text{AR}^\text{SYS}$.\label{fig:sys0d:qarsys}]{
        \includegraphics[width=0.32\linewidth]{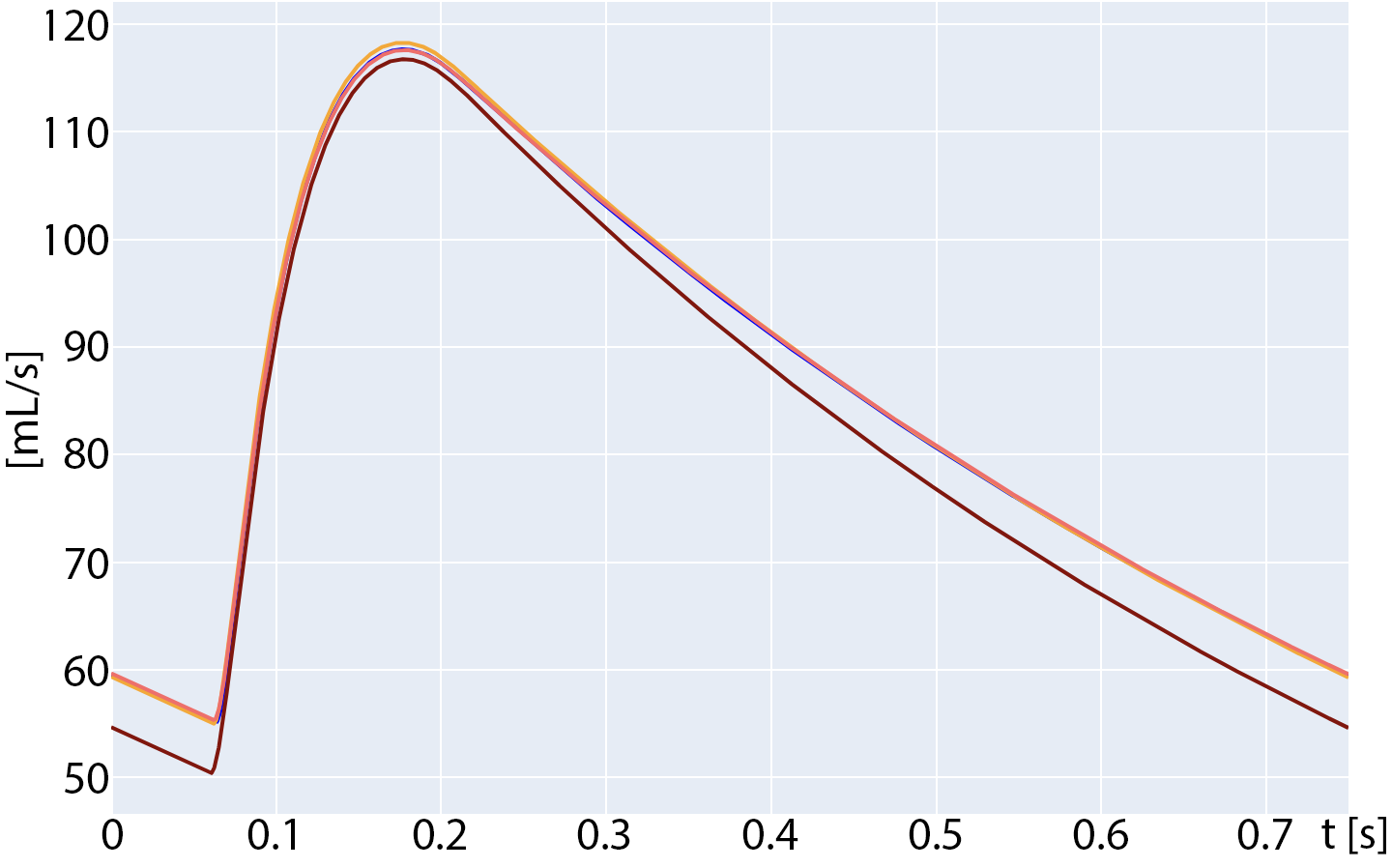}
    }
    \subfloat[$p_\text{AR}^\text{PUL}$.\label{fig:sys3d:parpul}]{
        \includegraphics[width=0.32\linewidth]{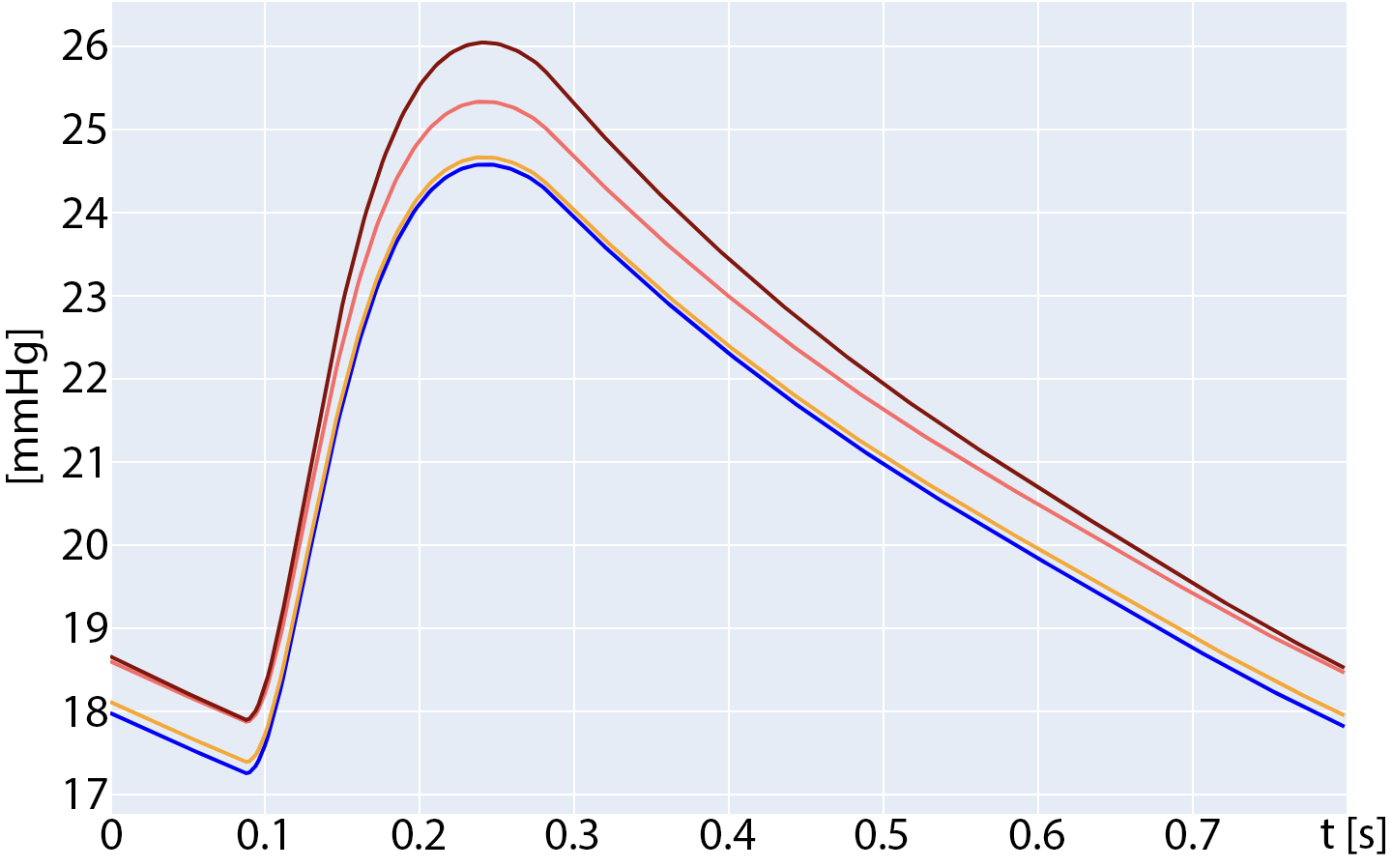}
    }
    \caption{Time-dependent variables from $(\mathscr{C}_\text{C})$ (on the left) and
    from the 3D--0D model (on the right).
    Mild, moderate and severe systemic hypertension 
    (in orange, light red and dark red) is compared with a 
    healthy individual (in blue).
    }
    \label{fig:sys-var}
\end{figure}

\begin{figure}[t!]
    \centering
    \subfloat{
        \rotatebox{90}{\qquad $t=0.26\text{s}$}
        \includegraphics[width=0.225\linewidth]{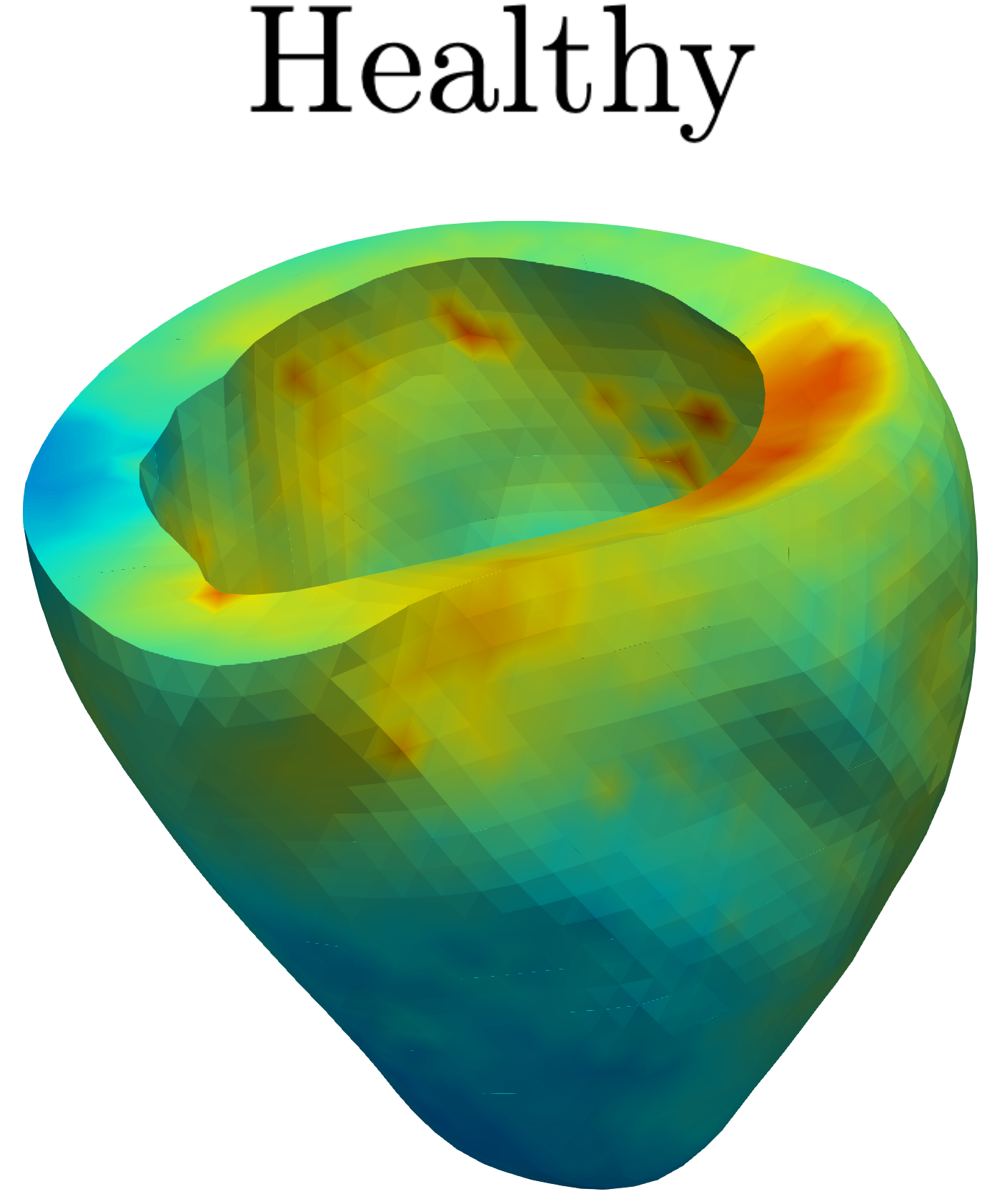}
    }
    \subfloat{
        \includegraphics[width=0.225\linewidth]{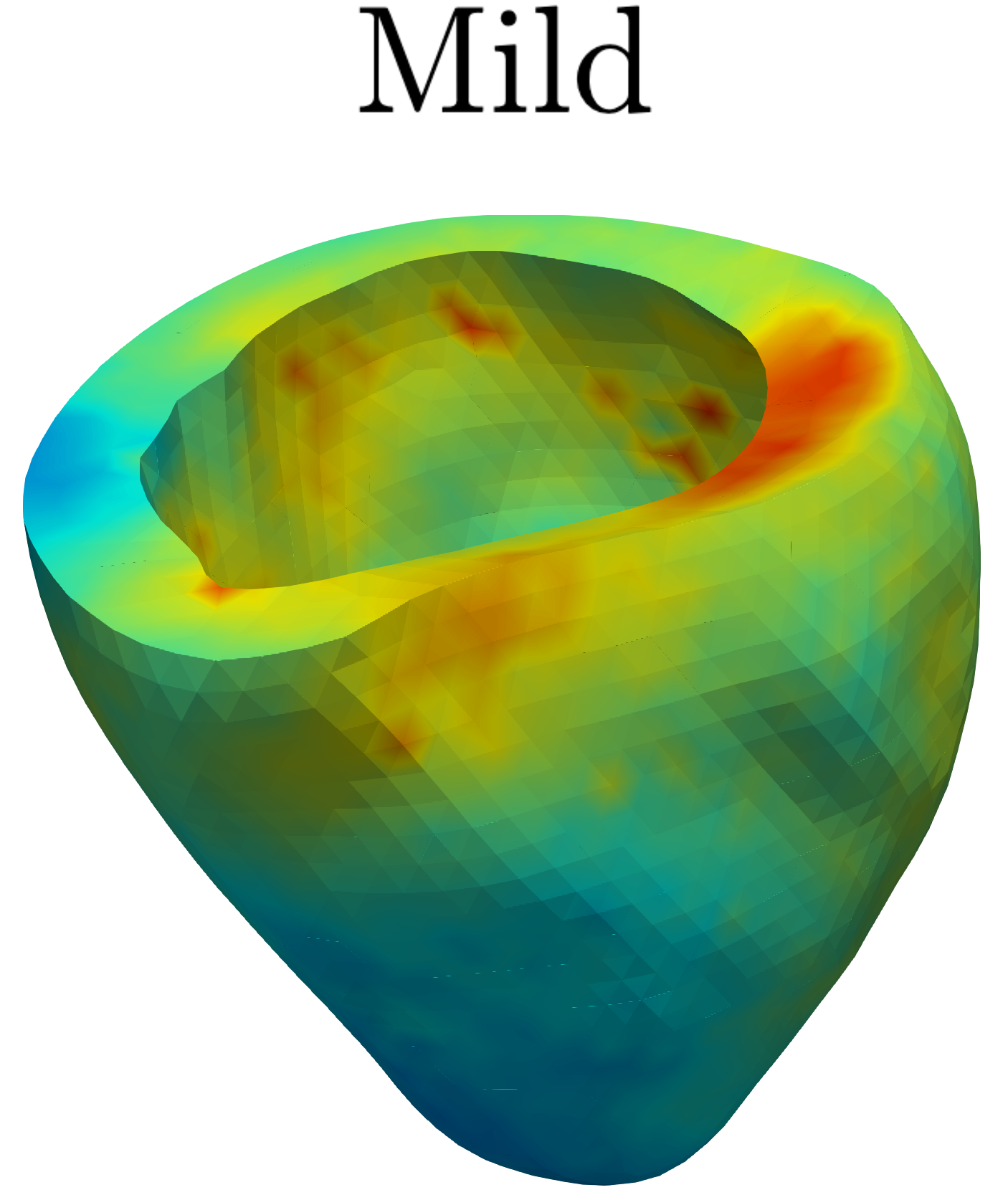}
    }
    \subfloat{
        \includegraphics[width=0.225\linewidth]{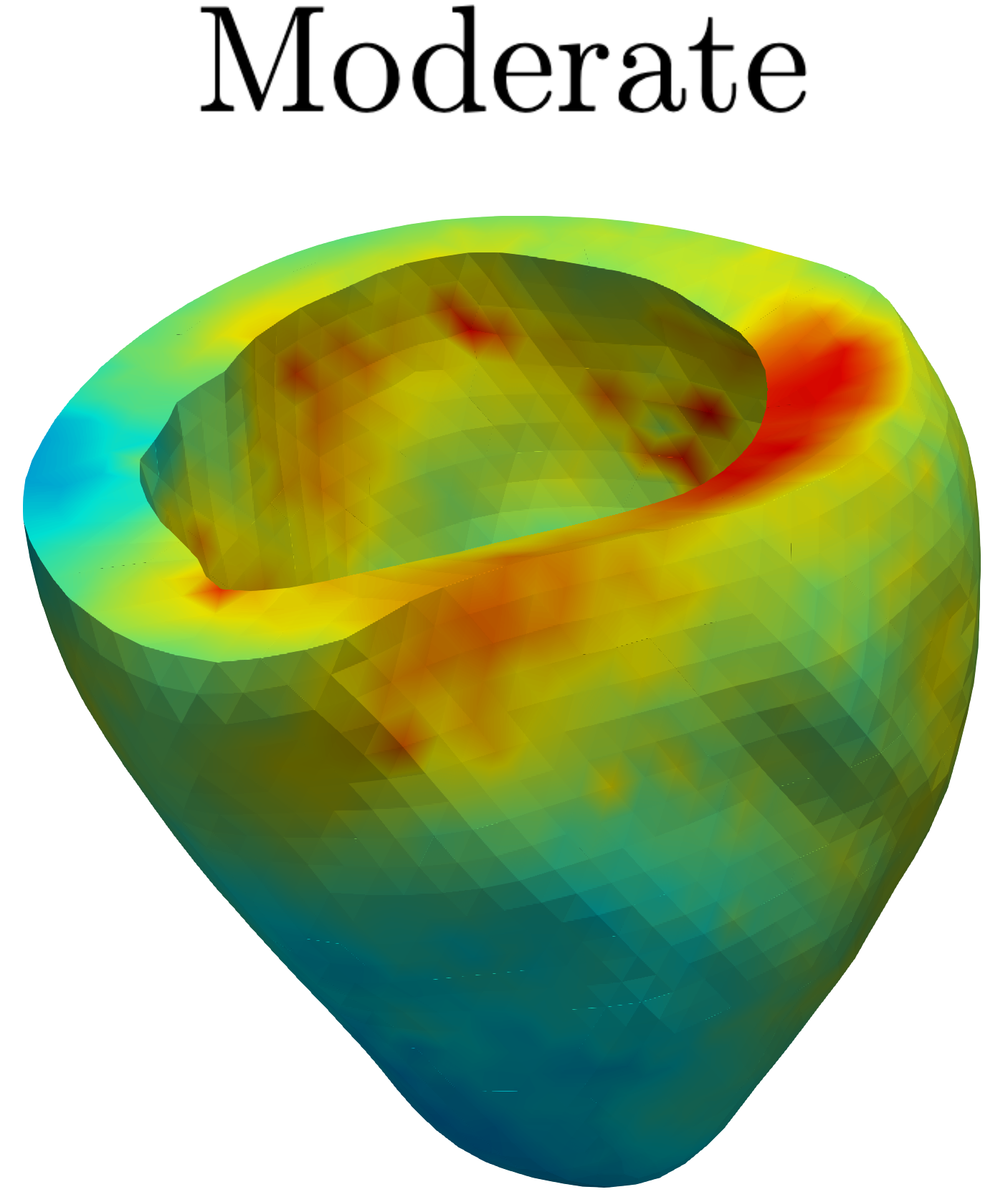}
    }
    \subfloat{
        \includegraphics[width=0.225\linewidth]{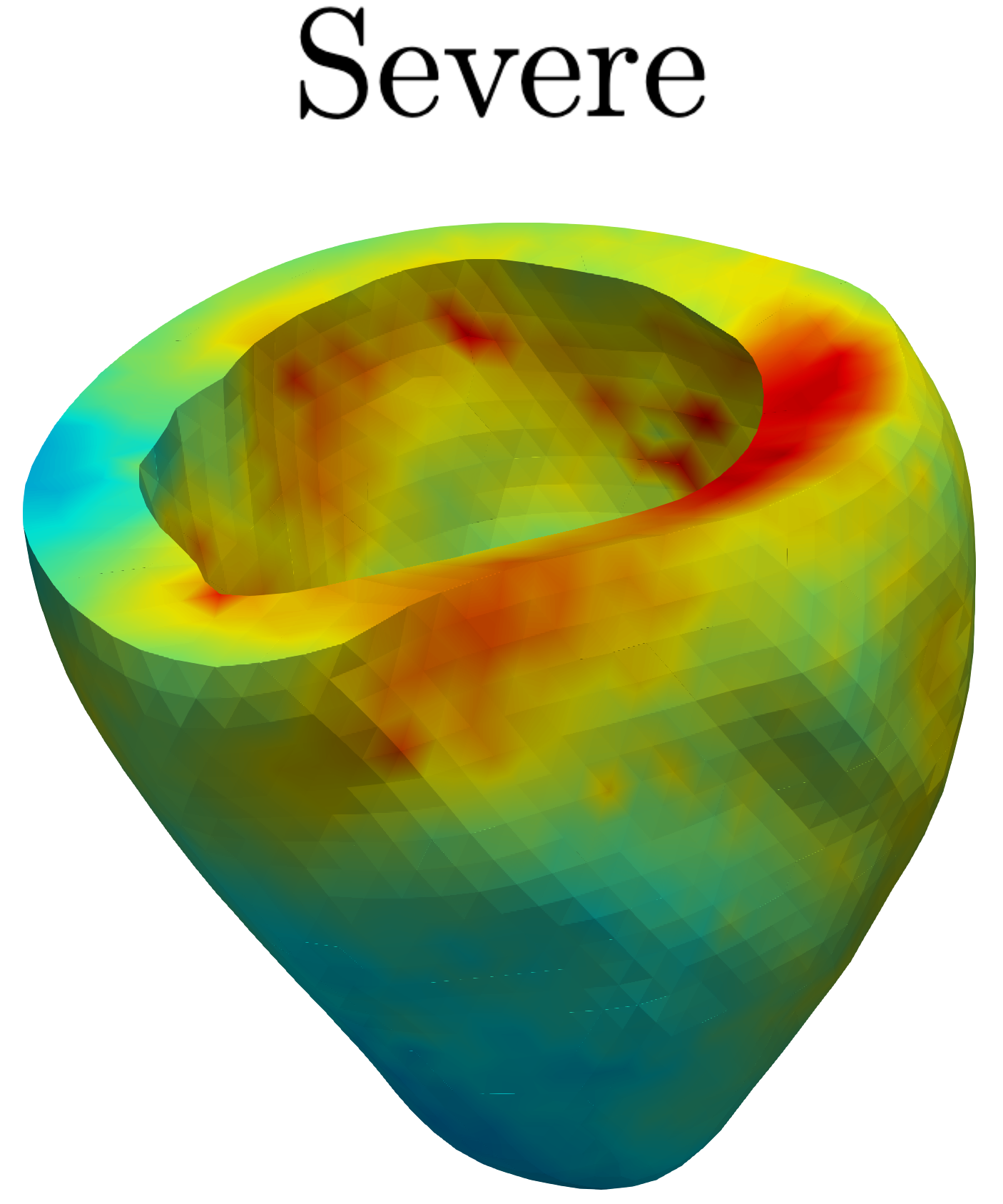}
    }\\
    \subfloat{
        \rotatebox{90}{\qquad $t=0.35\text{s}$}
        \includegraphics[width=0.225\linewidth]{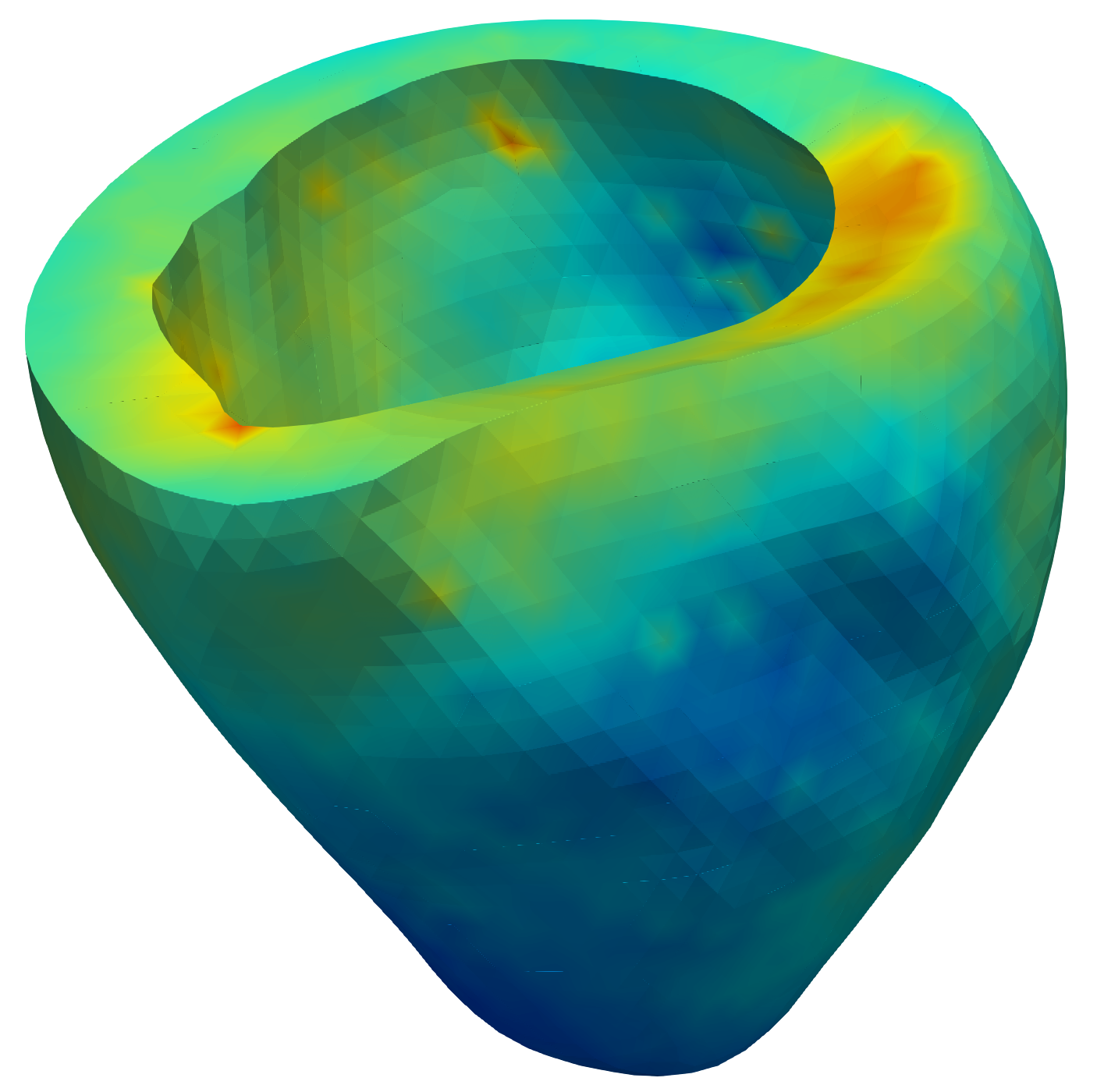}
    }
    \subfloat{
        \includegraphics[width=0.225\linewidth]{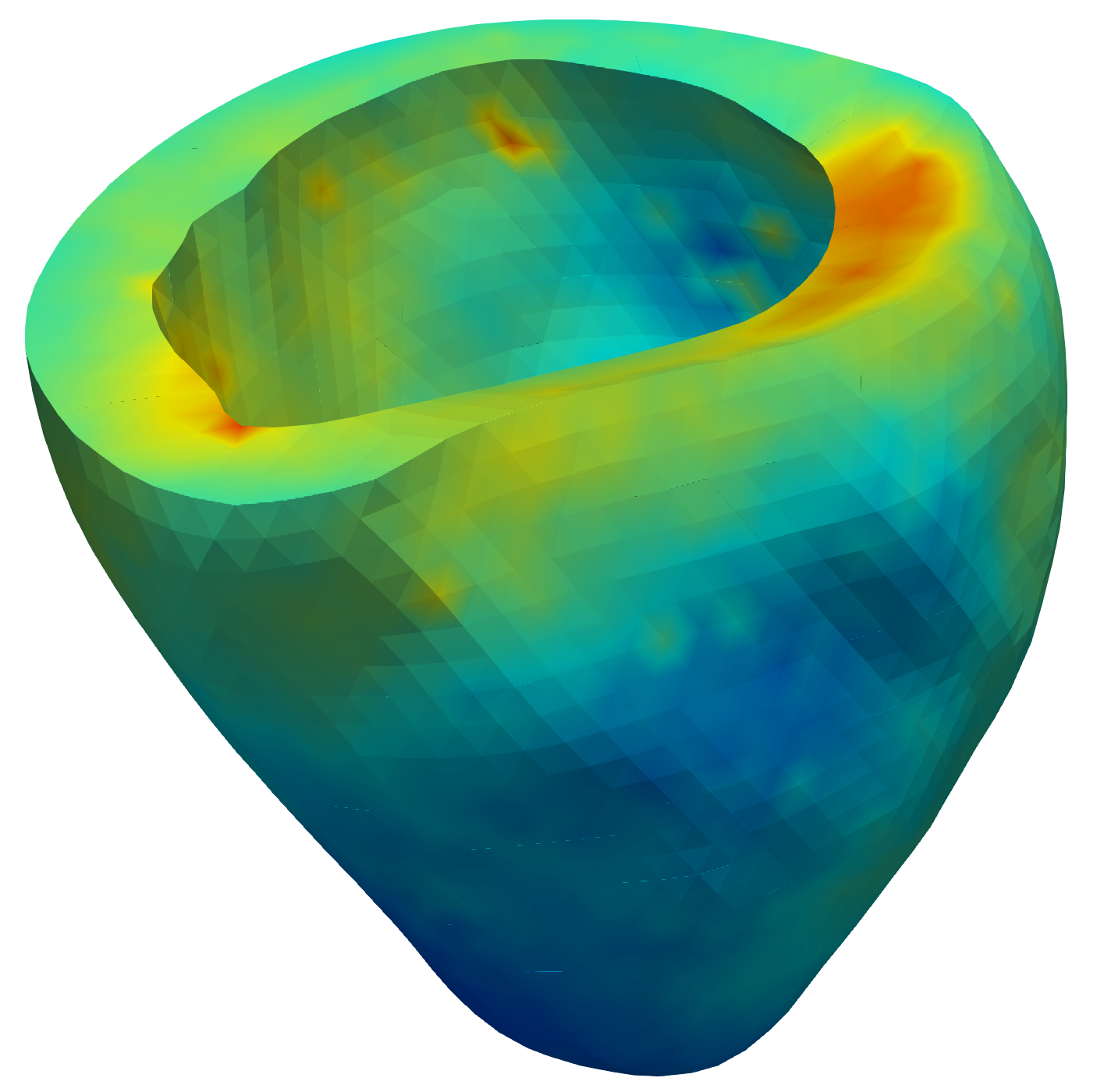}
    }
    \subfloat{
        \includegraphics[width=0.225\linewidth]{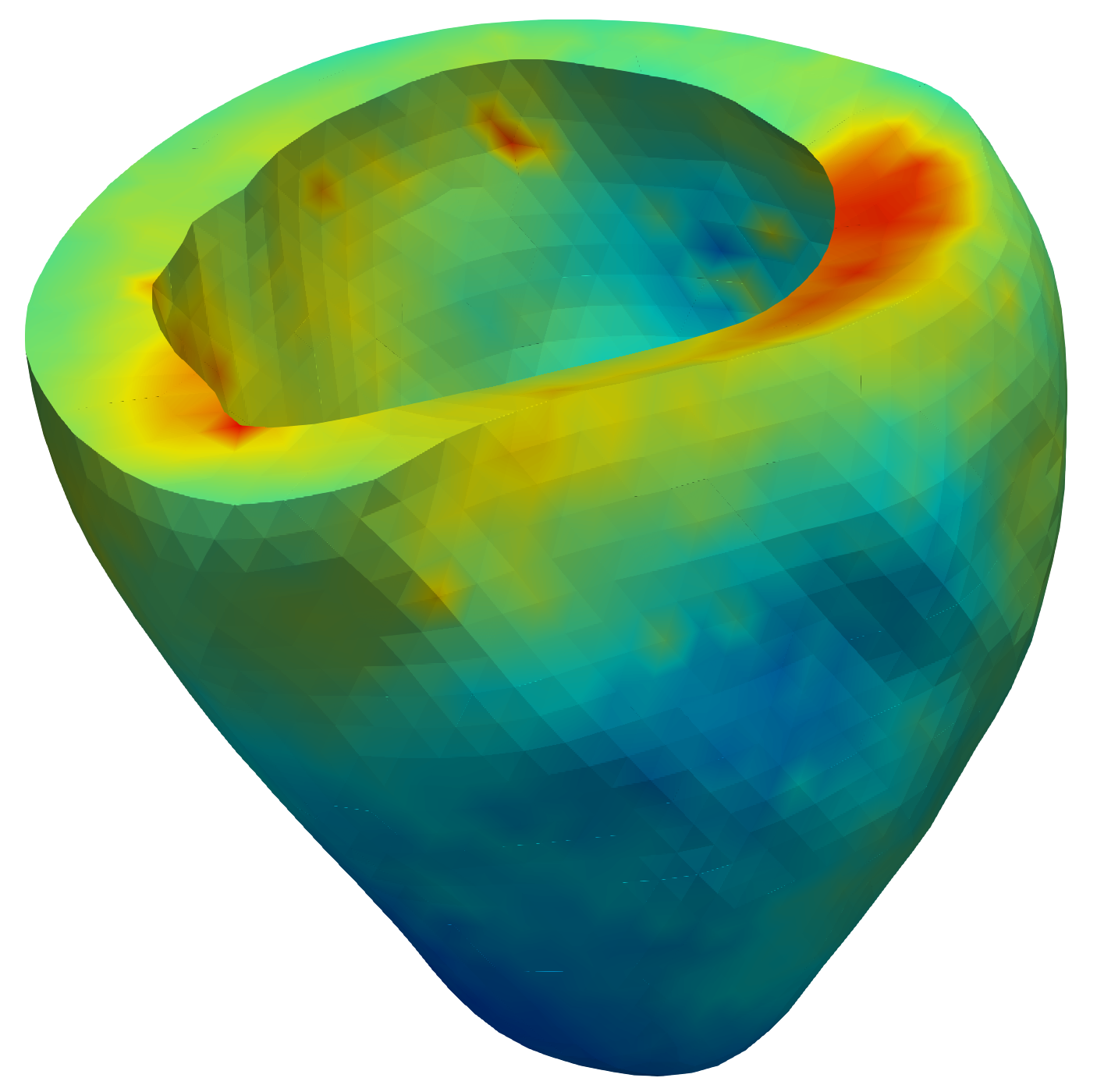}
    }
    \subfloat{
        \includegraphics[width=0.225\linewidth]{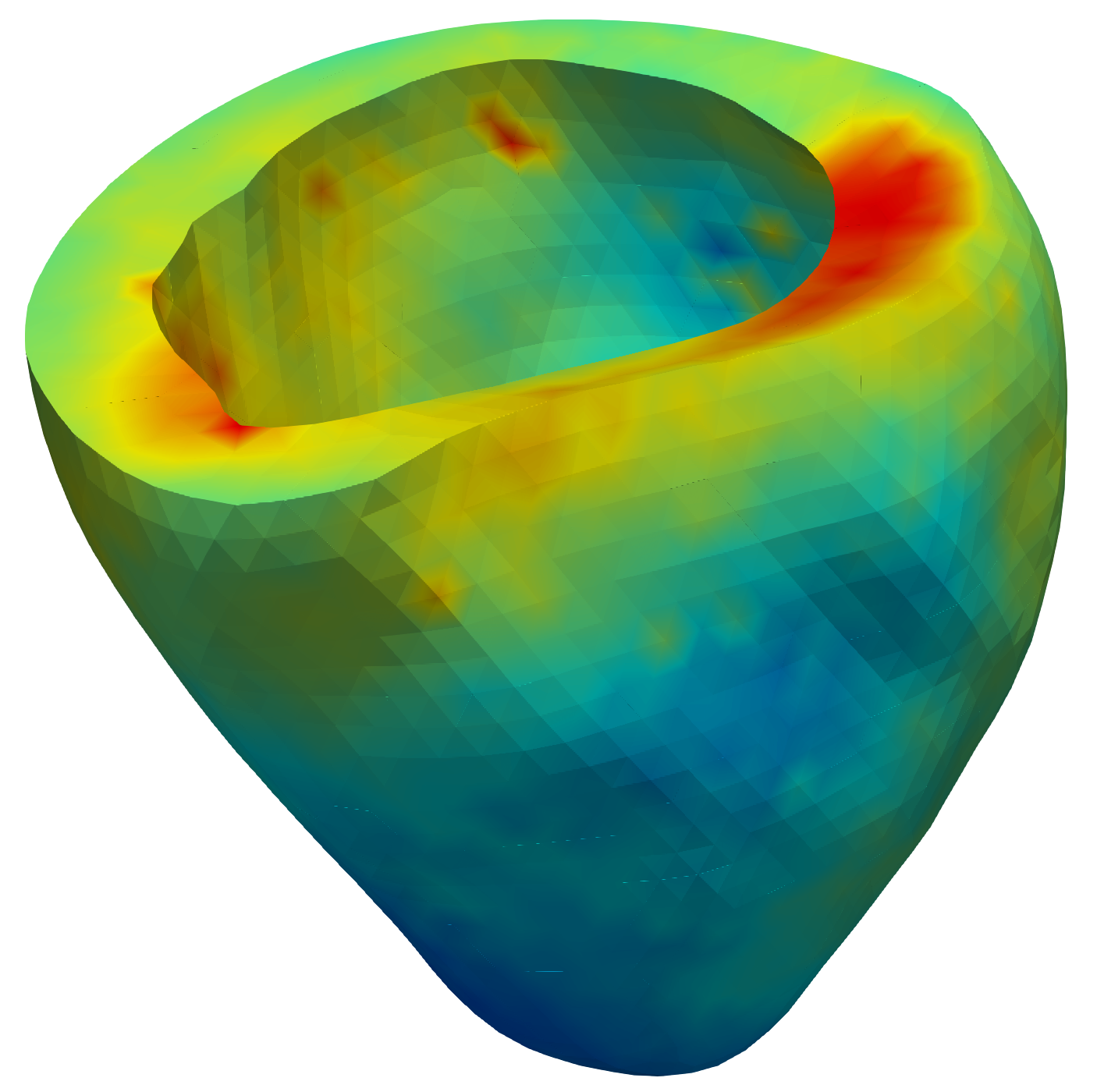}
    }\\
    \subfloat{
        \includegraphics[width=0.98\linewidth]{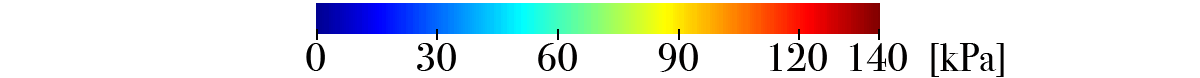}
    }
    \caption{Evolution 
    of the active tension $T_a$ in the left ventricle 
    during the last cardiac cycle.
    Each picture is displayed on the reference configuration $\Omega_0$.
    Mild, moderate and severe systemic hypertension is compared with a 
    healthy individual.
    }
    \label{fig:sys3d-ta}
\end{figure}

\begin{figure}[t!]
    \centering
    \subfloat{
        \rotatebox{90}{\qquad $t=0.15\text{s}$}
        \includegraphics[width=0.225\linewidth]{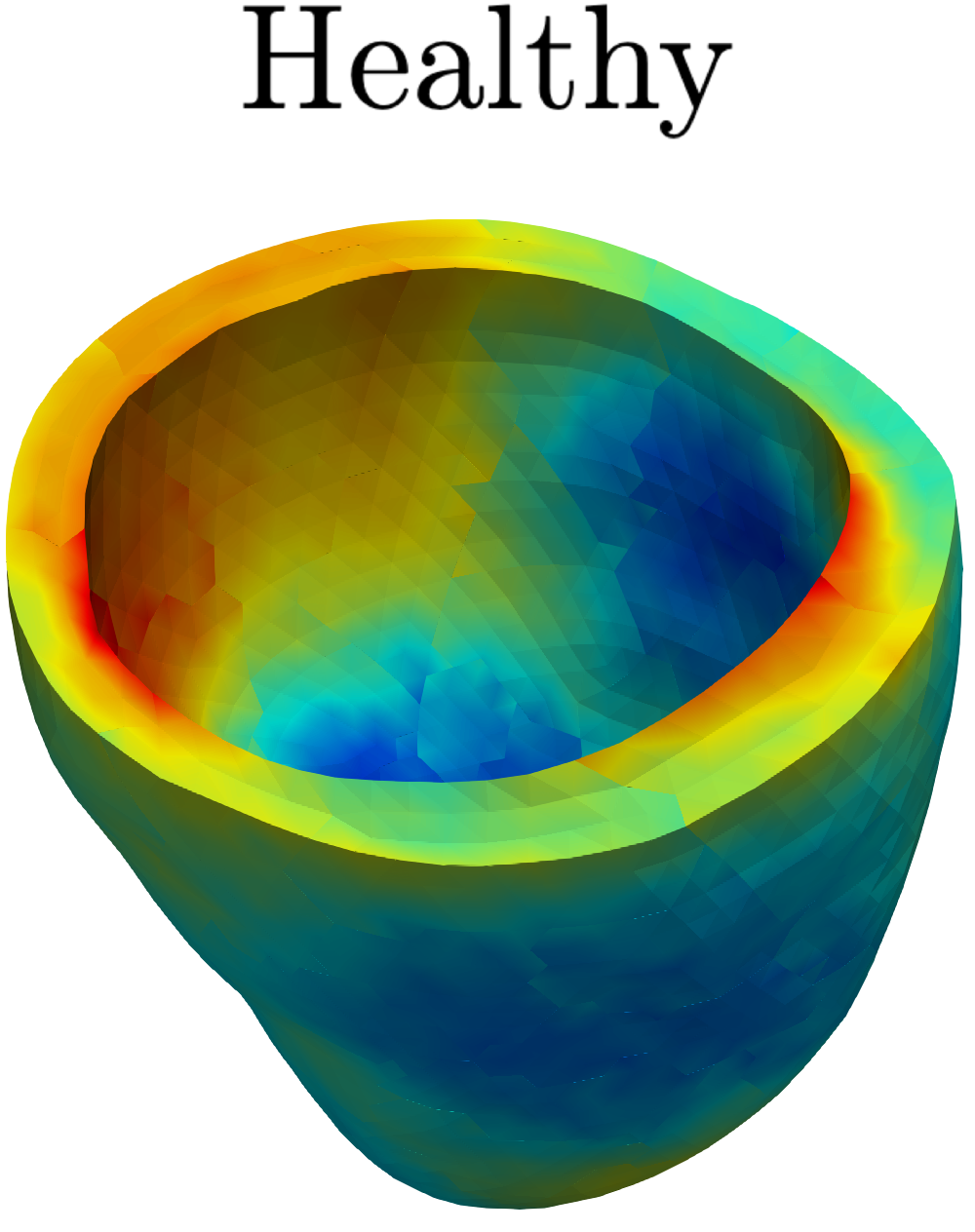}
    }
    \subfloat{
        \includegraphics[width=0.225\linewidth]{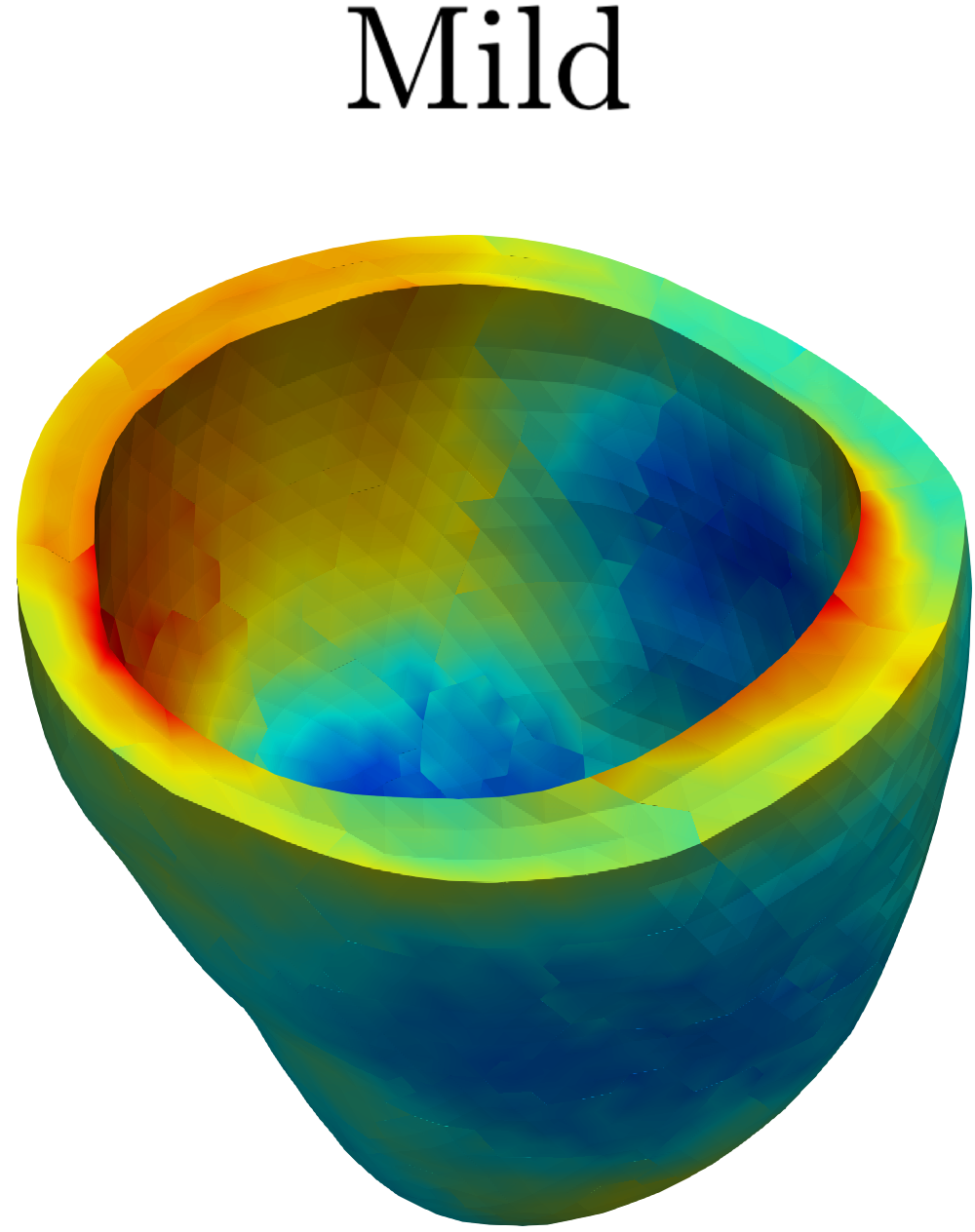}
    }
    \subfloat{
        \includegraphics[width=0.225\linewidth]{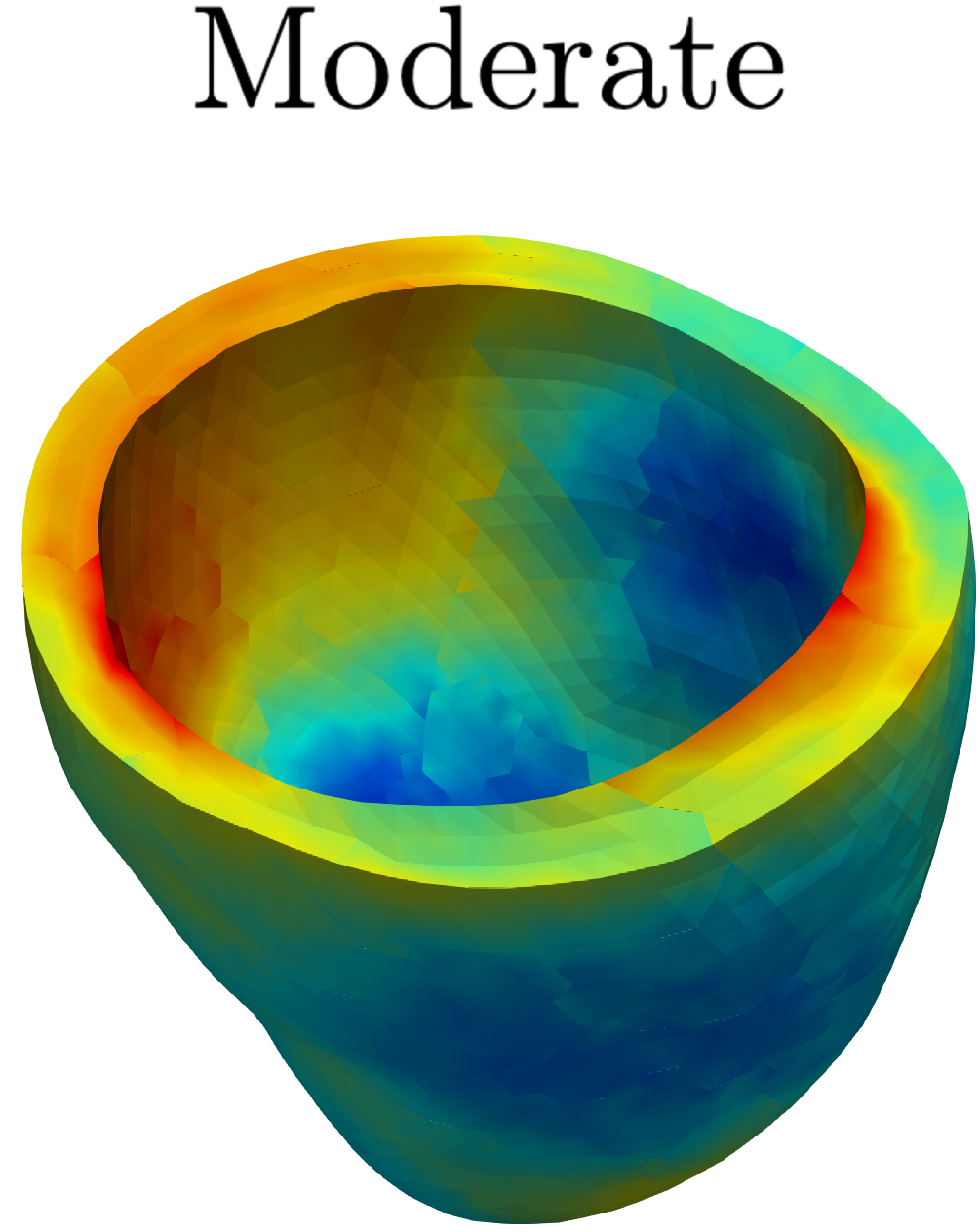}
    }
    \subfloat{
        \includegraphics[width=0.225\linewidth]{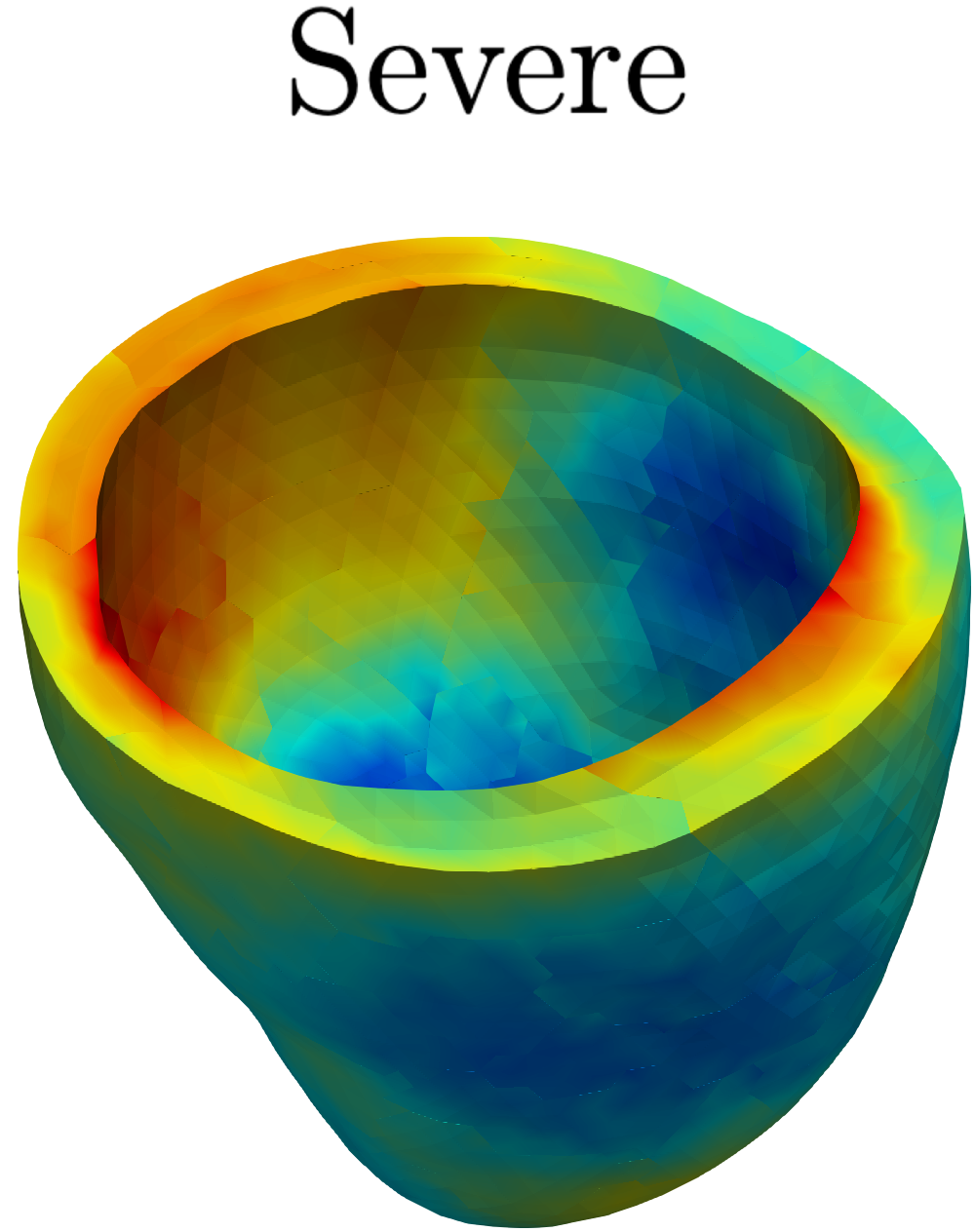}
    }\\
    \subfloat{
        \rotatebox{90}{\qquad $t=0.23\text{s}$}
        \includegraphics[width=0.225\linewidth]{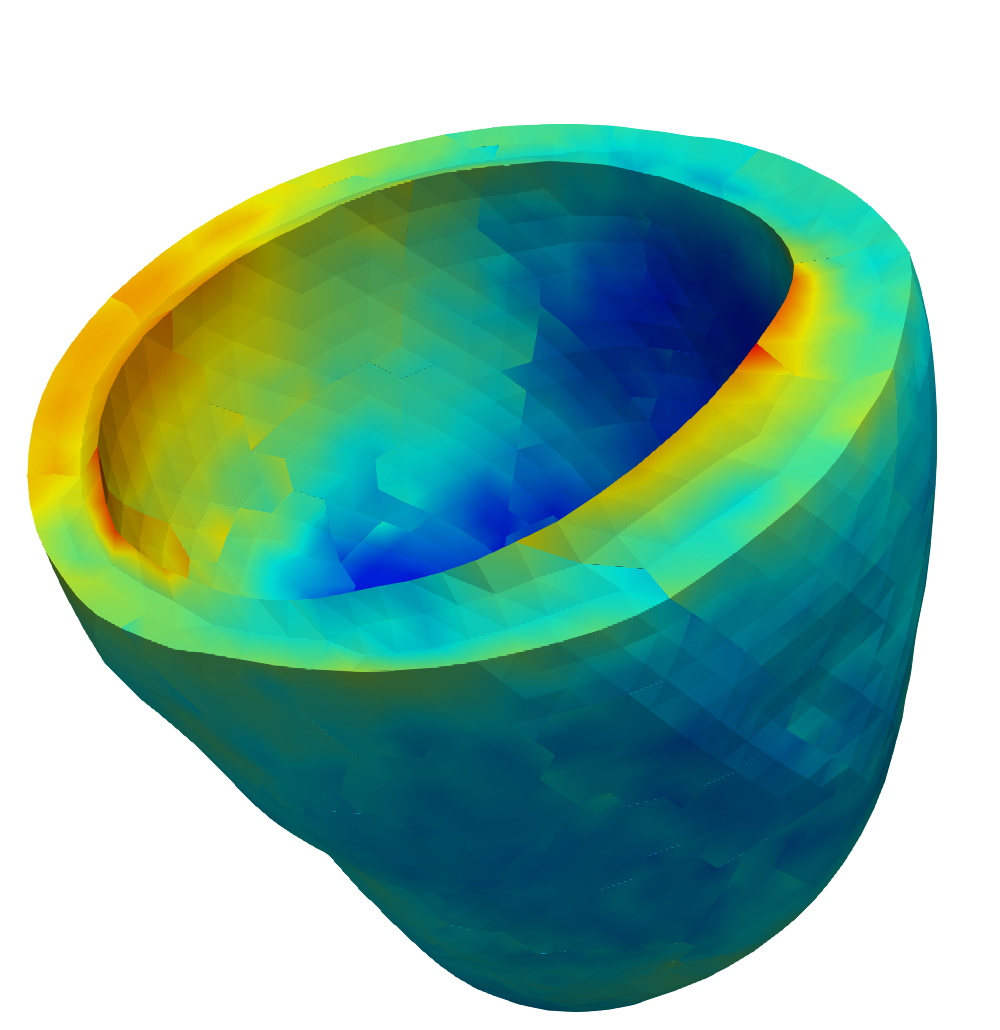}
    }
    \subfloat{
        \includegraphics[width=0.225\linewidth]{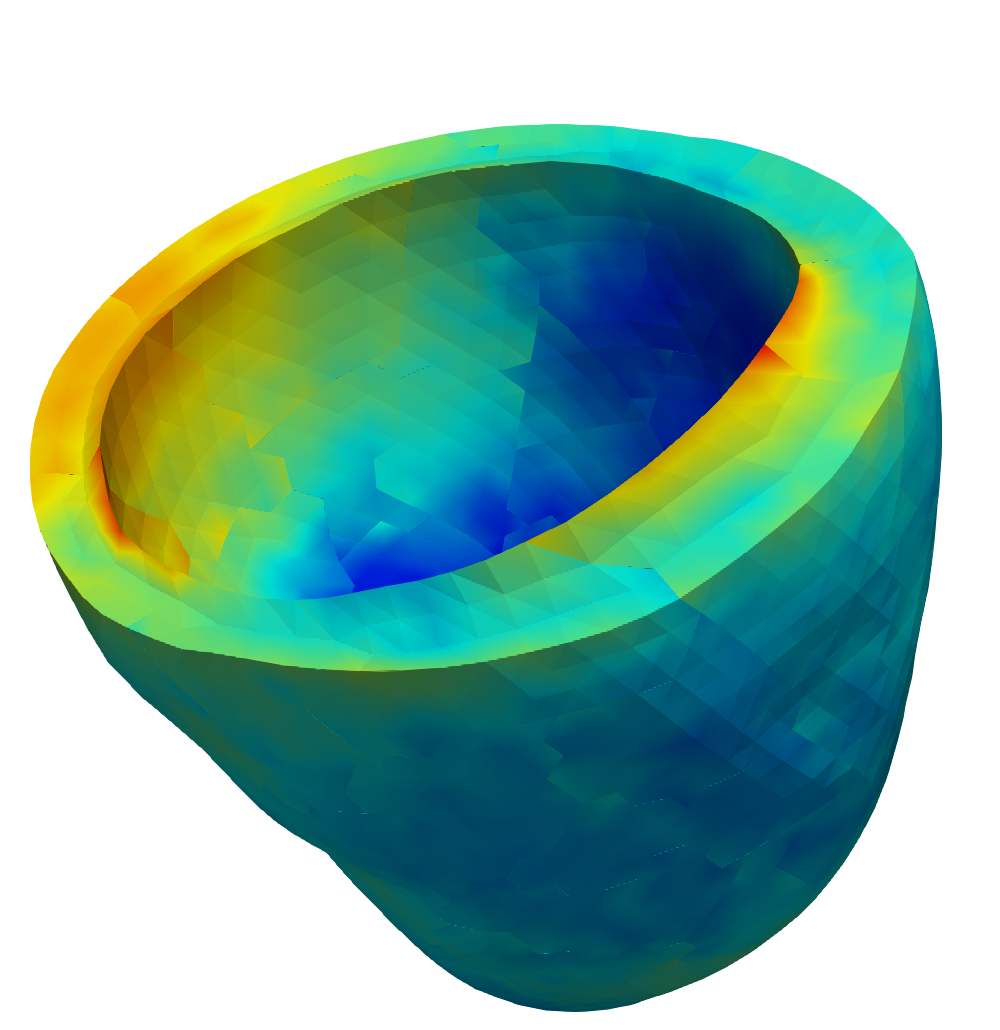}
    }
    \subfloat{
        \includegraphics[width=0.225\linewidth]{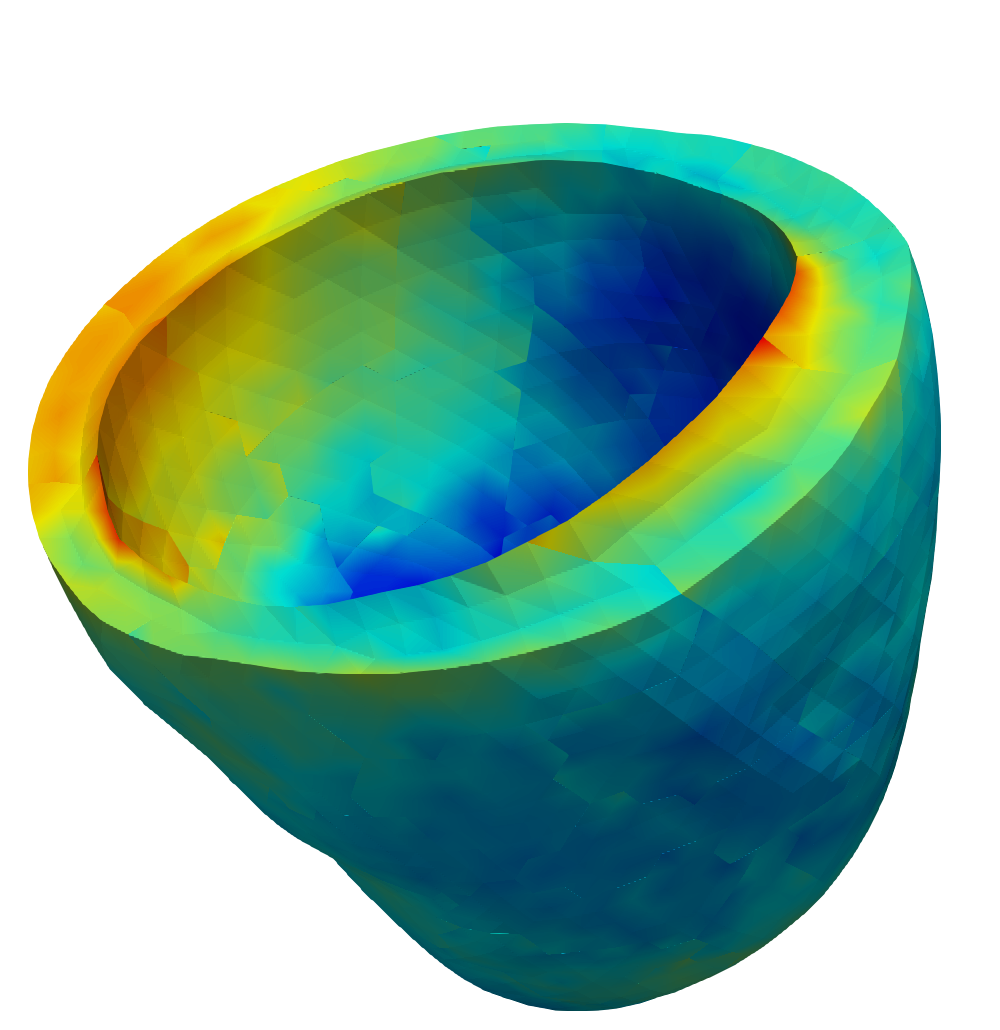}
    }
    \subfloat{
        \includegraphics[width=0.225\linewidth]{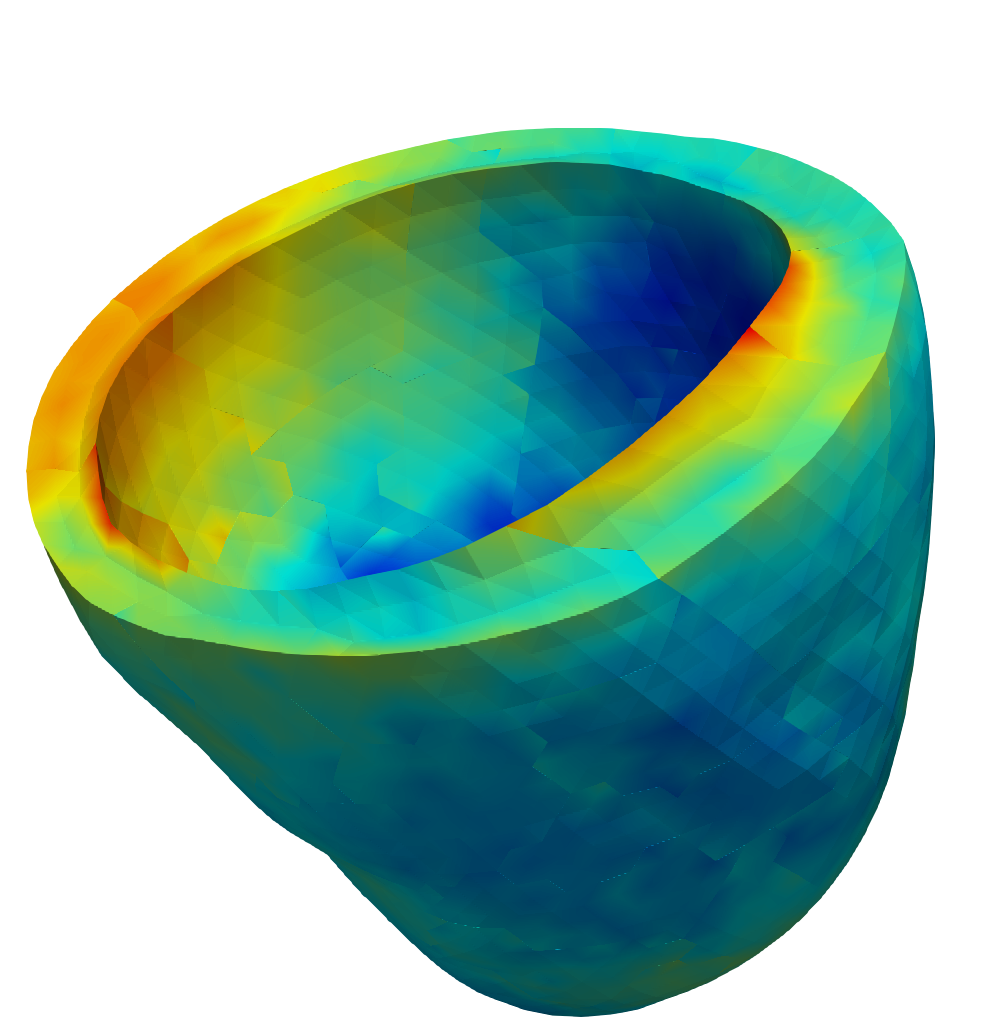}
    }\\
    \subfloat{
        \includegraphics[width=0.98\linewidth]{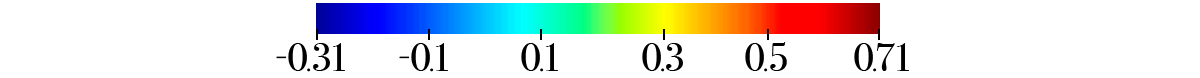}
    }
    \caption{Evolution of the Green-Lagrange strain trace $\operatorname{tr}(\mathbf{E})$ in 
    the left ventricle during the last cardiac cycle.
    Each picture is warped by the 
    displacement vector $\boldsymbol{d}$. Mild, moderate and severe systemic 
    hypertension is compared with a healthy individual.
    }
    \label{fig:sys3d-displ}
\end{figure}

Starting from $p_\text{LV}$, 
in Figures~\ref{fig:sys0d:lv-loop} and \ref{fig:sys0d:plv}
$(\mathscr{C}_\text{C})$ shows a marked increase in peak and end-systolic values, 
reflecting the increased afterload that forces the ventricle to generate higher 
pressures. The 3D--0D model confirms this trend but also reproduces a more realistic 
pressure curve, closer to patient data, 
with the distinct behaviors of pressure in the different phases more clearly visible,
together with an increase in both $p_\text{LV}$ and 
$V_\text{LV}$ 
in Figures~\ref{fig:sys3d:lv-loop} and \ref{fig:sys3d:plv}: 
this leads to a more consistent description of ventricular mechanics in hypertensive 
conditions \cite{vasc20}.
As shown in Table~\ref{table:sys-out}, 
the effect on $\text{LV}_\text{SW}$ is evident in both models, with a progressive 
raise as vascular resistance increases. In the 0D case, this increase represents 
the ventricle's attempt to compensate for the increase in afterload, even though it 
proves insufficient to sustain adequate flow. The 3D--0D model also predicts a higher 
$\text{LV}_\text{SW}$, but with a pressure-volume evolution that better matches 
clinical expectations, as the loop shows an upward shift in the minimal and maximal 
pressure, thus providing a more accurate picture of the increase in workload 
\cite{lv-hypertrophy}.
Considering $Q_\text{AV}$ 
in Figures~\ref{fig:sys0d:qav} and \ref{fig:sys3d:qav}, 
both models show a reduction with disease progression, due to compromised filling and 
limited ejection, with the 3D--0D model displaying more realistic peak values 
\cite{qav-syst}, due to 
the more accurate simulation of systolic ejection dynamics by the 3D ventricle.

As shown in Table~\ref{table:sys-out},
CI shows opposite trends in the two models: in $(\mathscr{C}_\text{C})$, it 
appears to increase slightly, mainly due to the greater pressure generation by the 
ventricle; however, this result is not consistent with clinical evidence
\cite{ci-syst}. The 3D--0D
model avoids this artefact, showing a more stable CI, which is more in line with 
the haemodynamic response expected in systemic hypertension \cite{ci-syst}.
Focusing on $p_\text{AR}^\text{SYS}$ 
in Figures~\ref{fig:sys0d:parsys} and \ref{fig:sys3d:parsys}, 
both models capture the increase in systolic and diastolic values. In 
$(\mathscr{C}_\text{C})$, this is reflected in an increase in $\text{SAP}_\text{max}$ 
and $\text{SAP}_\text{min}$ 
(Table~\ref{table:sys-out}), 
indicative of the heavier load imposed on the heart; the 3D--0D model confirms these 
trends, but with a slightly lower systolic peak. 
Quantitatively, $\text{SAP}_\text{max}$ increases steadily, while 
$\text{SAP}_\text{min}$ increases more moderately before showing a slight decrease in 
the most severe phase (Table~\ref{table:sys-out}).
With regard to $Q_\text{AR}^\text{SYS}$
in Figure~\ref{fig:sys0d:qarsys}, 
$(\mathscr{C}_\text{C})$ suggests initial stability followed by a decline, which could 
lead to a reduction in peripheral perfusion. 
Furthermore, in Figure~\ref{fig:sys3d:parpul} the 3D--0D model shows that 
$p_\text{AR}^\text{PUL}$ remains nearly 
unchanged in mild severity, consistently with clinical evidence indicating that 
systemic hypertension primarily affects the systemic circulation, with pulmonary 
pressures increasing only in advanced stages due to backward transmission, i.e. 
the rise in left heart pressures being passively transmitted to the pulmonary vessels.

Finally, active tension in Figure~\ref{fig:sys3d-ta} shows a notable peak increase 
during isovolumetric contraction and ejection, reflecting the need for higher 
contractile force to counteract the increased afterload; in addition, spatially 
resolved active tension peaks in the top portion of the left ventricle at 
$t = 0.26\,\text{s}$ and $t = 0.35\,\text{s}$, suggesting regional stress, 
potentially due to local hypertrophy or contractile heterogeneity 
\cite{active-tension-various-cond}. 
The trace of the Green-Lagrange strain tensor 
$\mathbf{E} = \dfrac{1}{2}(\mathbf{C}-\mathbf{I})$ in Figure~\ref{fig:sys3d-displ} 
shows only minor changes, with a slight increase in 
$\left|\operatorname{tr}(\mathbf{E})\right|$ maximal values. This indicates a marginal increase in 
ventricular dilation and a slightly stronger contraction during the cardiac cycle, 
consistent with the higher arterial pressure. 
The overall volumetric deformation remains close to the physiological case.

\subsection{Pulmonary Hypertension}
\label{subsec:pul-hyper}%

\begin{table}[t]
    \centering 
    \begin{tabular}{l l l l l l l l l l}
    \hline
    \multicolumn{2}{c}{} & \multicolumn{4}{c}{0D model} & \multicolumn{4}{c}{3D--0D model}\\
    \cmidrule(lr){3-6} \cmidrule(lr){7-10} 
    Output & Unit & Healthy & Mild & Moderate & Severe & Healthy & Mild & Moderate & Severe \\

    \hline
    $\text{RV}_\text{IEDV}$    & $\text{mL}/\text{m}^2$ & 68.2  &  70.3  &   73.4   &   78.1   & 71.5  &  73.6  &   78.6   &   86.0\\
    $\text{RV}_\text{IESV}$    & $\text{mL}/\text{m}^2$ & 32.6  & 35.0 &   38.7  &  44.8  & 34.2  & 36.2 &   40.2  &  46.9\\
    $\text{RV}_\text{EF}$     & \% & 52.2  & 50.2  &   47.3  &   42.6   & 52.2  & 50.9  &   48.9  &   45.5\\
    $\text{RV}_\text{SW}$    & $\text{mmHg}\cdot \text{L}$   & 1.2 & 1.3 &   1.5   &  1.6   & 1.5 & 1.6 &   1.8   &  2.0\\ 
    $\max{\nabla p_\text{rAV}}$ & mmHg & 17.3 & 19.6 &   22.6   &  27.6  & 26.4 & 28.4 & 31.0 & 35.1\\
        \hline
    \end{tabular}
    \caption{Variations of 
    selected time-independent 
    outputs in pulmonary hypertension, from both $(\mathscr{C}_\text{C})$ and 3D--0D model.
    }
    \label{table:pul-out}
\end{table}

\begin{figure}[p]
    \centering
    \subfloat[RV PV loop.\label{fig:pul0d:rv-loop}]{
        \includegraphics[height=0.25\linewidth]{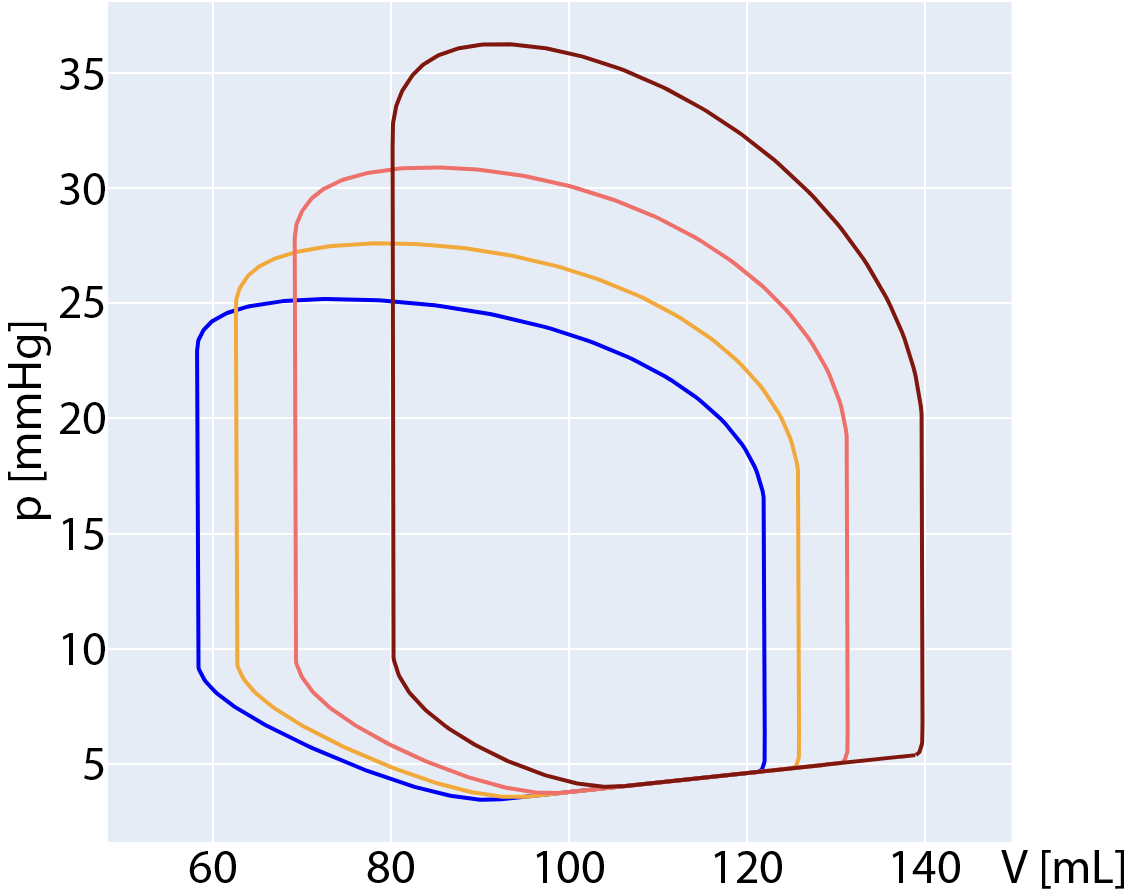}
    }
    \subfloat[RV PV loop.\label{fig:pul3d:rv-loop}]{
        \includegraphics[height=0.25\linewidth]{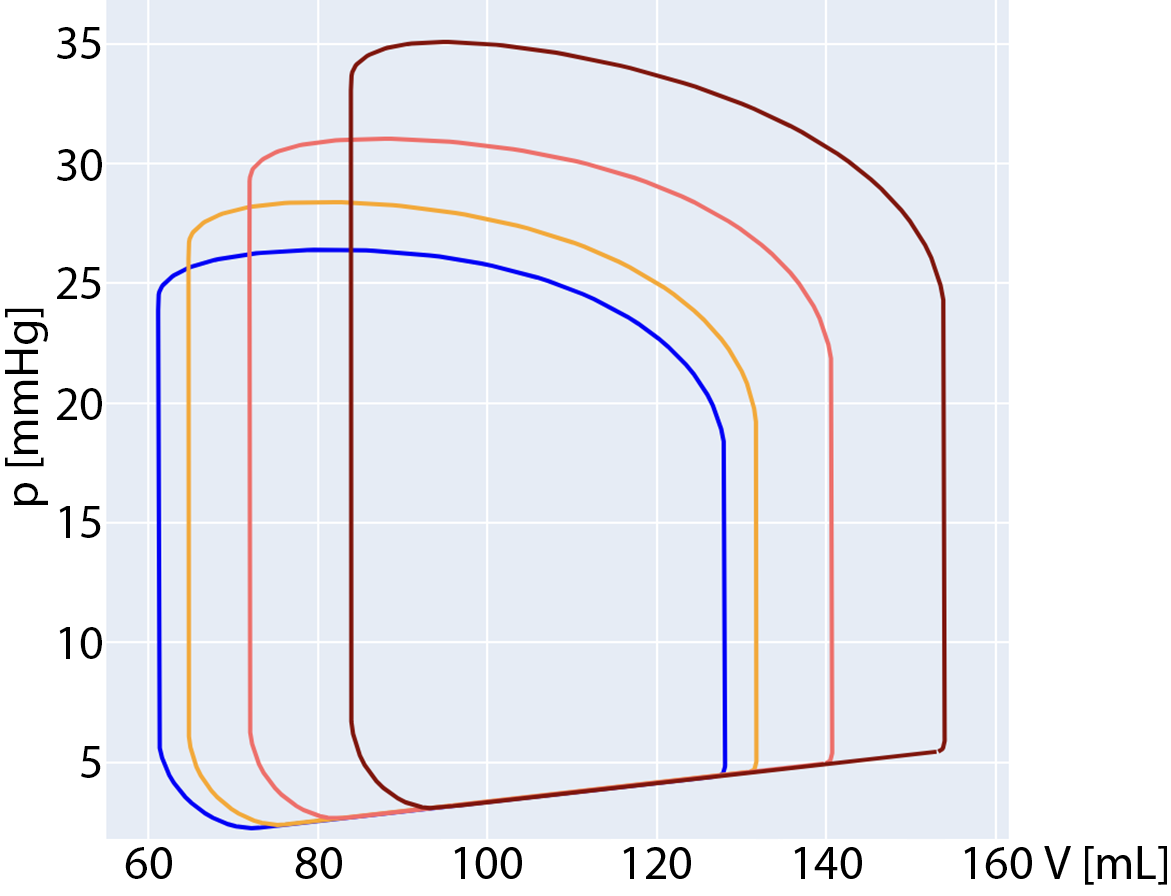}
    }\\
    \subfloat[$Q_\text{PV}$.\label{fig:pul0d:qpv}]{
        \includegraphics[width=0.32\linewidth]{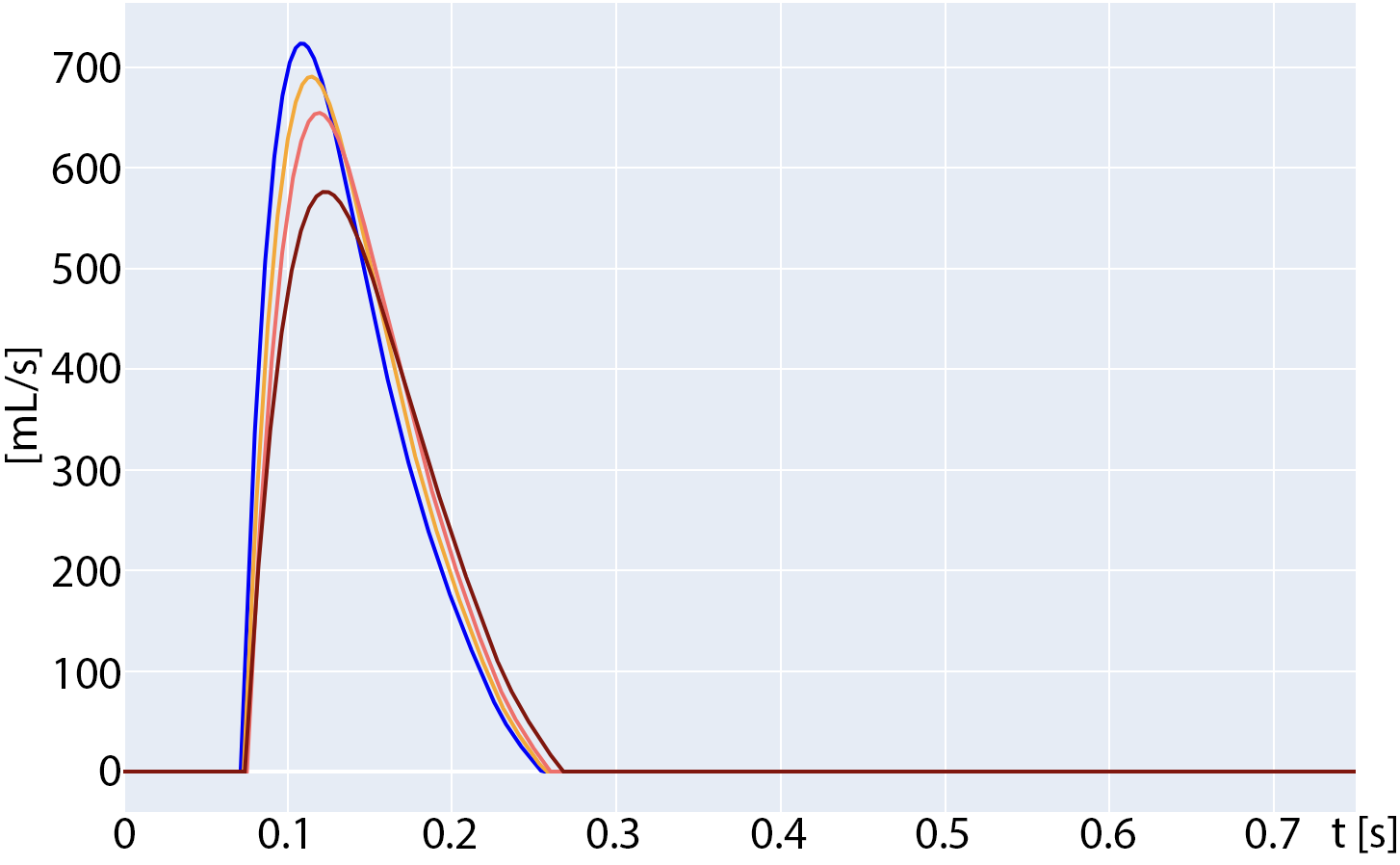}
    }
    \subfloat[$Q_\text{PV}$.\label{fig:pul3d:qpv}]{
        \includegraphics[width=0.32\linewidth]{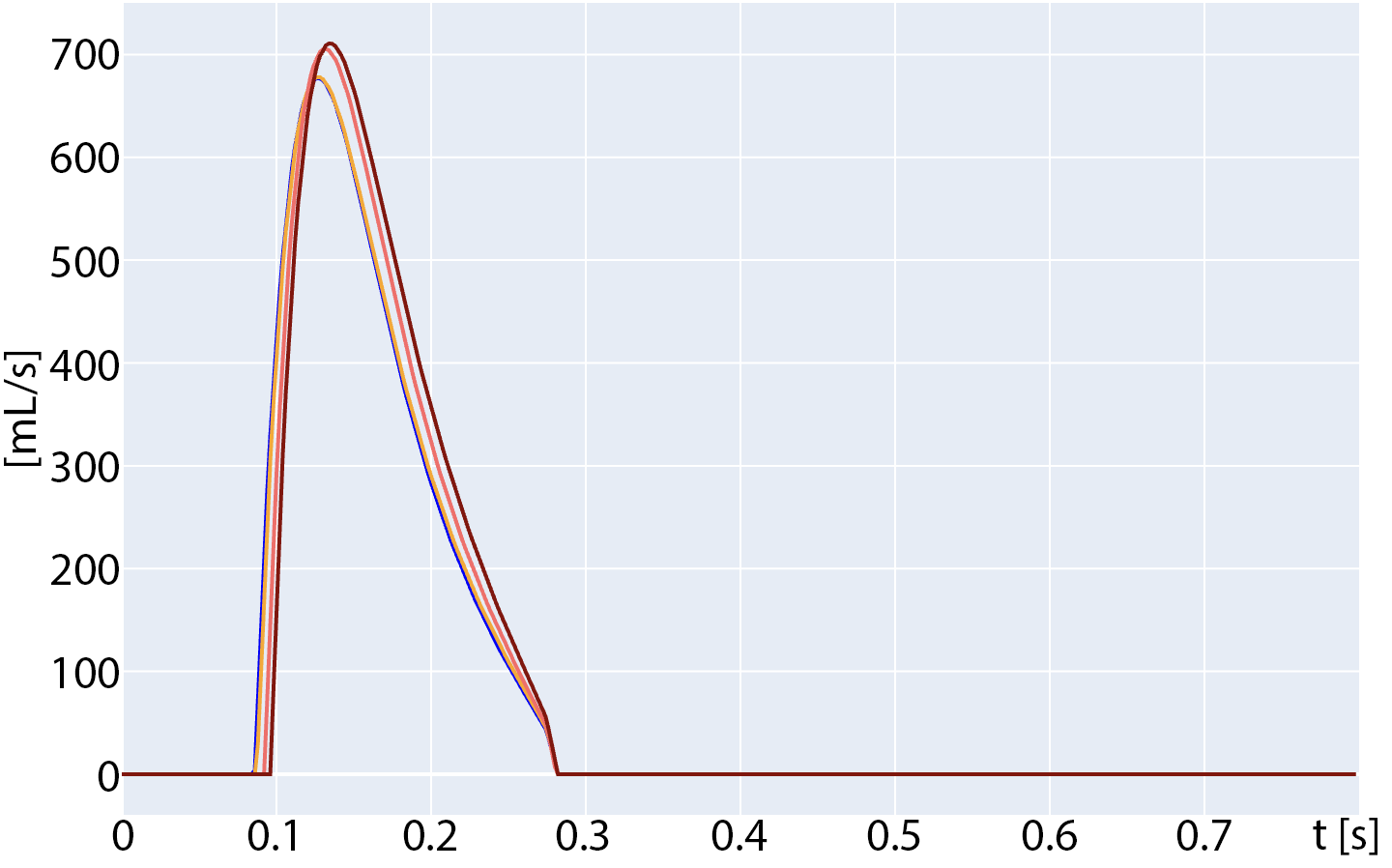}
    }\\\hspace{0.112\linewidth}
    \subfloat[$p_\text{AR}^\text{PUL}$.\label{fig:pul0d:parpul}]{
        \includegraphics[width=0.32\linewidth]{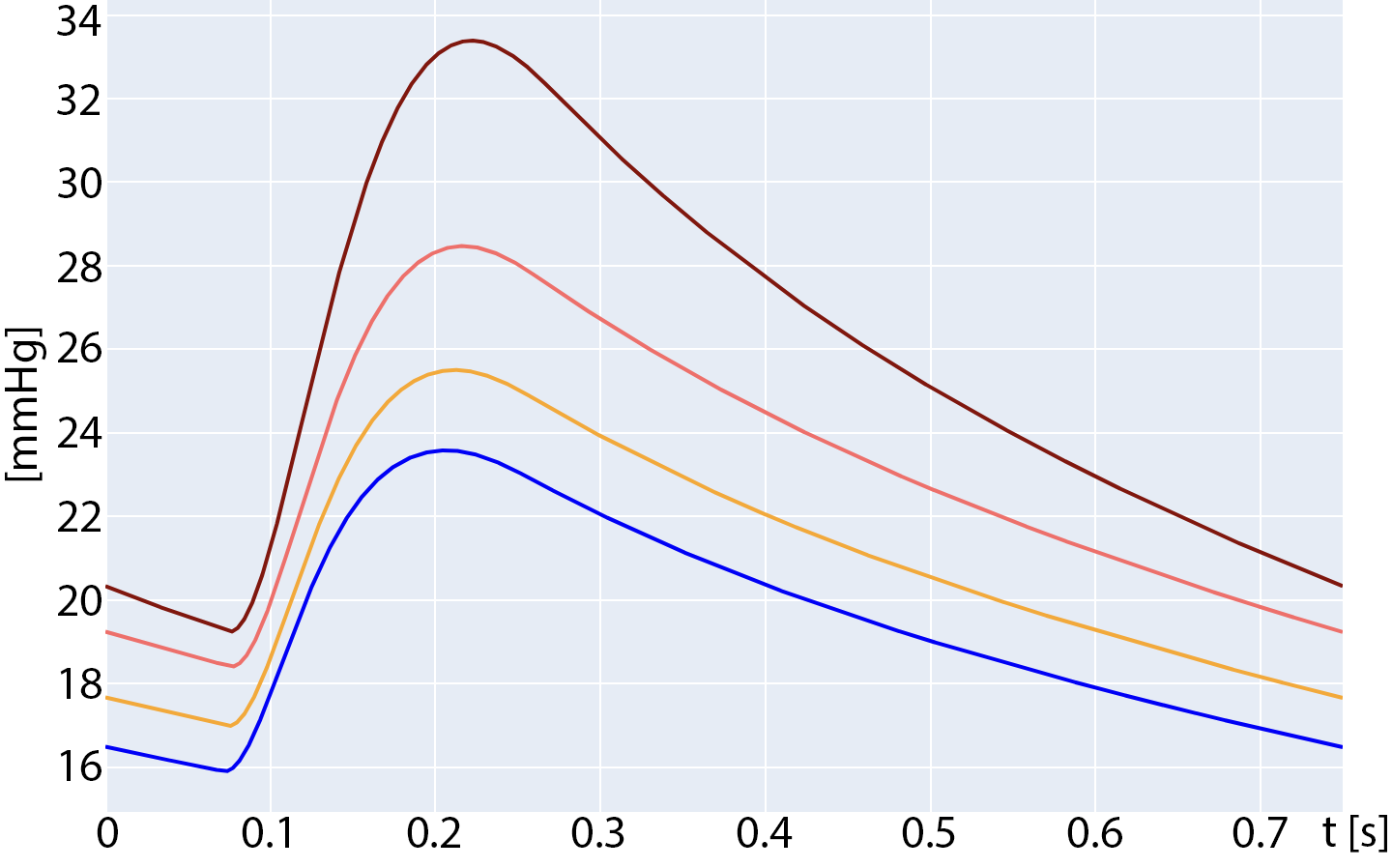}
    }
    \subfloat[$p_\text{AR}^\text{PUL}$.\label{fig:pul3d:parpul}]{
        \includegraphics[width=0.32\linewidth]{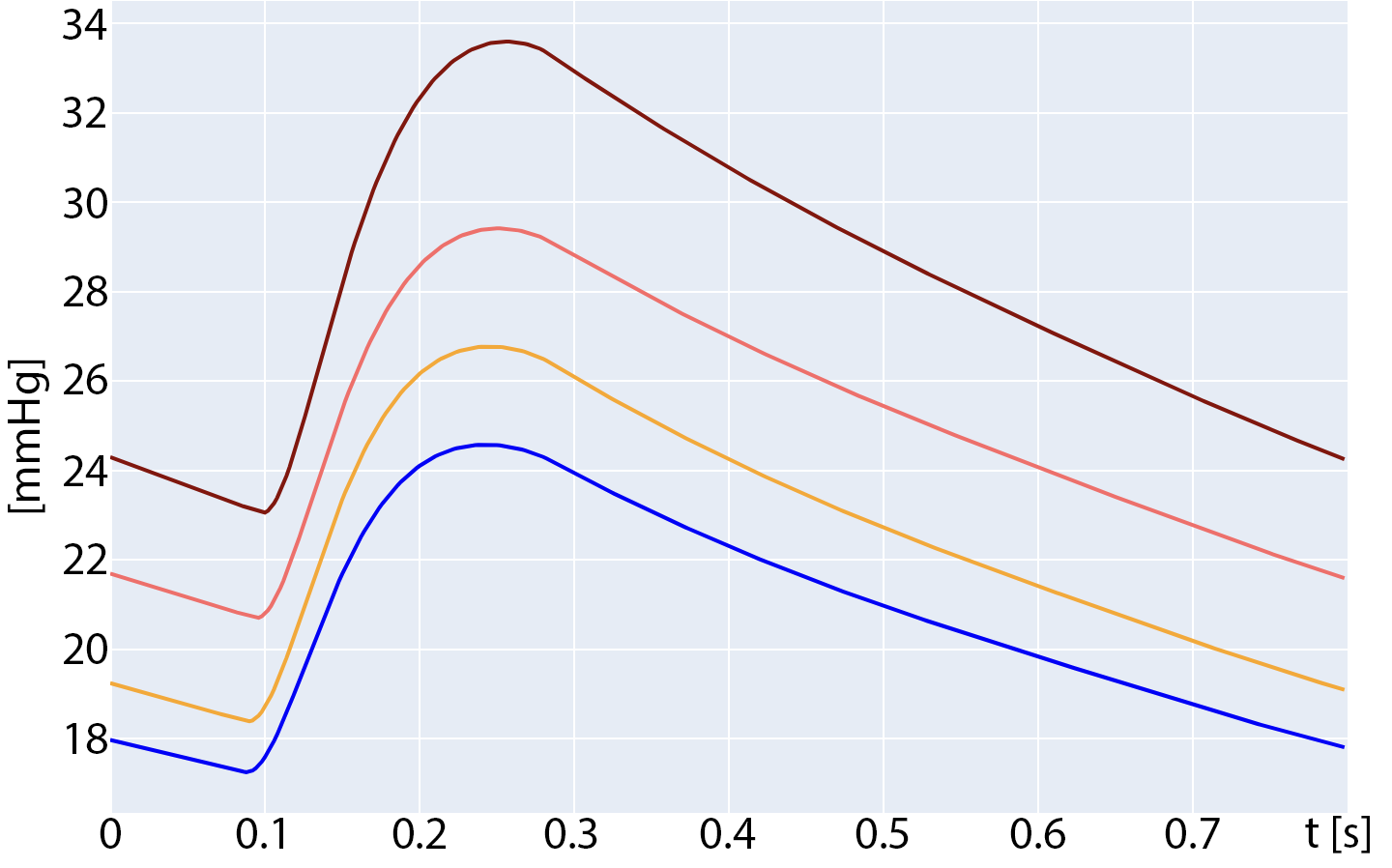}
    }
    \begin{subfigure}[b]{0.12\textwidth}
        \raisebox{0.6\height}{\includegraphics[width=\textwidth]{legend_variables.png}}
    \end{subfigure}\\
    \subfloat[$Q_\text{AR}^\text{PUL}$.\label{fig:pul0d:qarpul}]{
        \includegraphics[width=0.32\linewidth]{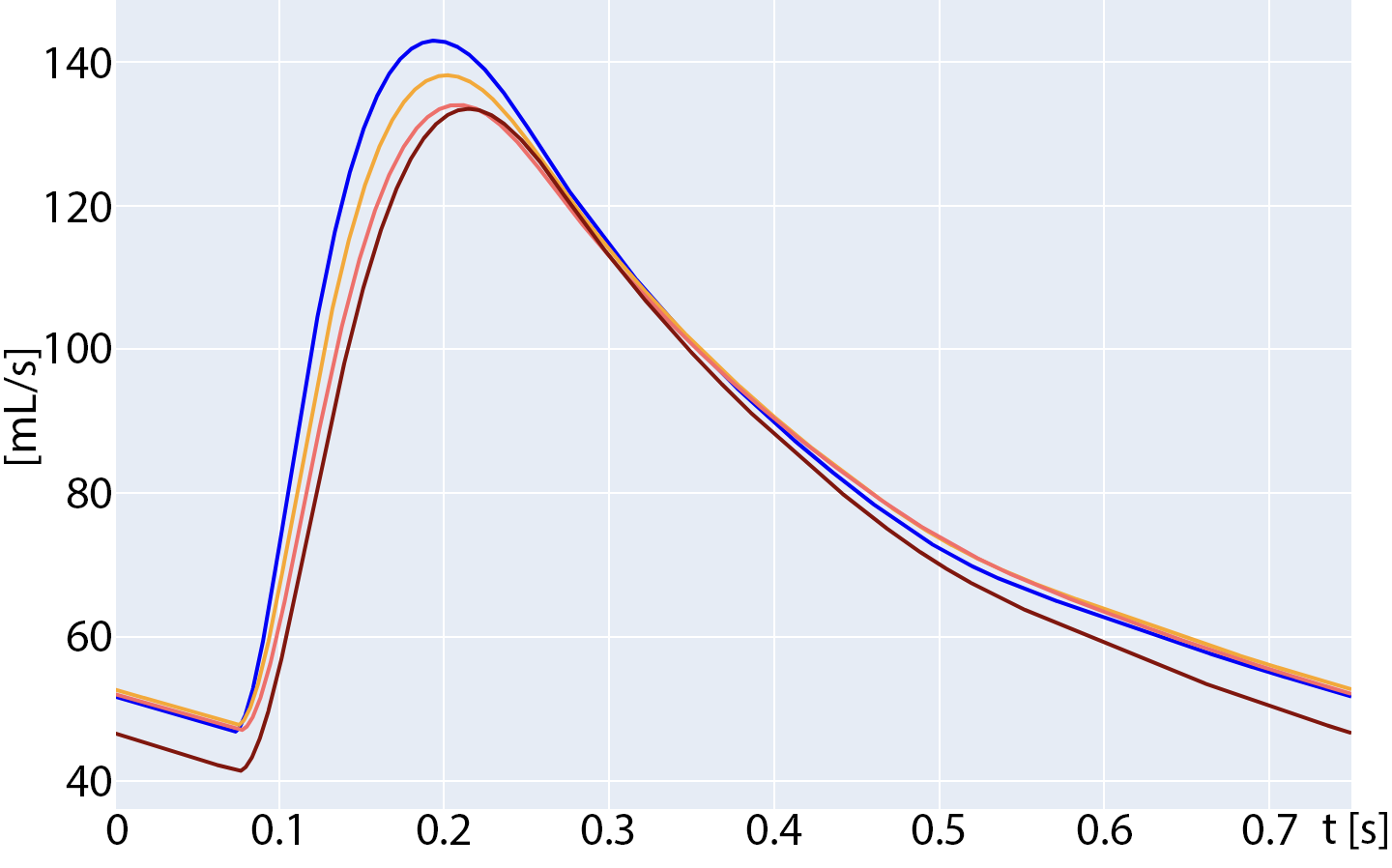}
    }
    \subfloat[$Q_\text{AR}^\text{PUL}$.\label{fig:pul3d:qarpul}]{
        \includegraphics[width=0.32\linewidth]{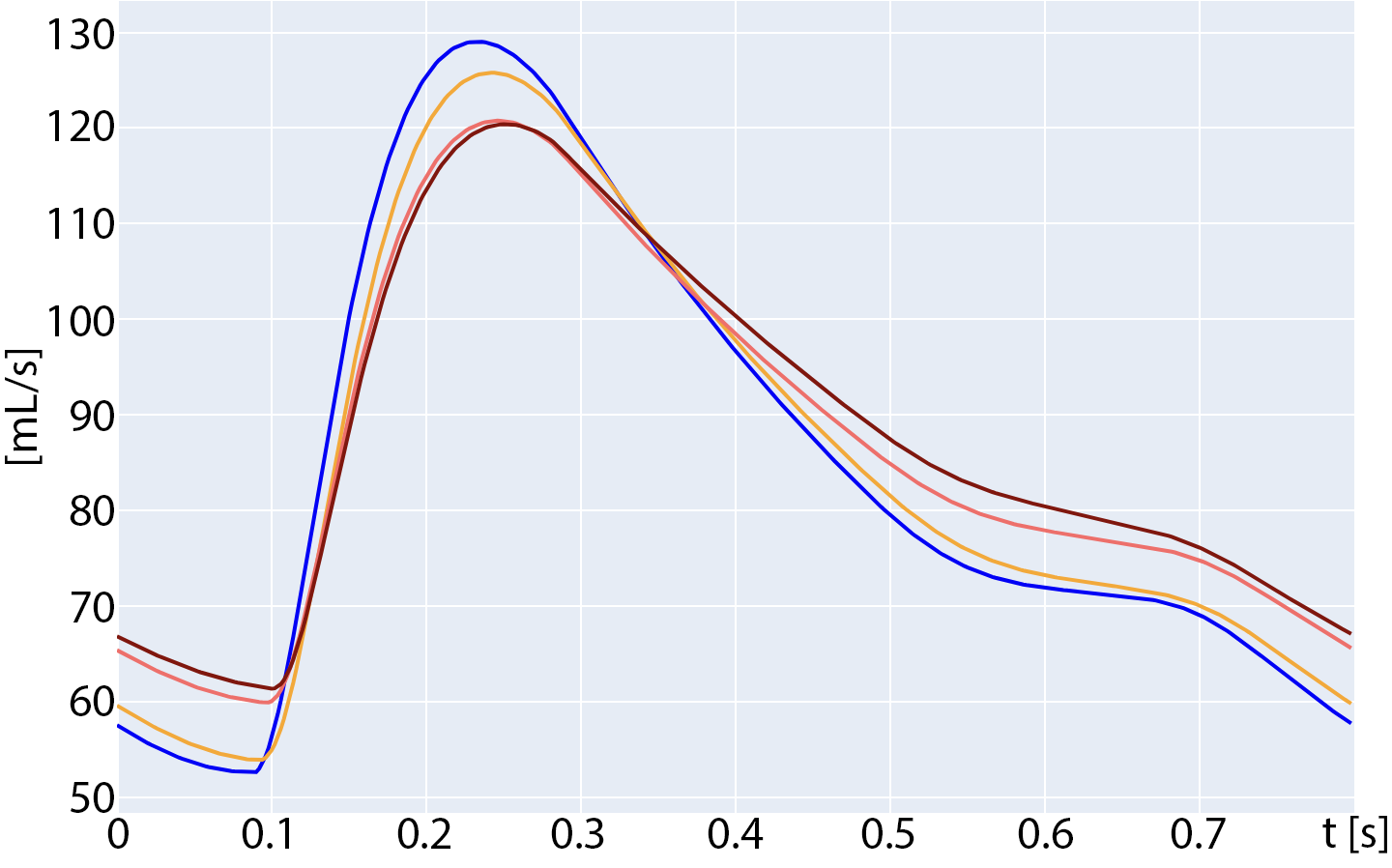}
    }\\
    \subfloat[$Q_\text{VEN}^\text{SYS}$.\label{fig:pul0d:qvensys}]{
        \includegraphics[width=0.32\linewidth]{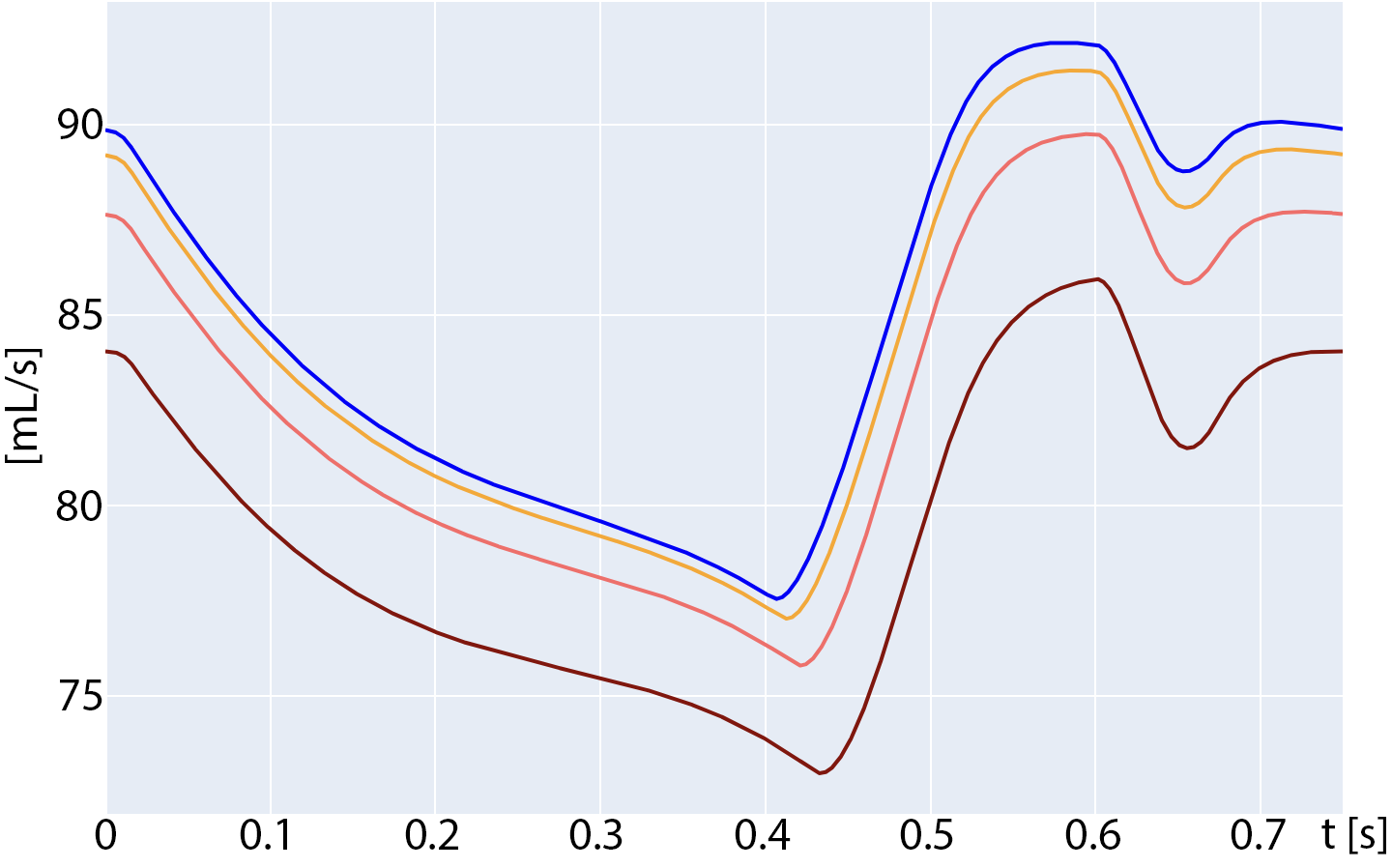}
    }
    \subfloat[$Q_\text{VEN}^\text{SYS}$.\label{fig:pul3d:qvensys}]{
        \includegraphics[width=0.32\linewidth]{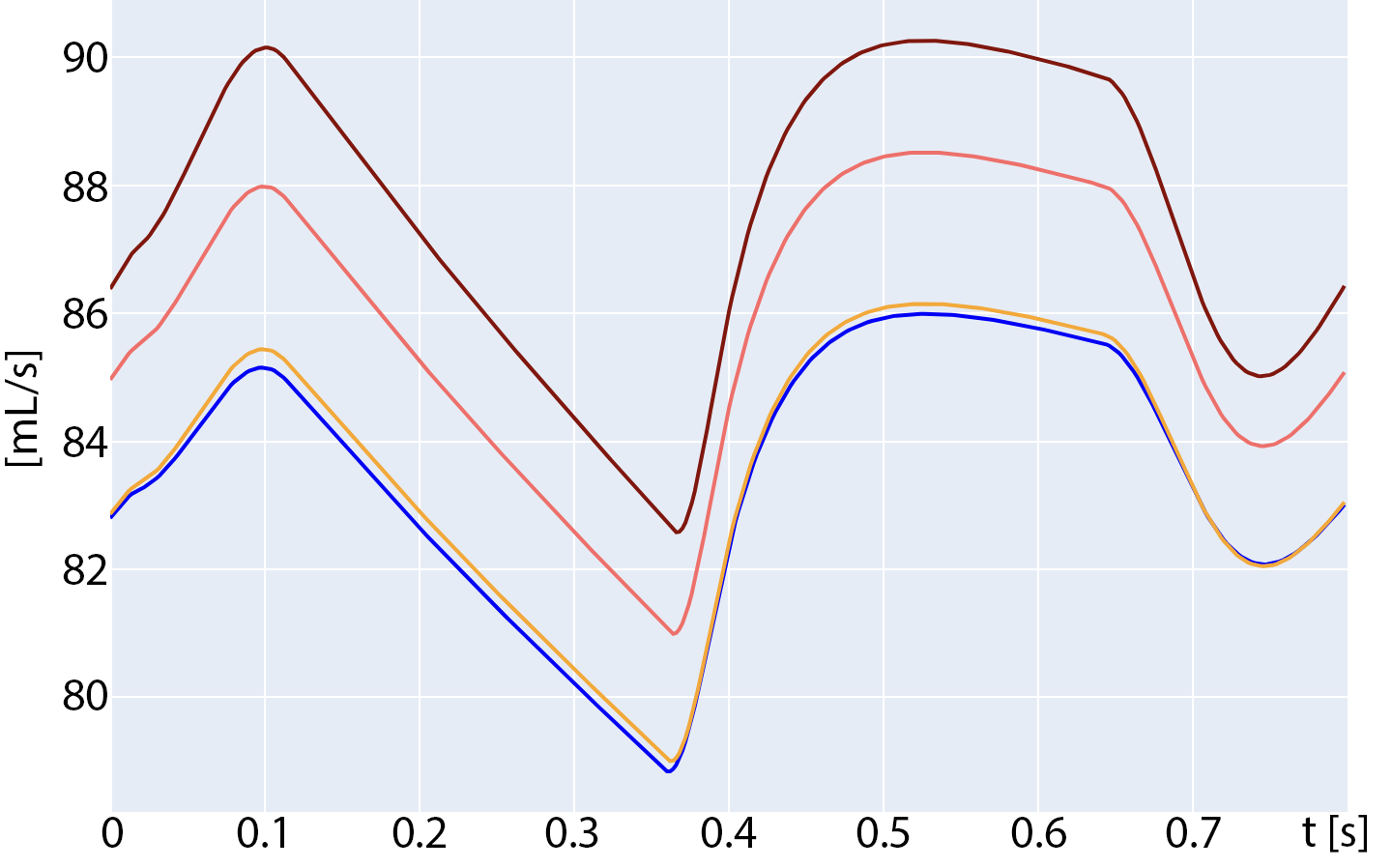}
    }
    \caption{Time-dependent variables from $(\mathscr{C}_\text{C})$ (on the left) and
    from the 3D--0D model (on the right).
    Mild, moderate and severe pulmonary hypertension 
    (in orange, light red and dark red) is compared with a 
    healthy individual (in blue).
    }
    \label{fig:pul-var}
\end{figure}

Considering the right ventricular PV loop
in Figures~\ref{fig:pul0d:rv-loop} and \ref{fig:pul3d:rv-loop}, 
both end-systolic and end-diastolic pressures increase in $(\mathscr{C}_\text{C})$, 
and the ventricle becomes stiffer. The 3D--0D loop shows a similar overall trend, but 
the 3D representation allows for a more physiologically accurate description of ejection 
dynamics and regional mechanics \cite{rv-loop-pul-hyper}.

Looking at right ventricular volumes, 
as shown in Table~\ref{table:pul-out}, 
both models detect a marked increase in $\text{RV}_\text{IEDV}$, reflecting higher 
filling pressures and ventricular dilatation. In $(\mathscr{C}_\text{C})$, 
$\text{RV}_\text{IESV}$ also increases, leading to a reduction in $\text{RV}_\text{EF}$ 
and impaired $Q_\text{PV}$
in Figure~\ref{fig:pul0d:qpv}. 
The 3D--0D model confirms the increase in $\text{RV}_\text{IESV}$ and the decrease in 
$\text{RV}_\text{EF}$, highlighting the progressive inefficiency of right ventricular 
emptying. In particular, unlike the 0D model, $Q_\text{PV}$ remains relatively stable 
in the 3D--0D case
(Figure~\ref{fig:pul3d:qpv}), 
suggesting that the right ventricle compensates through increased pressure generation to 
preserve stroke volume despite increased resistance.

As shown in Table~\ref{table:pul-out}, 
$\text{RV}_\text{SW}$ increases in both models, reflecting the increased energy demand 
imposed by pulmonary hypertension. While $(\mathscr{C}_\text{C})$ quantifies this 
increase as a compensatory response to maintain flow, the 3D--0D model integrates this 
increase with the spatial distribution of ventricular contraction, providing a more 
realistic assessment of workload \cite{rv-sw-pul-hyper}.
$p_{\text{AR}}^{\text{PUL}}$ increases in both models
(Figures~\ref{fig:pul0d:parpul} and \ref{fig:pul3d:parpul}), 
demonstrating the afterload imposed on the right ventricle;
$Q_{\text{AR}}^{\text{PUL}}$ decreases in 0D due to limited ventricular ejection 
(Figure~\ref{fig:pul0d:qarpul}), 
whereas in the 3D--0D model, the flow curve is flatter and more prolonged, reflecting 
a delayed peak and a more physiological right ventricular ejection under pulmonary 
hypertension.
\cite{qarpul-pul-hyper}
(Figure~\ref{fig:pul3d:qarpul}). 
In Figures~\ref{fig:pul0d:qvensys} and \ref{fig:pul3d:qvensys}, 
$Q_\text{VEN}^\text{SYS}$ is compromised in the 0D model, potentially causing congestion, 
while in 3D--0D it increases as a compensatory mechanism to attenuate pressure 
variations, although this may ultimately fail in cases of further aggravated conditions.
Finally, $\max \nabla p_{\text{rAV}}$ rises in both models
(Table~\ref{table:pul-out}), 
reflecting reduced filling efficiency and reduced overall cardiac performance. 
Interventricular interactions may further compromise left ventricular filling, an 
effect present in 0D and less pronounced in 3D--0D due to more accurate spatial 
resolution.

As pulmonary hypertension mainly affects right ventricular dynamics, only marginal 
effects are expected on the left ventricle in the absence of left-sided dysfunction
\cite{pulhyper:guidelines,qarpul-pul-hyper};
then, no figures related to 3D left ventricular variables are included.

\subsection{Renovascular Hypertension with Secondary Pulmonary Hypertension}
\label{subsec:reno-hyper}%

\begin{table}[t!]
    \centering 
    \begin{tabular}{l l l l l l l l l l}
    \hline
    \multicolumn{2}{c}{} & \multicolumn{4}{c}{0D model} & \multicolumn{4}{c}{3D--0D model}\\
    \cmidrule(lr){3-6} \cmidrule(lr){7-10} 
    Output & Unit & Healthy & Mild & Moderate & Severe & Healthy & Mild & Moderate & Severe \\
    \hline 
    $\text{LV}_\text{IEDV}$ &$\text{mL}/\text{m}^2$ &59.7 &61.6 &69.0 &79.6     &69.9 &71.4 &75.4 &82.3\\
    $\text{LV}_\text{IESV}$ &$\text{mL}/\text{m}^2$ &23.9 &25.8 &30.7 &38.7     &29.5 &31.0 &32.8 &34.4\\
    $\text{RV}_\text{IEDV}$ &$\text{mL}/\text{m}^2$ &68.2 &69.0 &75.3 & 83.7        &71.5 &72.4 &79.7 & 93.8\\
    $\text{RV}_\text{IESV}$ &$\text{mL}/\text{m}^2$ &32.6 &33.5 &37.2 &43.2         &34.2 &34.8 &38.5 &45.2\\
    $\text{LV}_\text{SW}$    & $\text{mmHg}\cdot \text{L}$   & 5.9 & 6.4 &   8.4   &  11.6   & 6.3 & 6.8 &   7.8   &  9.5\\
    $\text{RV}_\text{SW}$    & $\text{mmHg}\cdot \text{L}$   & 1.1 & 1.3 &   1.6   &  2.1   & 1.5 & 1.5 &   1.9   &  2.7 \\ 
    \hline
    \end{tabular}
    \caption{Variations of 
    selected time-independent 
    outputs in renovascular hypertension with secondary pulmonary 
    hypertension, from both $(\mathscr{C}_\text{C})$ and 3D--0D model.
    }
    \label{table:reno-out}
\end{table}

\begin{figure}[p]
    \centering
    \subfloat[LV PV loop.\label{fig:reno0d:lv-loop}]{
        \includegraphics[height=0.24\linewidth]{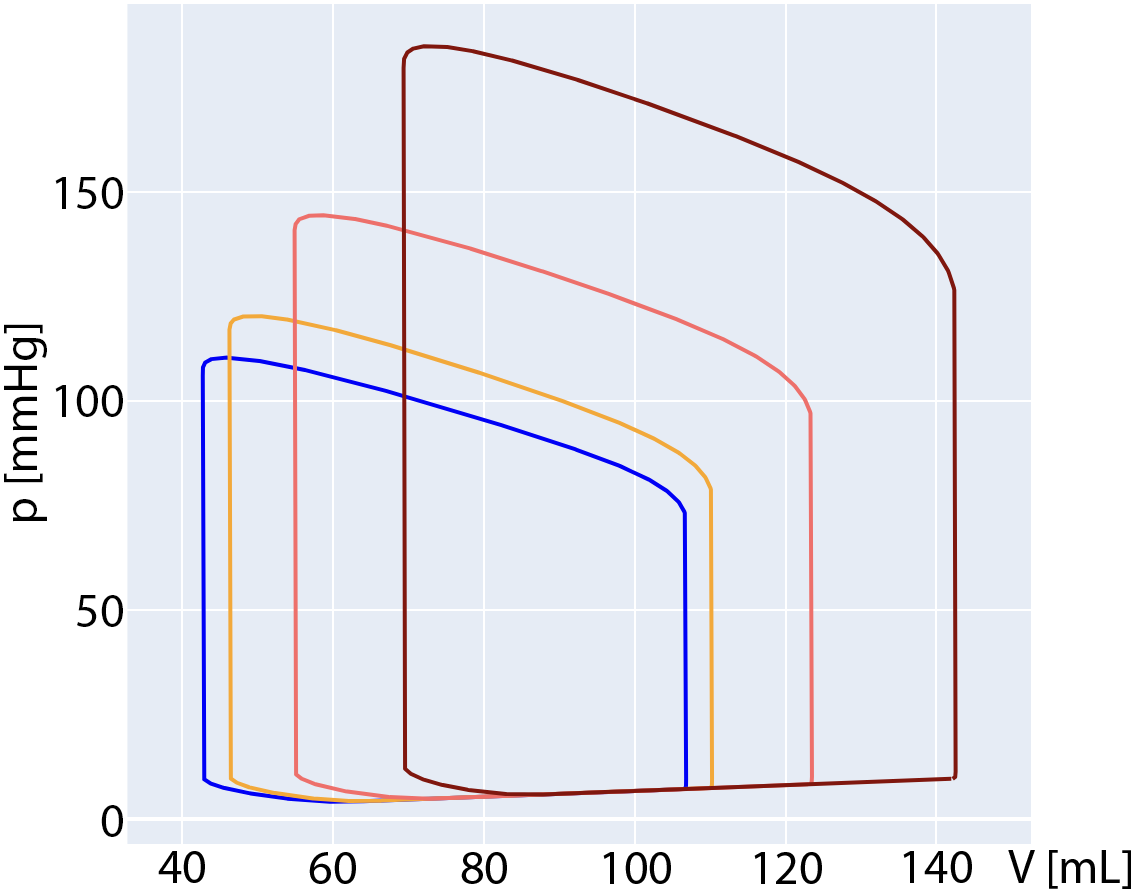}
    }
    \subfloat[LV PV loop.\label{fig:reno3d:lv-loop}]{
        \includegraphics[height=0.24\linewidth]{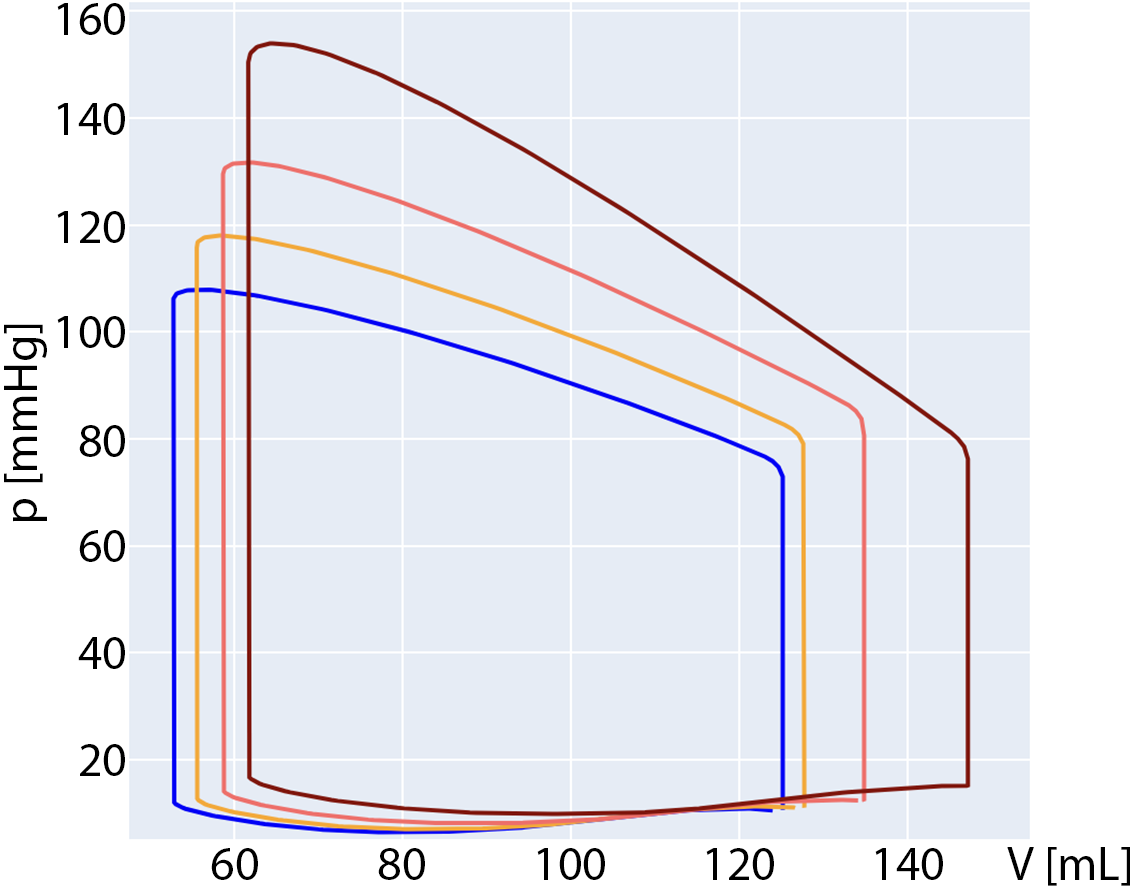}
    }\\
    \subfloat[RV PV loop.\label{fig:reno0d:rv-loop}]{
        \includegraphics[height=0.24\linewidth]{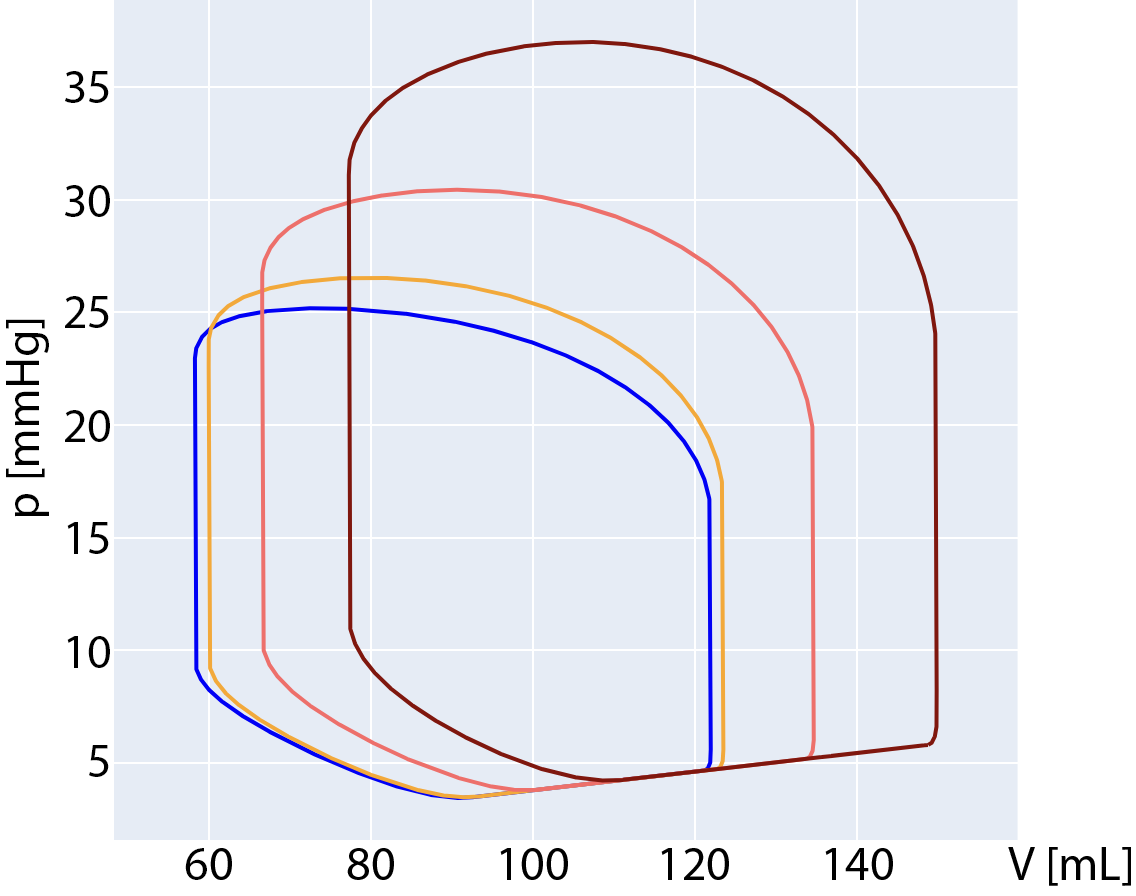}
    }
    \subfloat[RV PV loop.\label{fig:reno3d:rv-loop}]{
        \includegraphics[height=0.24\linewidth]{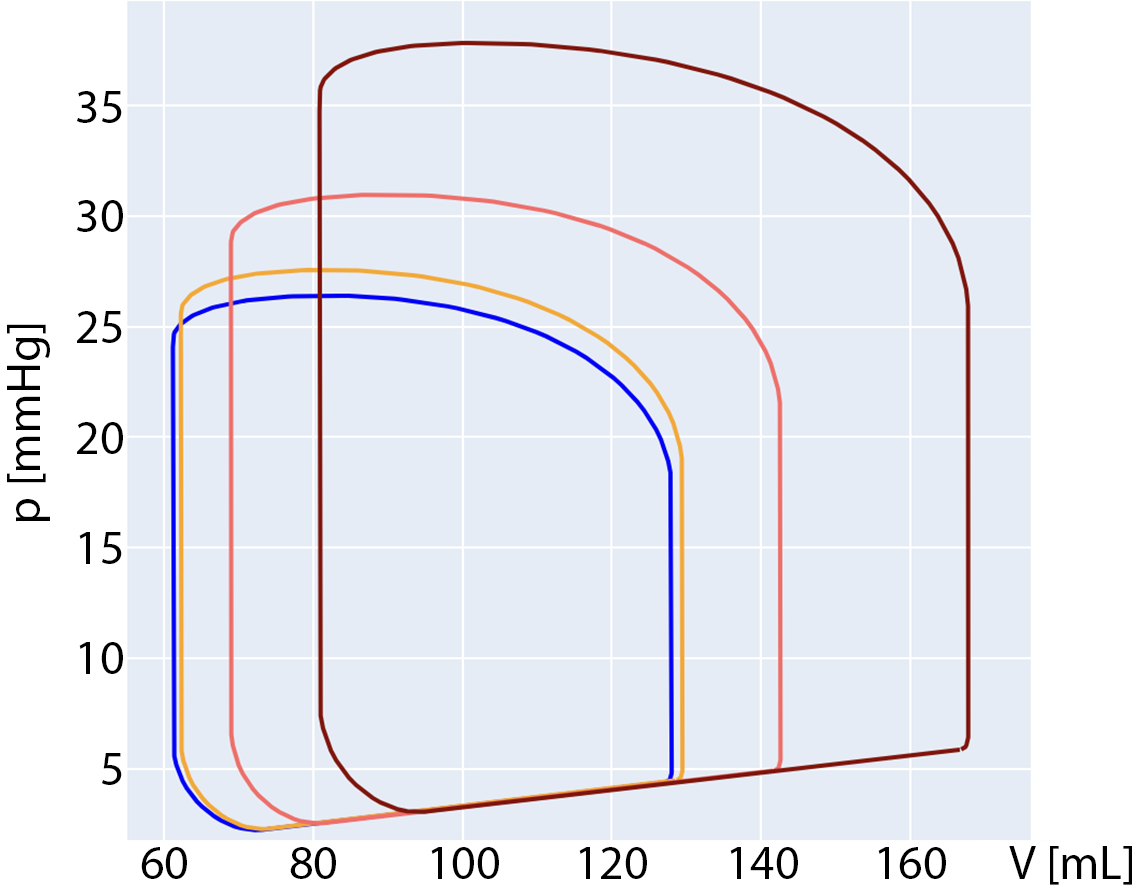}
    }\\\hspace{0.112\linewidth}
    \subfloat[$p_\text{AR}^\text{SYS}$.\label{fig:reno0d:parsys}]{
        \includegraphics[width=0.31\linewidth]{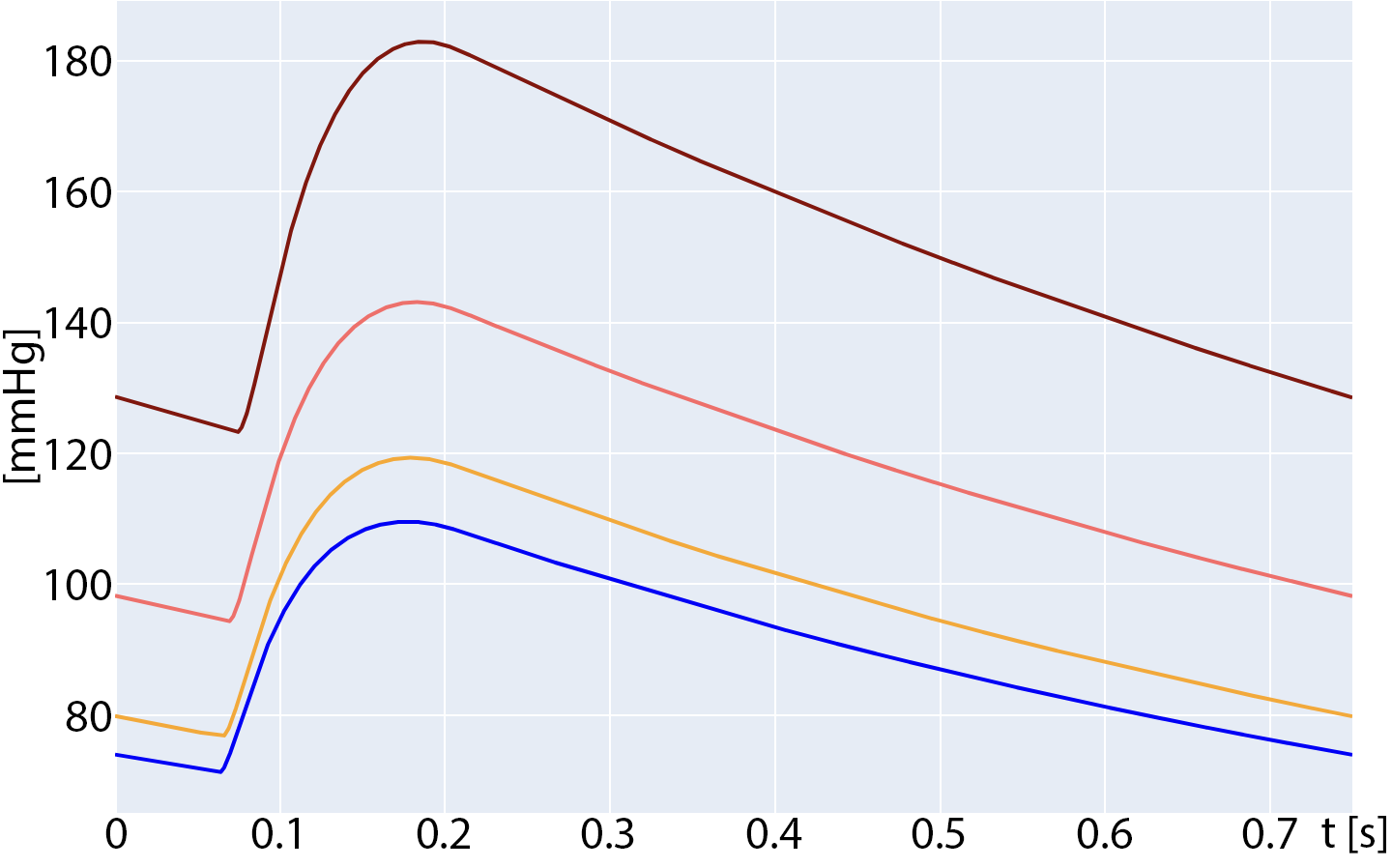}
    }
    \subfloat[$p_\text{AR}^\text{SYS}$.\label{fig:reno3d:parsys}]{
        \includegraphics[width=0.31\linewidth]{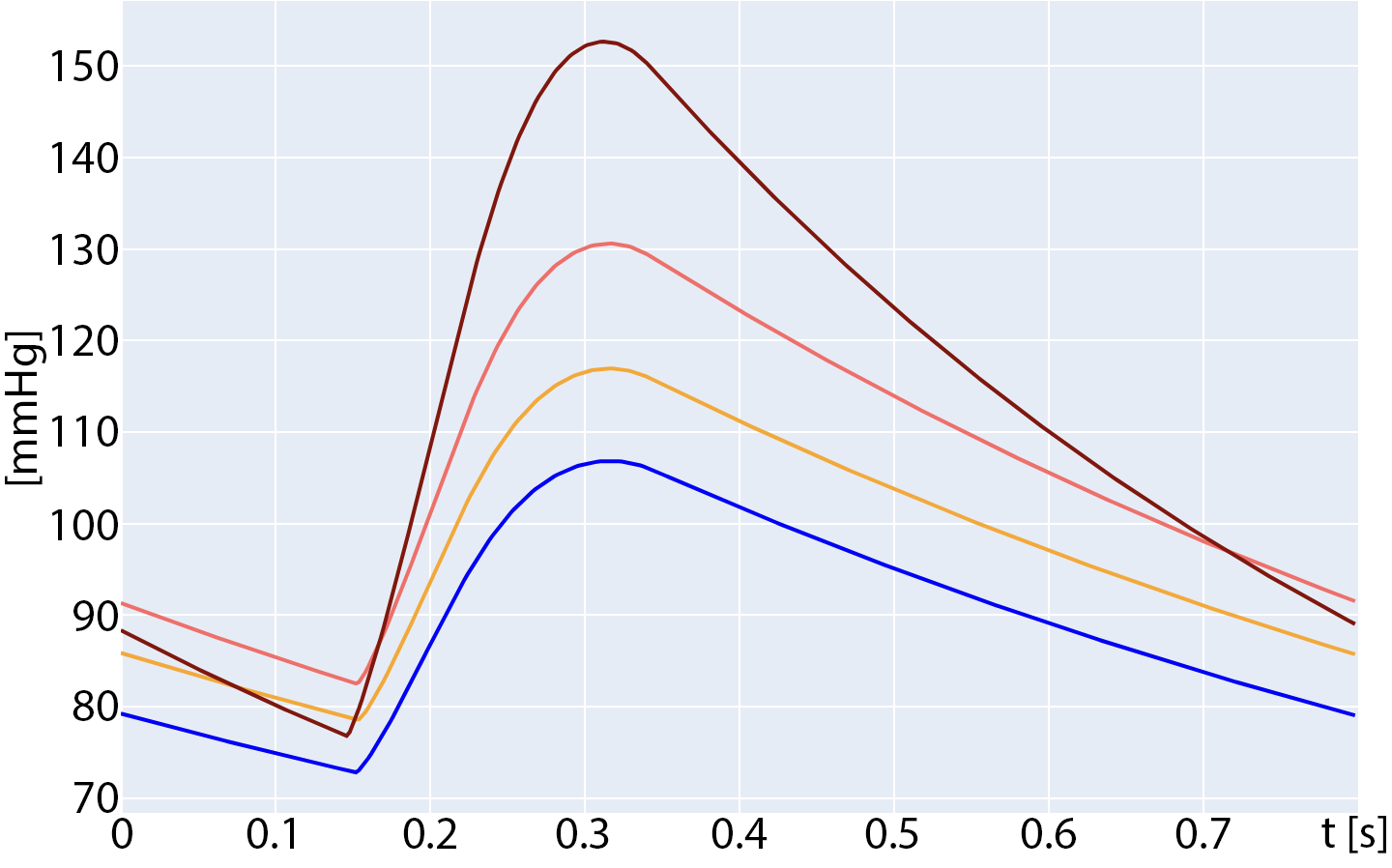}
    }
    \begin{subfigure}[b]{0.12\textwidth}
        \raisebox{0.6\height}{\includegraphics[width=\textwidth]{legend_variables.png}}
    \end{subfigure}\\
    \subfloat[$Q_\text{AR}^\text{SYS}$.\label{fig:reno0d:qarsys}]{
        \includegraphics[width=0.31\linewidth]{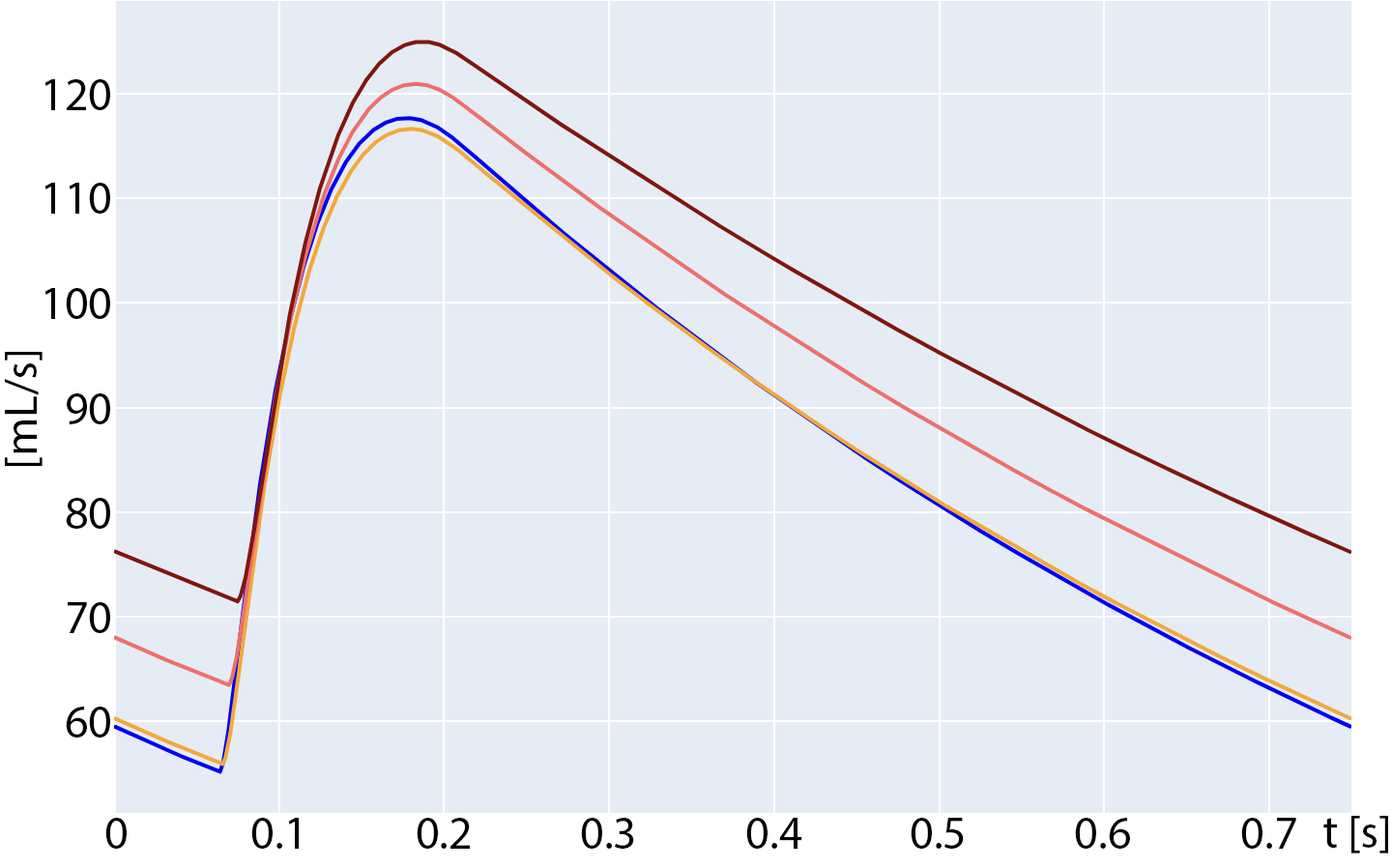}
    }
    \subfloat[$Q_\text{AR}^\text{SYS}$.\label{fig:reno3d:qarsys}]{
        \includegraphics[width=0.31\linewidth]{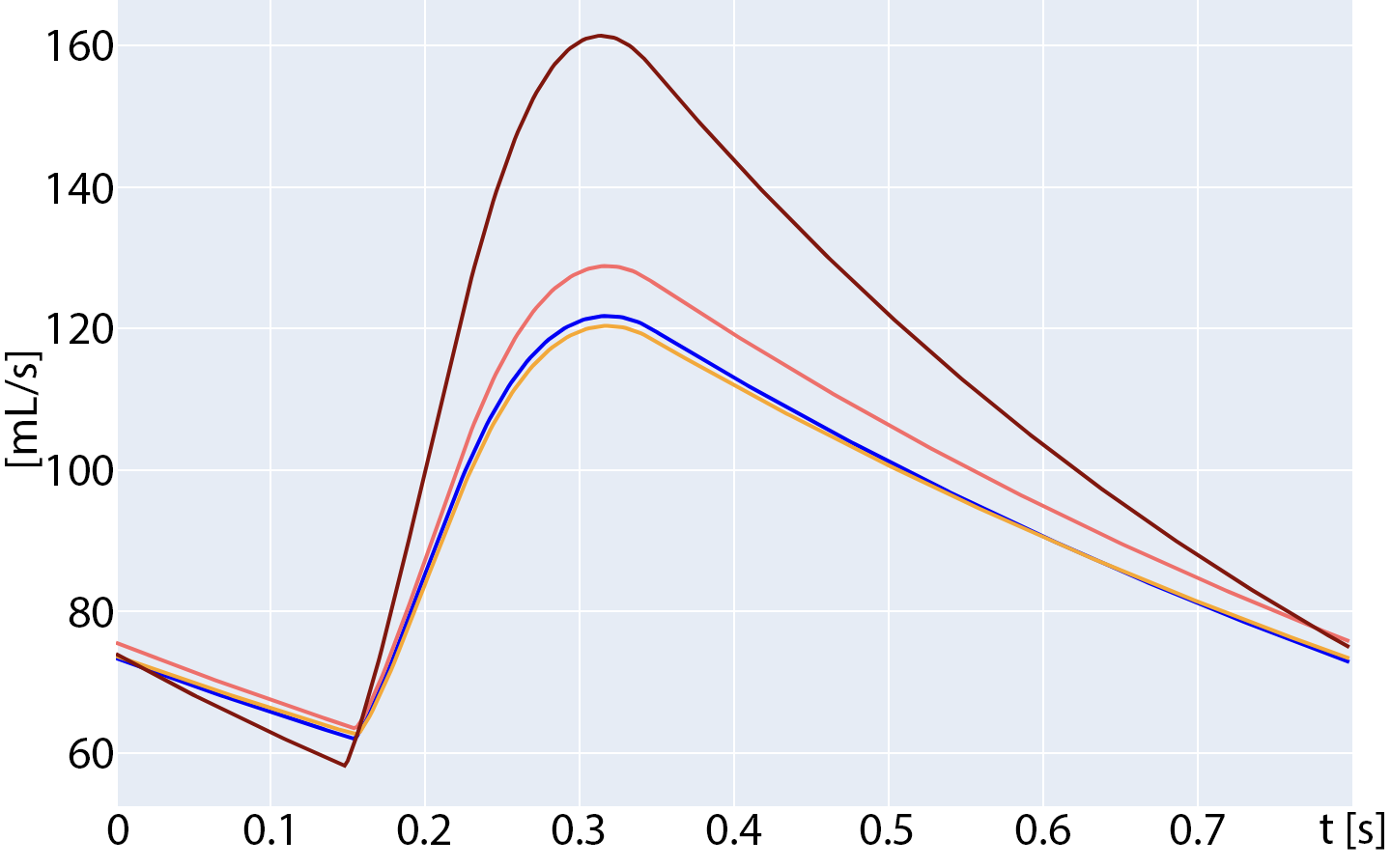}
    }\\
    \subfloat[$p_\text{AR}^\text{PUL}$.\label{fig:reno0d:parpul}]{
        \includegraphics[width=0.31\linewidth]{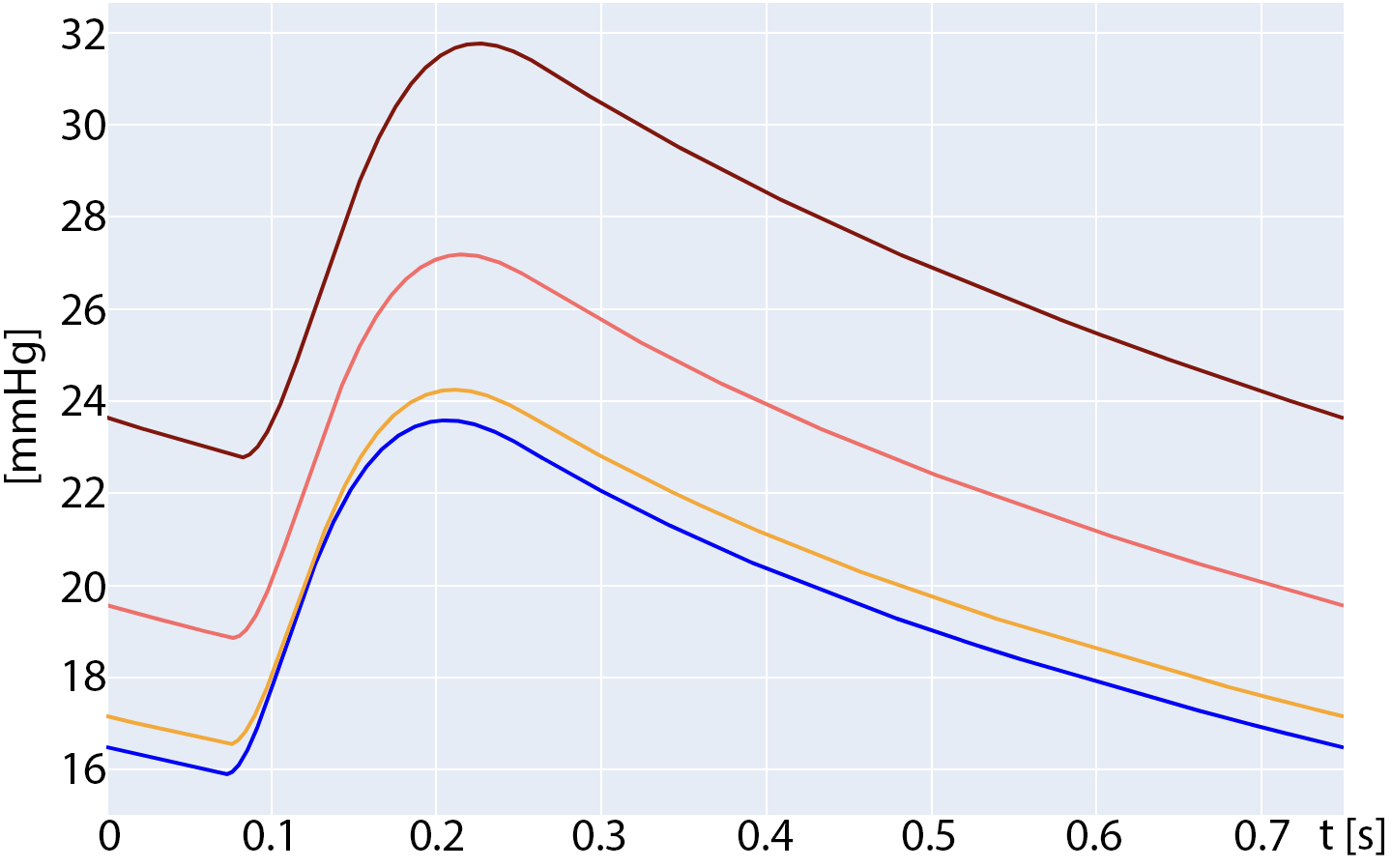}
    }
    \subfloat[$p_\text{AR}^\text{PUL}$.\label{fig:reno3d:parpul}]{
        \includegraphics[width=0.31\linewidth]{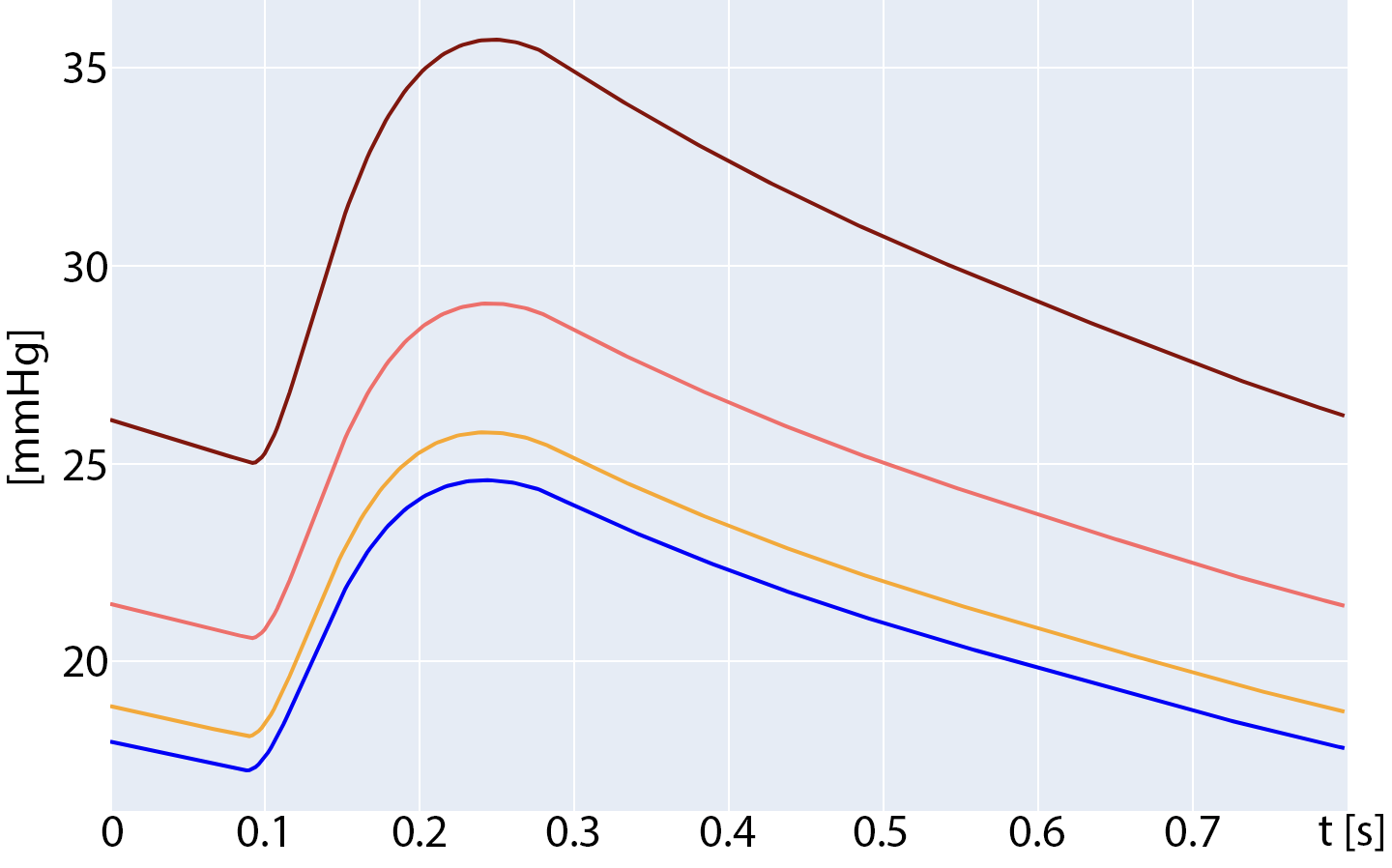}
    }
    \caption{Time-dependent variables from $(\mathscr{C}_\text{C})$ (on the left) and
    from the 3D--0D model (on the right).
    Mild, moderate and severe 
    renovascular hypertension 
    with secondary pulmonary hypertension 
    (in orange, light red and dark red) is compared with a 
    healthy individual (in blue).
    }
    \label{fig:reno-var}
\end{figure}

\begin{figure}[t!]
    \centering
    \subfloat[$p_\text{C}^\text{PUL}$.\label{fig:reno0d:pcpul}]{
        \includegraphics[width=0.31\linewidth]{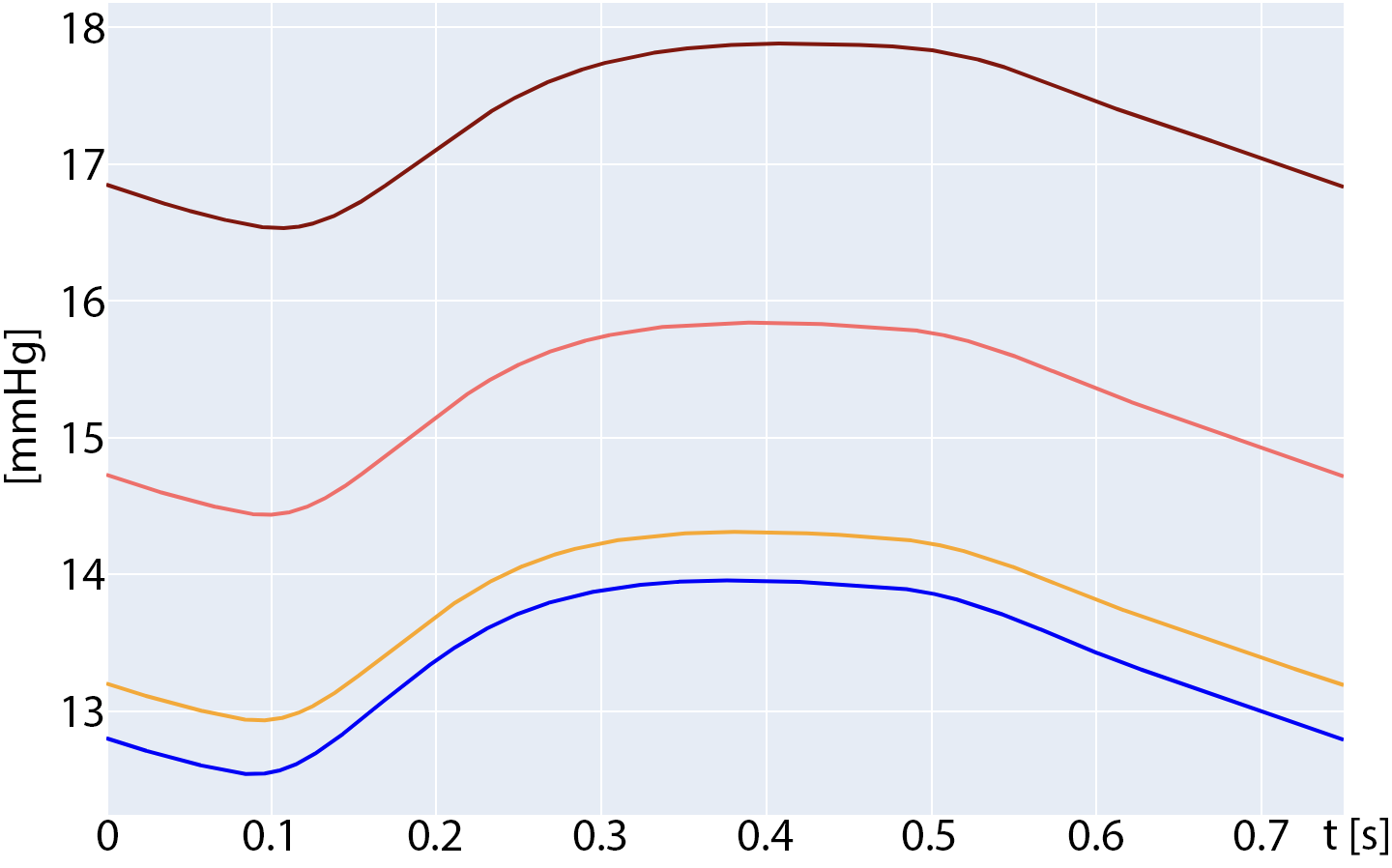}
    }
    \subfloat[$Q_\text{C}^\text{PUL}$.\label{fig:reno3d:qcpul}]{
        \includegraphics[width=0.31\linewidth]{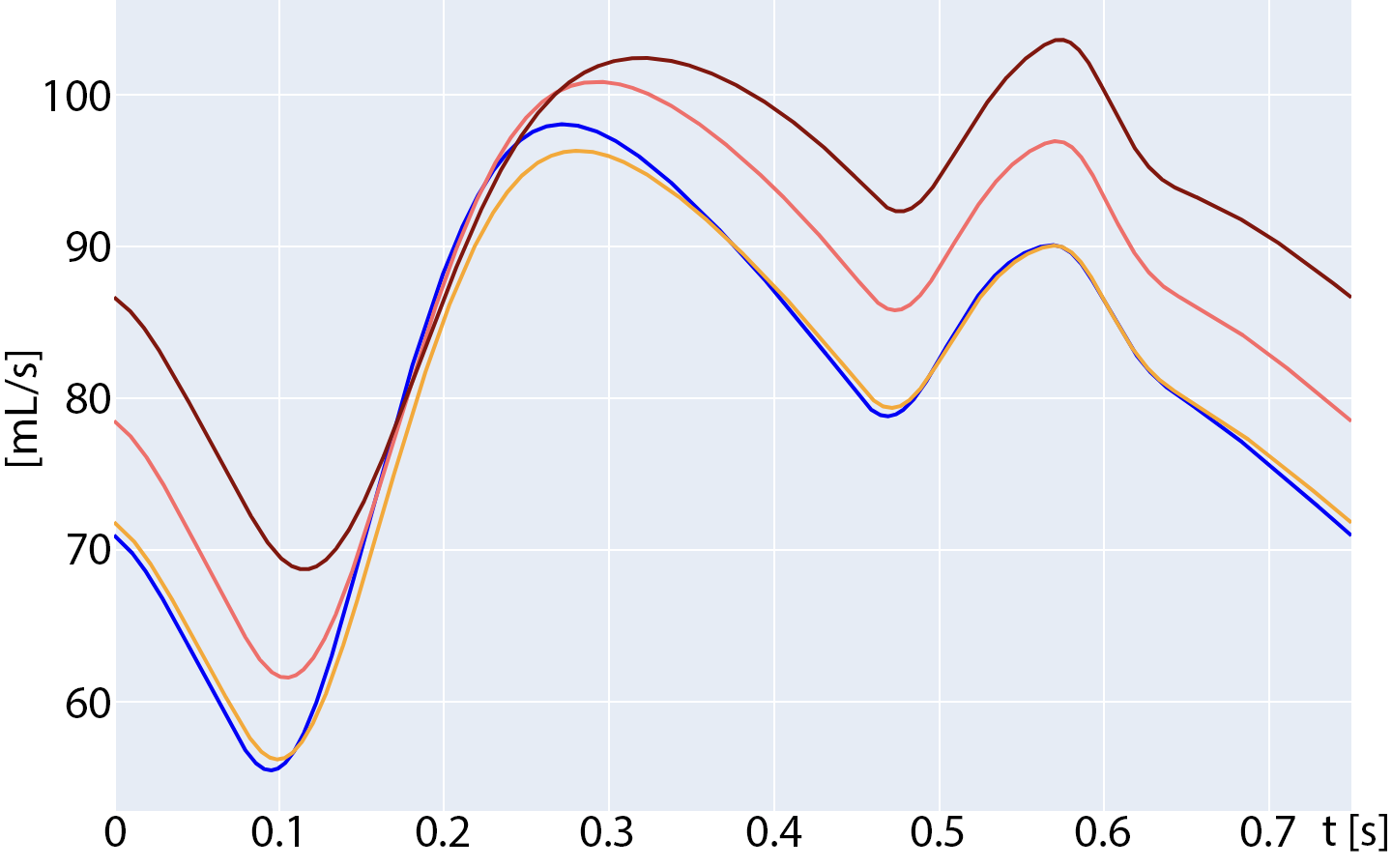}
    }\\\hspace{0.112\linewidth}
    \subfloat[$p_\text{C}^\text{SYS}$.\label{fig:reno0d:pcsys}]{
        \includegraphics[width=0.31\linewidth]{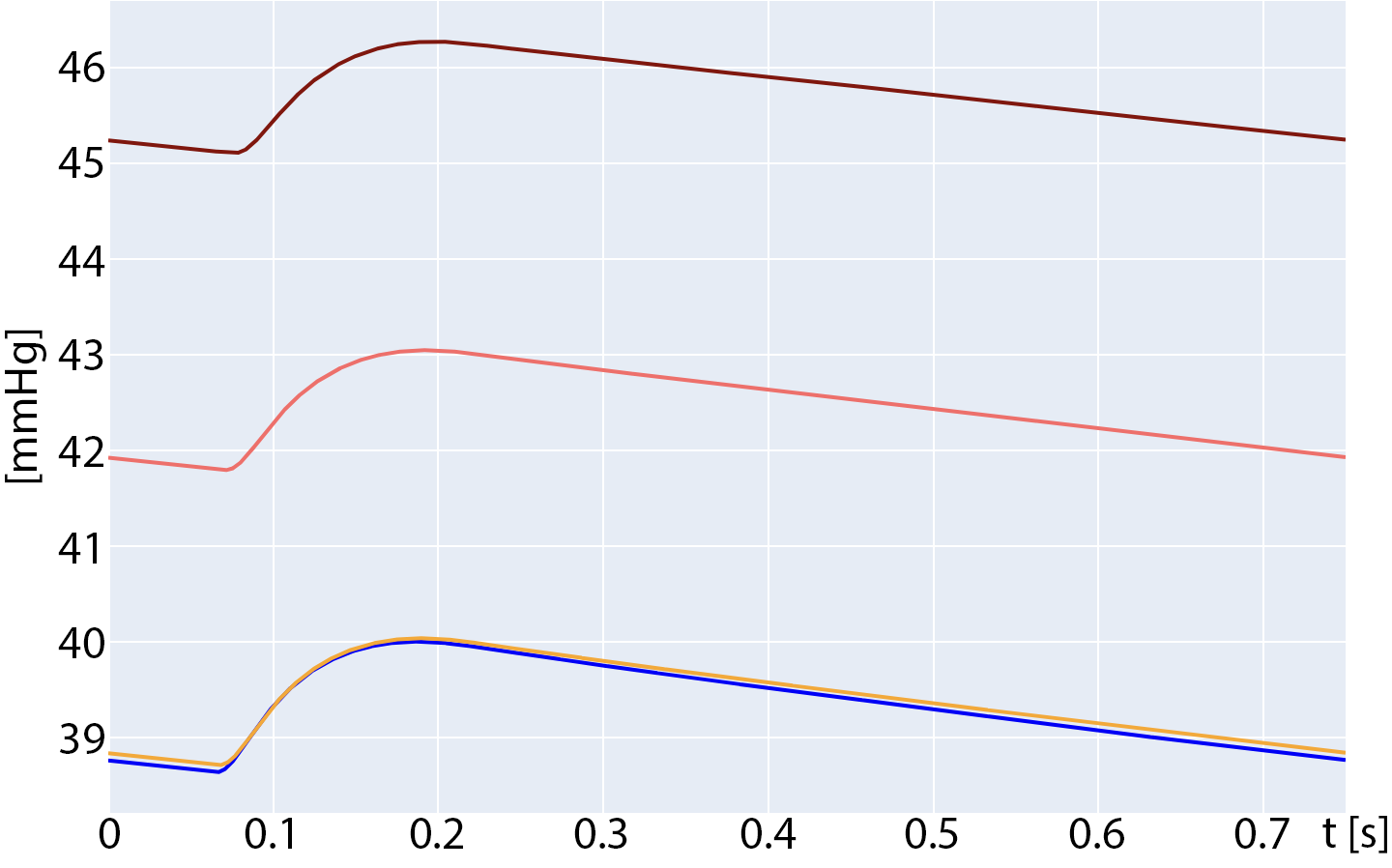}
    }
    \subfloat[$Q_\text{C}^\text{SYS}$.\label{fig:reno3d:qcsys}]{
        \includegraphics[width=0.31\linewidth]{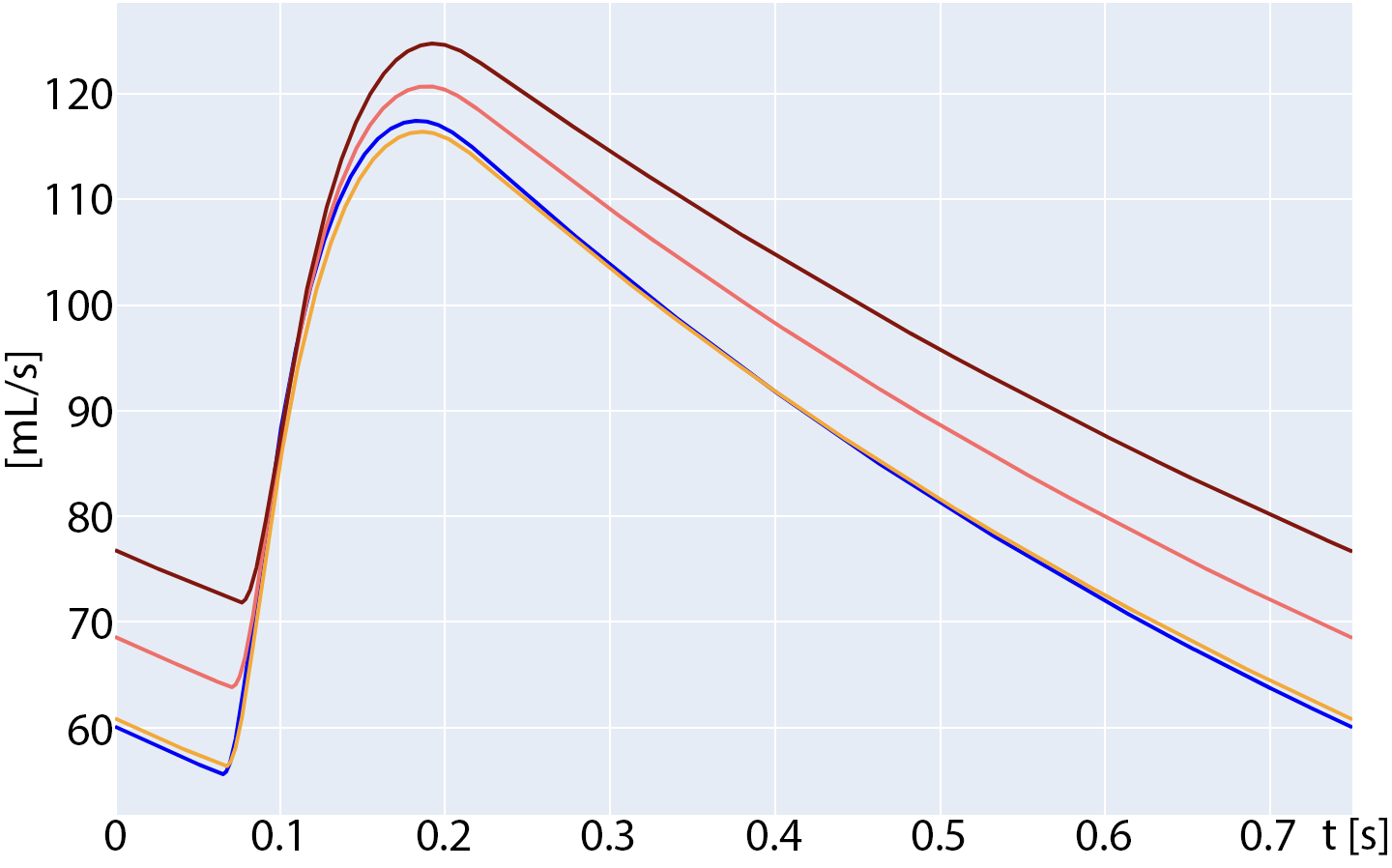}
    }
    \begin{subfigure}[b]{0.12\textwidth}
        \raisebox{0.6\height}{\includegraphics[width=\textwidth]{legend_variables.png}}
    \end{subfigure}
    \caption{Time-dependent capillary variables from $(\mathscr{C}_\text{C})$.
    Mild, moderate and severe 
    renovascular hypertension
    with secondary pulmonary hypertension 
    (in orange, light red and dark red) is compared with a 
    healthy individual (in blue).
    }
    \label{fig:reno0d-cap}
\end{figure}

\begin{figure}[t!]
    \centering
    \subfloat{
        \rotatebox{90}{\qquad $t=0.26\text{s}$}
        \includegraphics[width=0.225\linewidth]{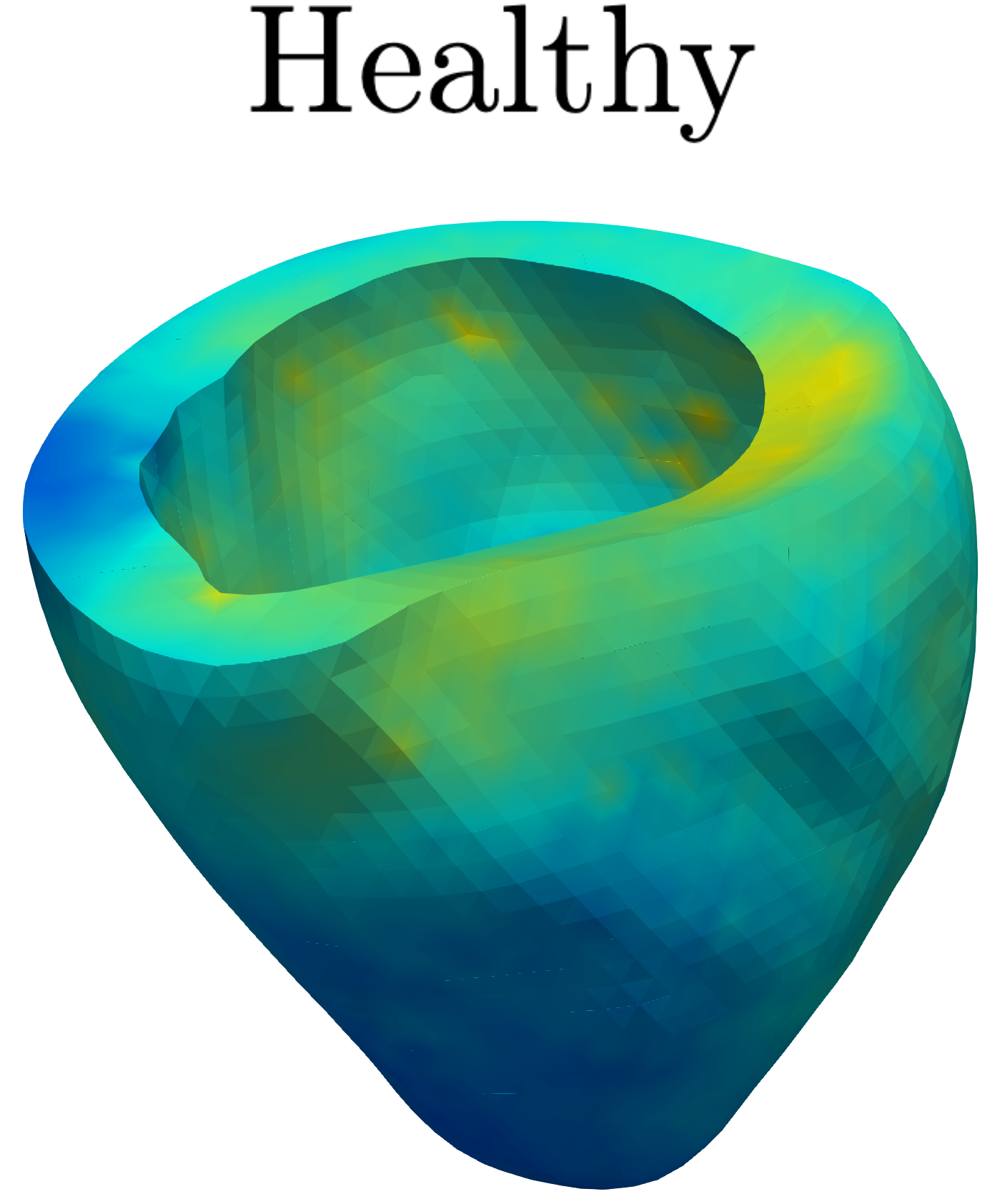}
    }
    \subfloat{
        \includegraphics[width=0.225\linewidth]{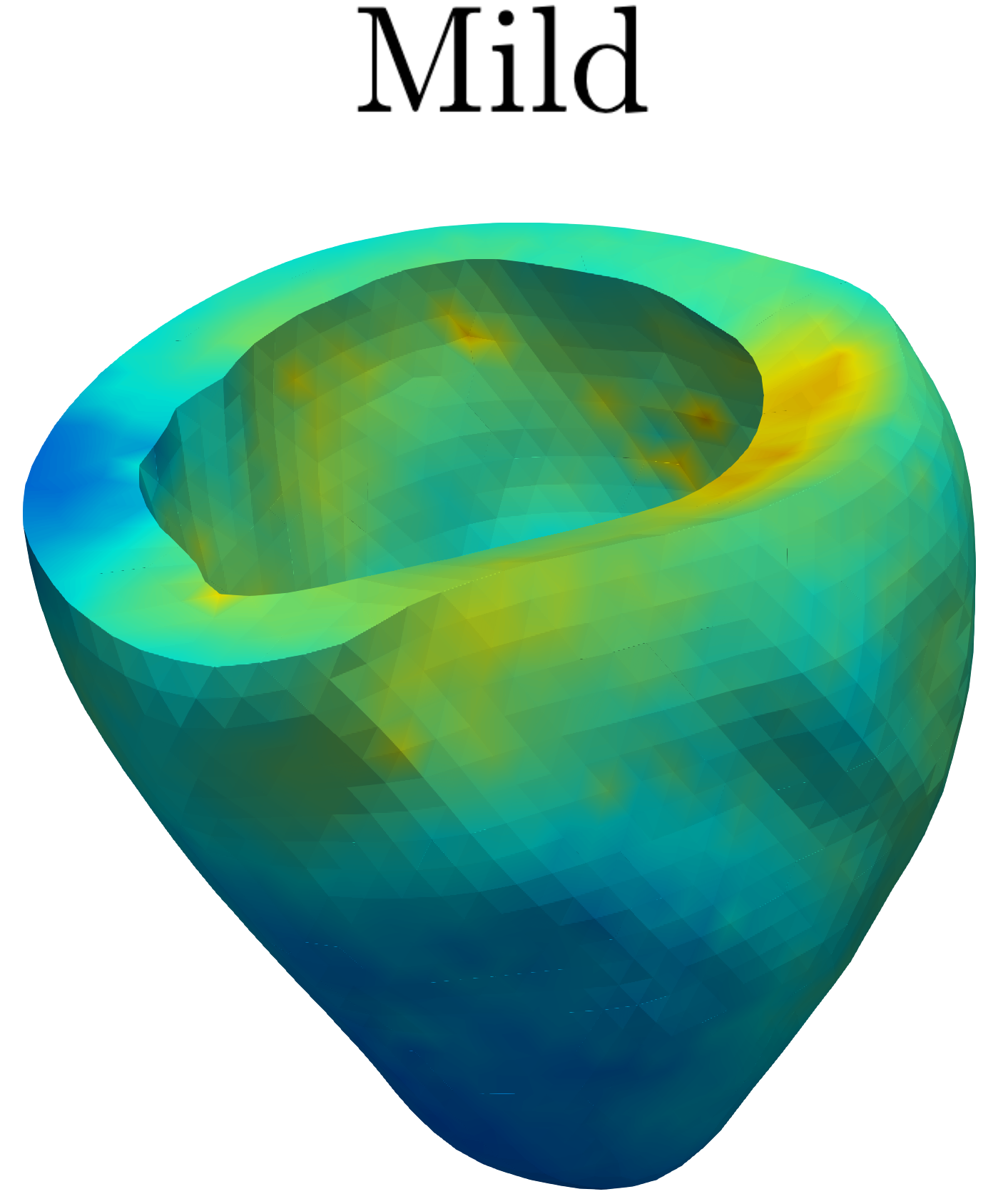}
    }
    \subfloat{
        \includegraphics[width=0.225\linewidth]{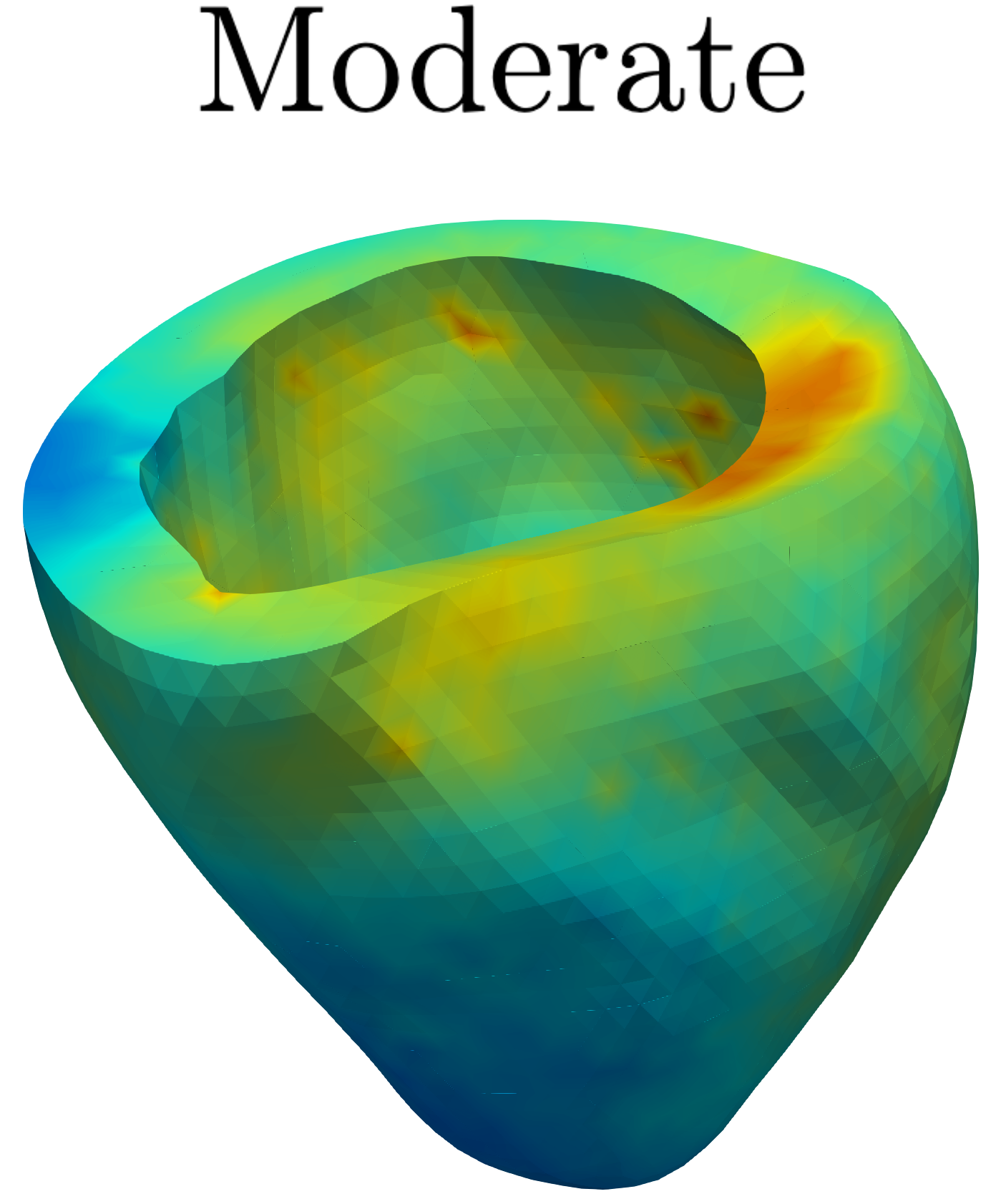}
    }
    \subfloat{
        \includegraphics[width=0.225\linewidth]{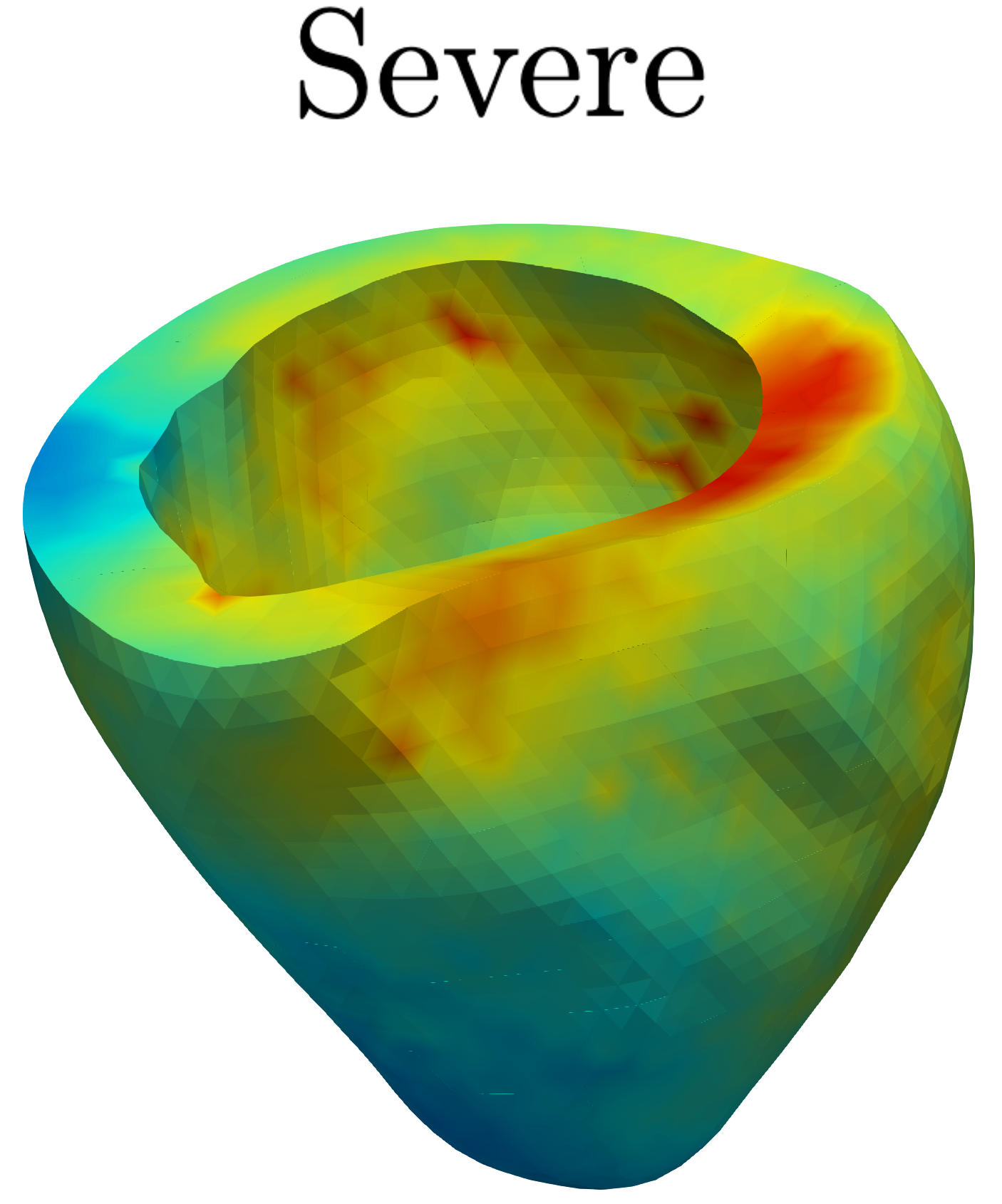}
    }\\
    \subfloat{
        \rotatebox{90}{\qquad $t=0.35\text{s}$}
        \includegraphics[width=0.225\linewidth]{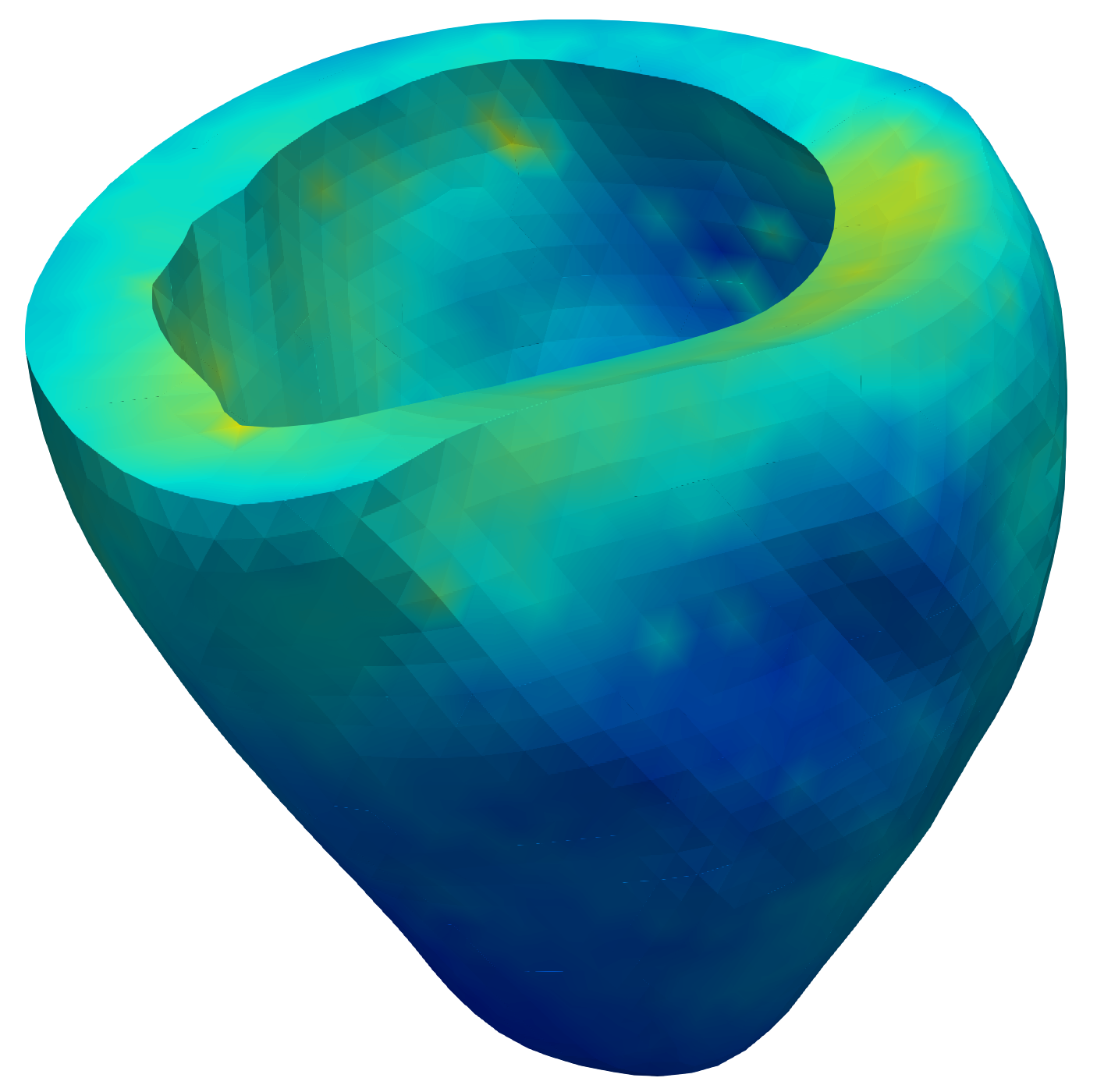}
    }
    \subfloat{
        \includegraphics[width=0.225\linewidth]{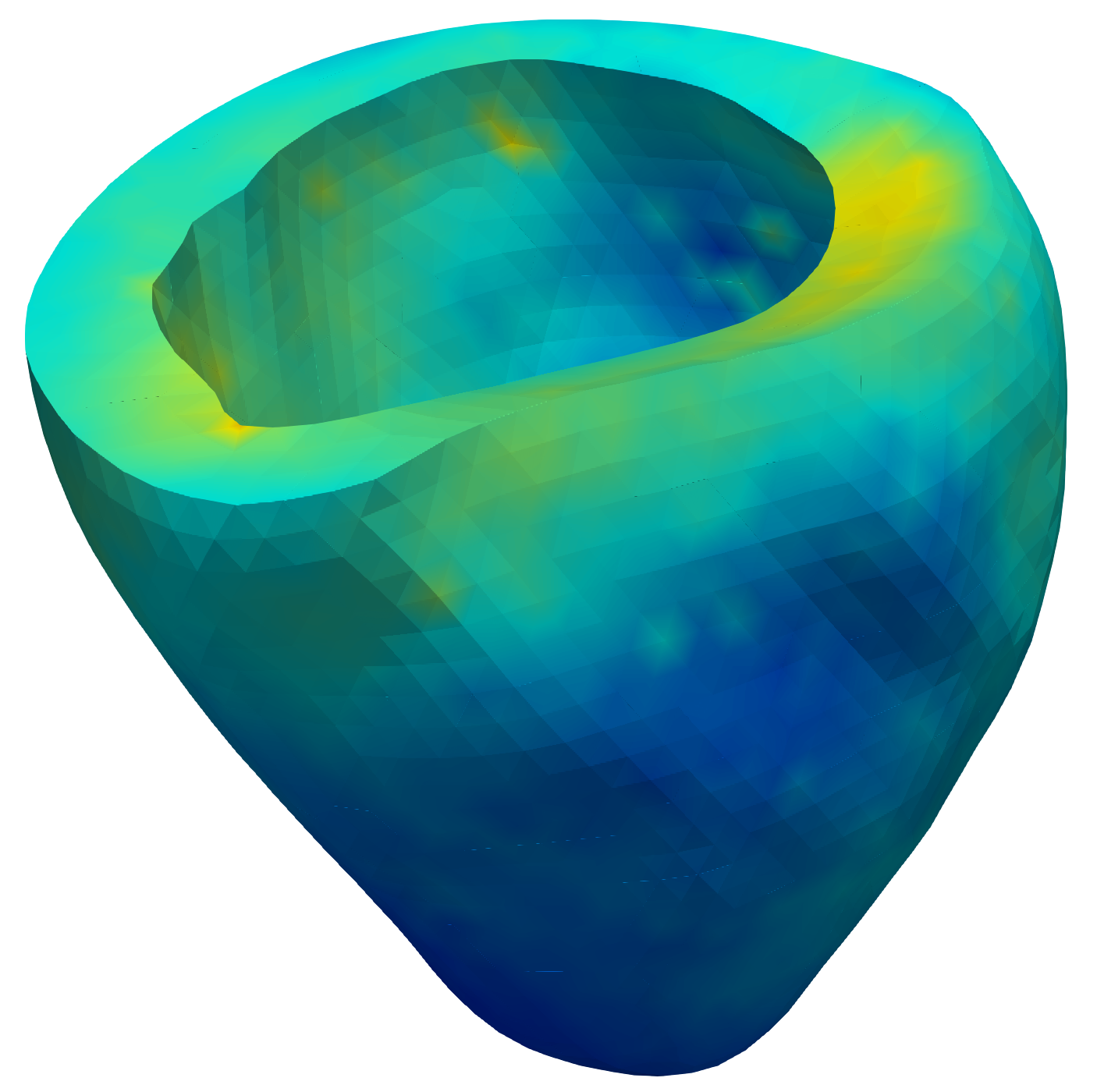}
    }
    \subfloat{
        \includegraphics[width=0.225\linewidth]{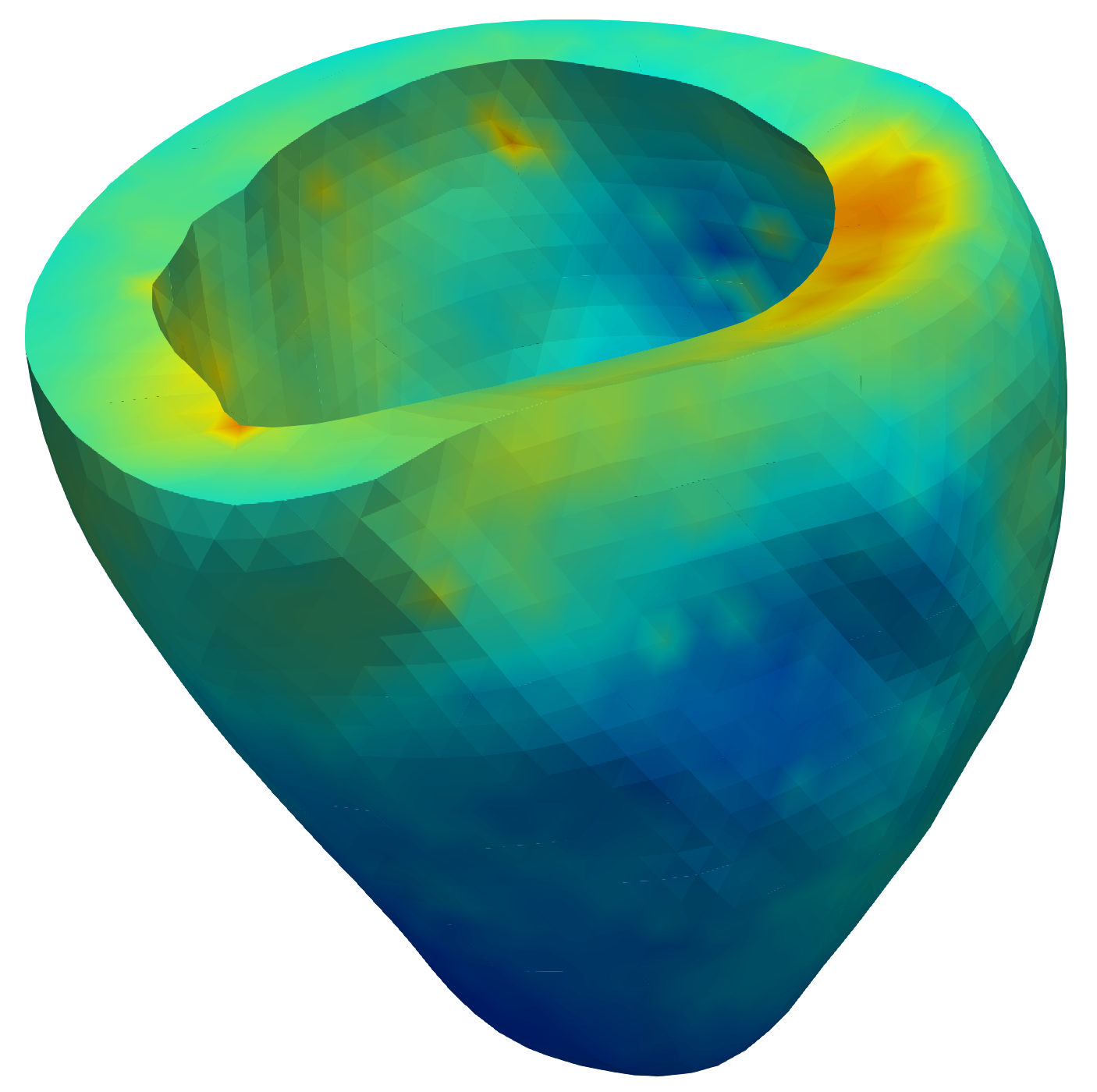}
    }
    \subfloat{
        \includegraphics[width=0.225\linewidth]{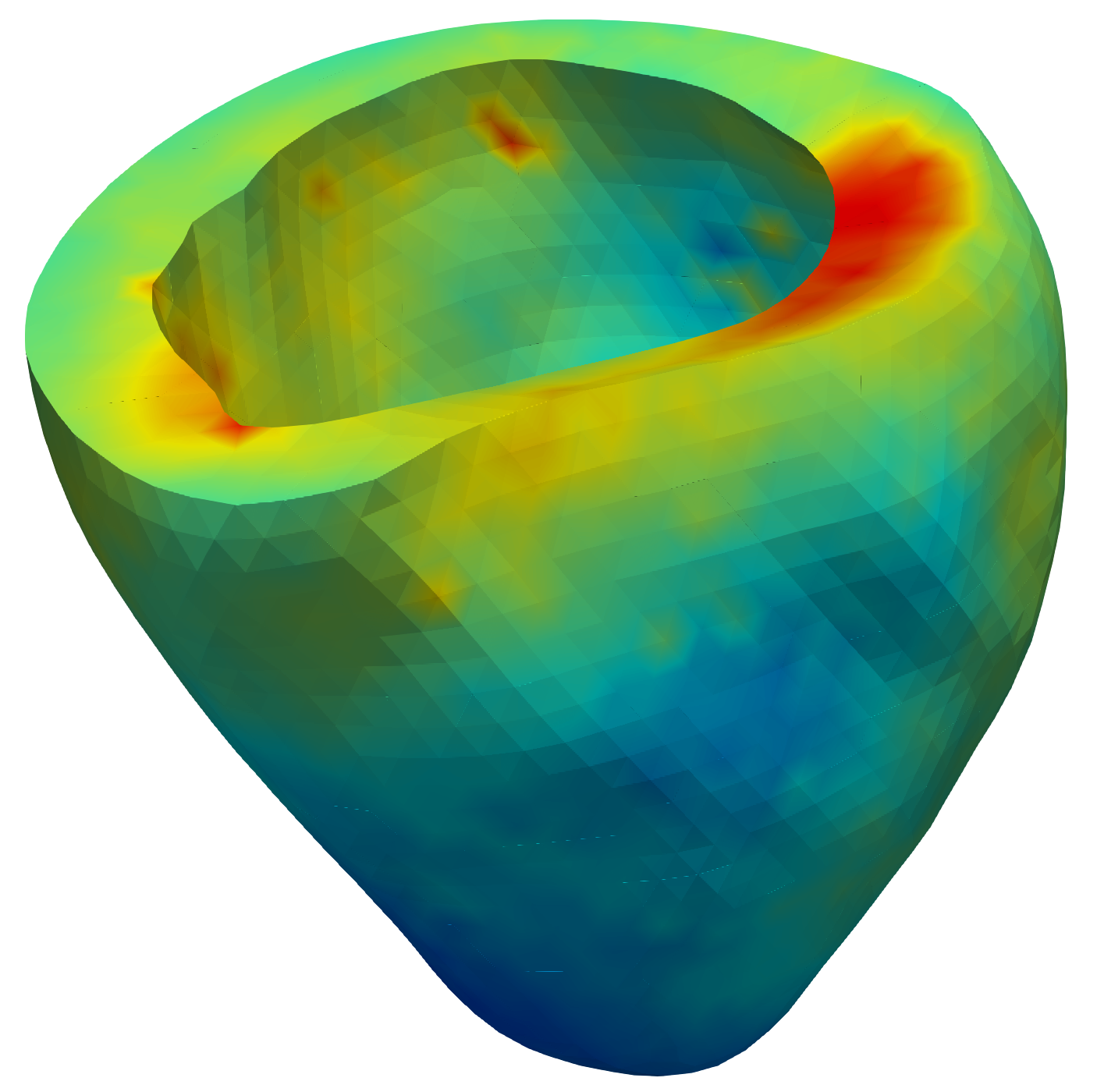}
    }\\
    \subfloat{
        \includegraphics[width=0.98\linewidth]{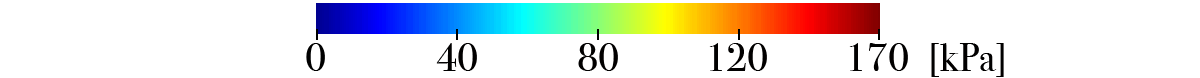}
    }
    \caption{Evolution of the active tension $T_a$ in the left ventricle 
    during the last cardiac cycle.
    Each picture is displayed on the reference configuration $\Omega_0$. 
    Mild, moderate and severe renovascular hypertension with 
    secondary pulmonary hypertension is compared with a healthy individual.
    }
    \label{fig:reno3d-ta}
\end{figure}

\begin{figure}[t!]
    \centering
    \subfloat{
        \rotatebox{90}{\qquad $t=0.15\text{s}$}
        \includegraphics[width=0.225\linewidth]{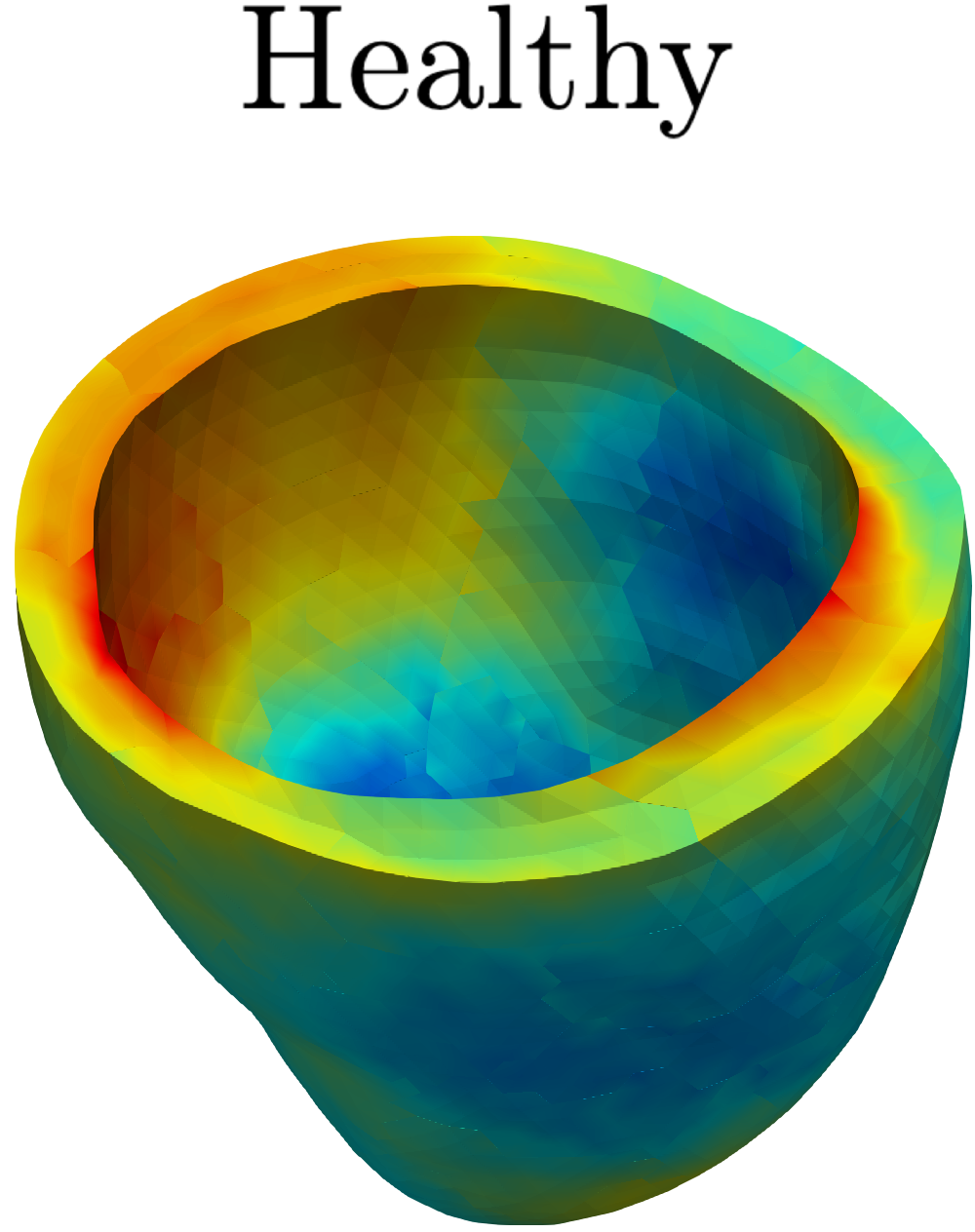}
    }
    \subfloat{
        \includegraphics[width=0.225\linewidth]{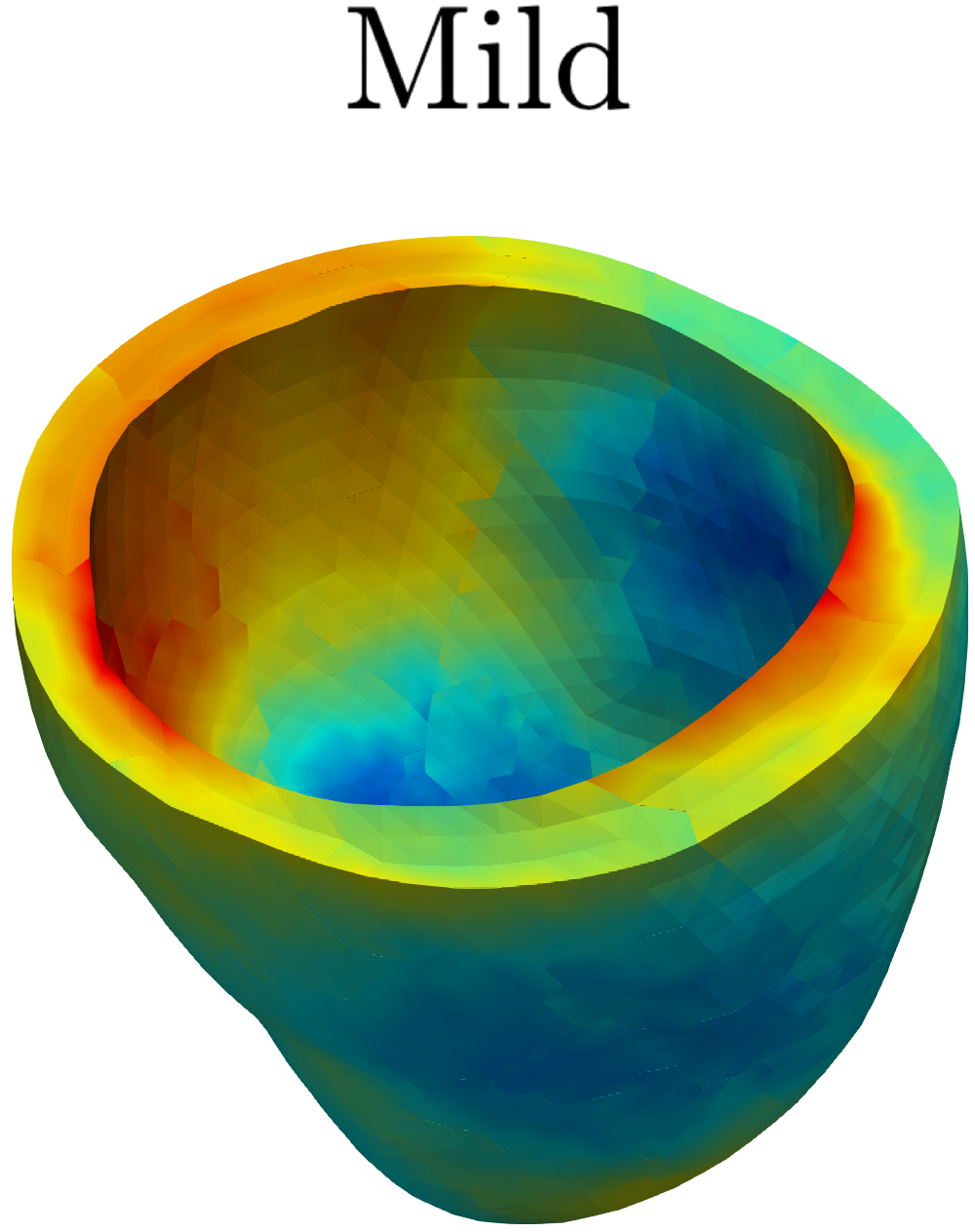}
    }
    \subfloat{
        \includegraphics[width=0.225\linewidth]{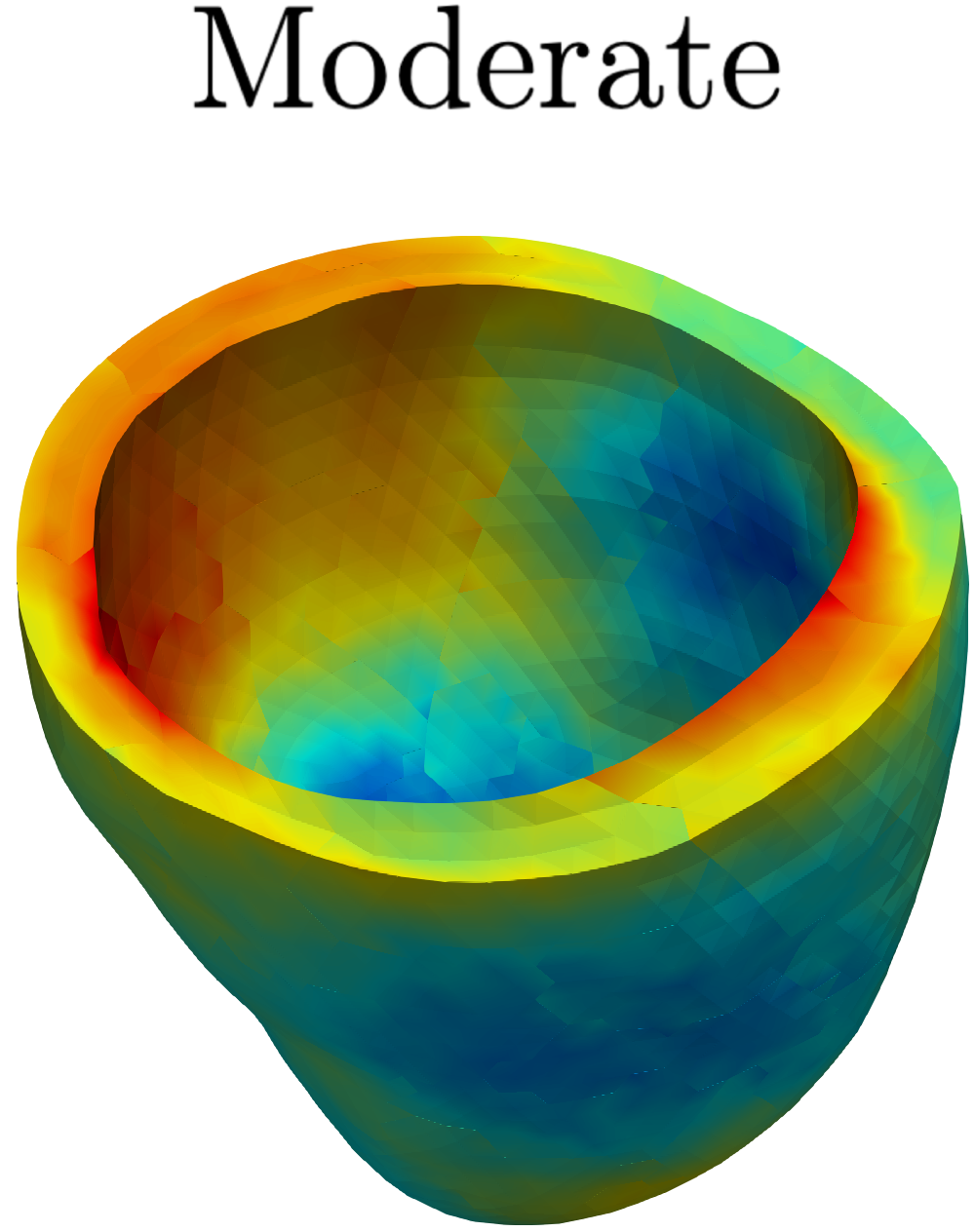}
    }
    \subfloat{
        \includegraphics[width=0.225\linewidth]{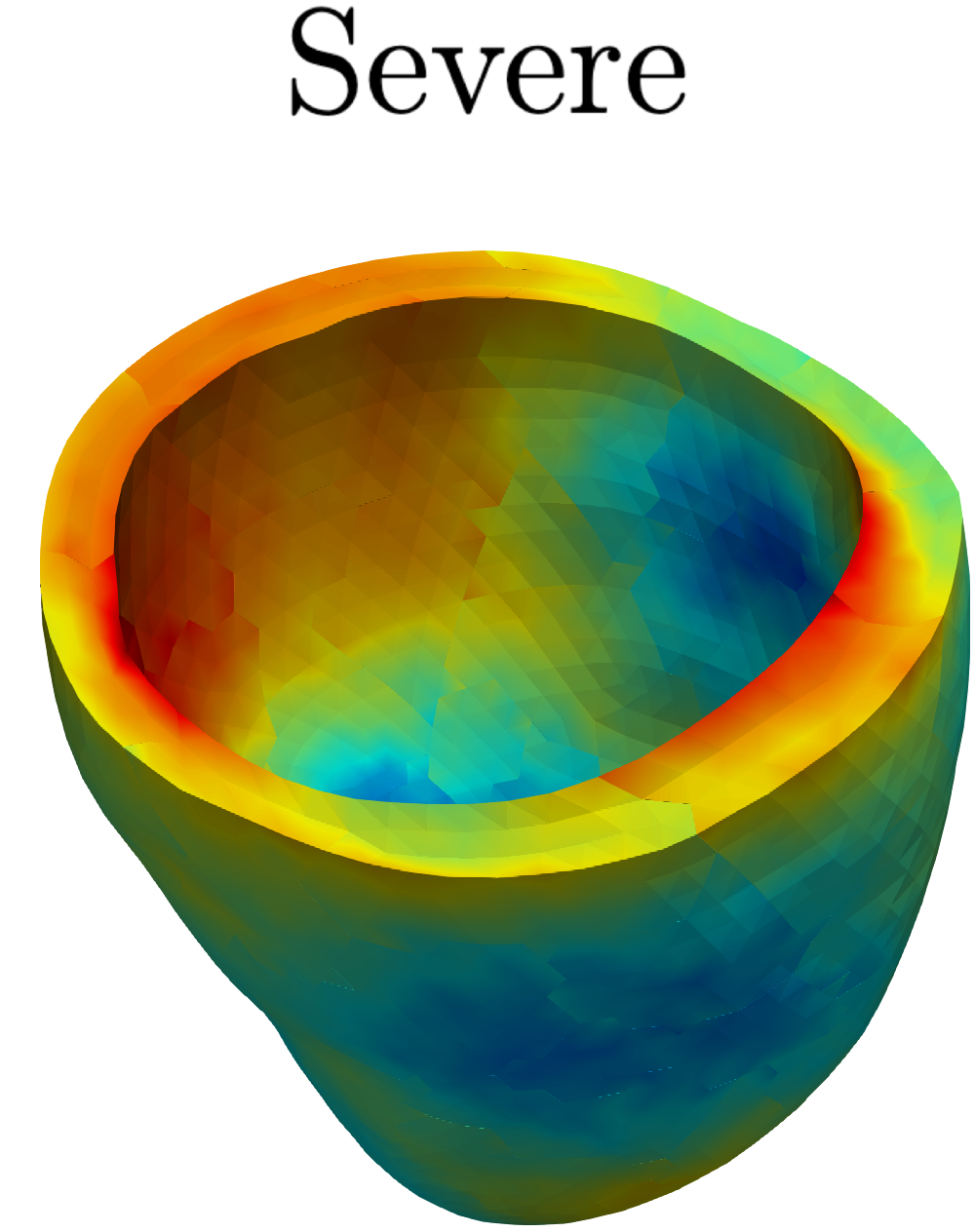}
    }\\
    \subfloat{
        \rotatebox{90}{\qquad $t=0.23\text{s}$}
        \includegraphics[width=0.225\linewidth]{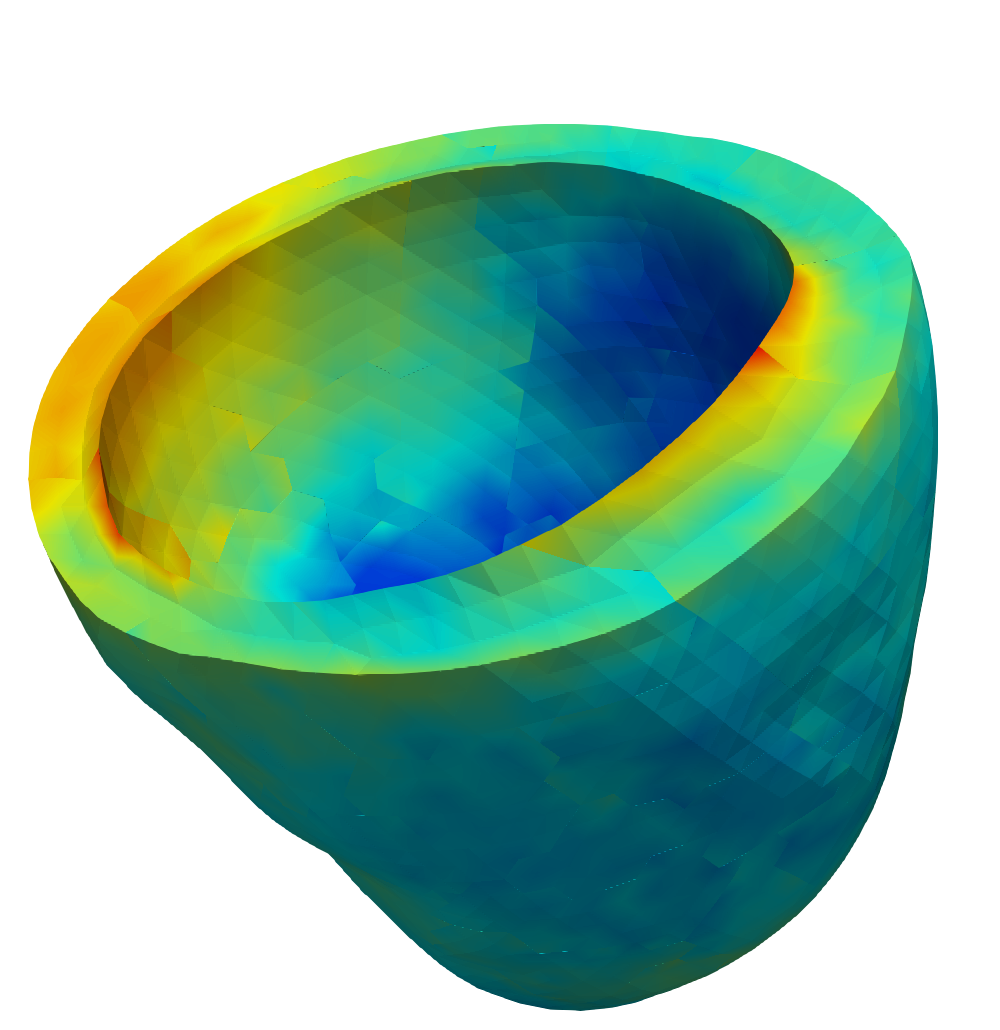}
    }
    \subfloat{
        \includegraphics[width=0.225\linewidth]{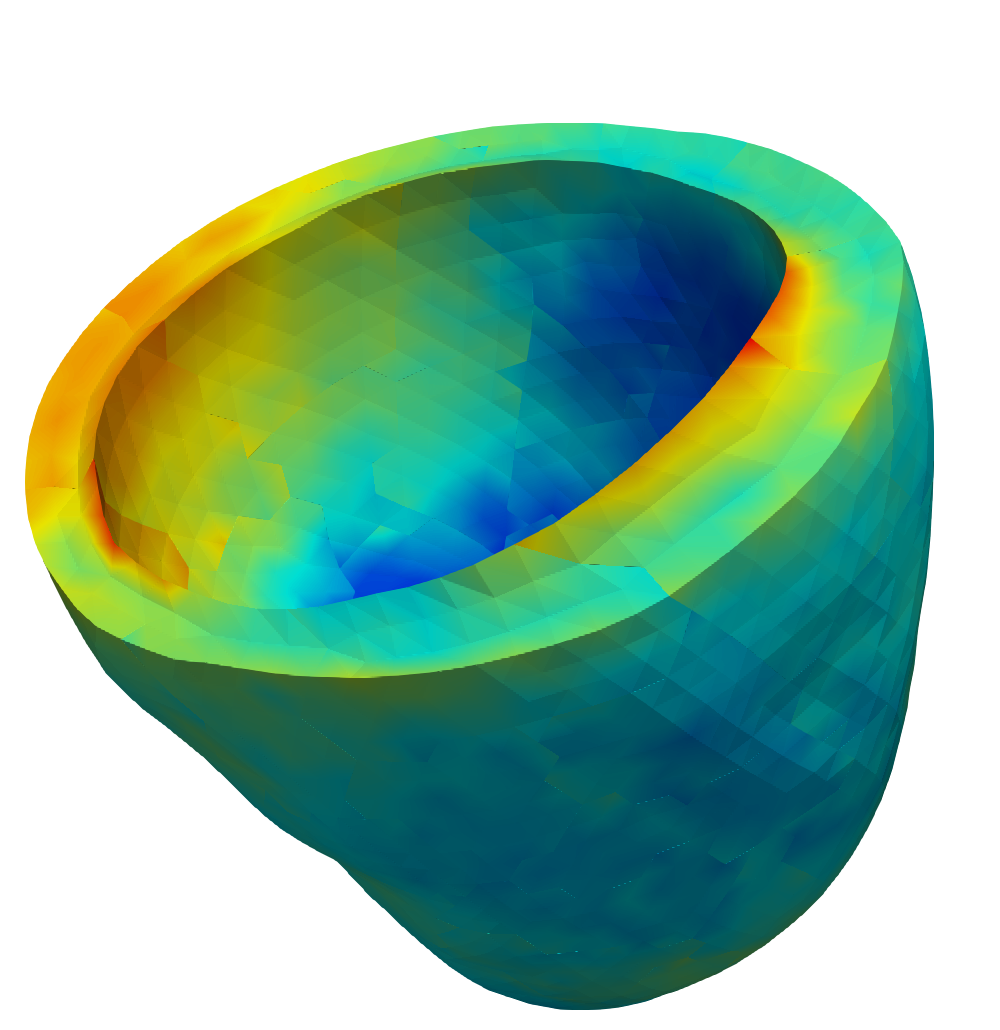}
    }
    \subfloat{
        \includegraphics[width=0.225\linewidth]{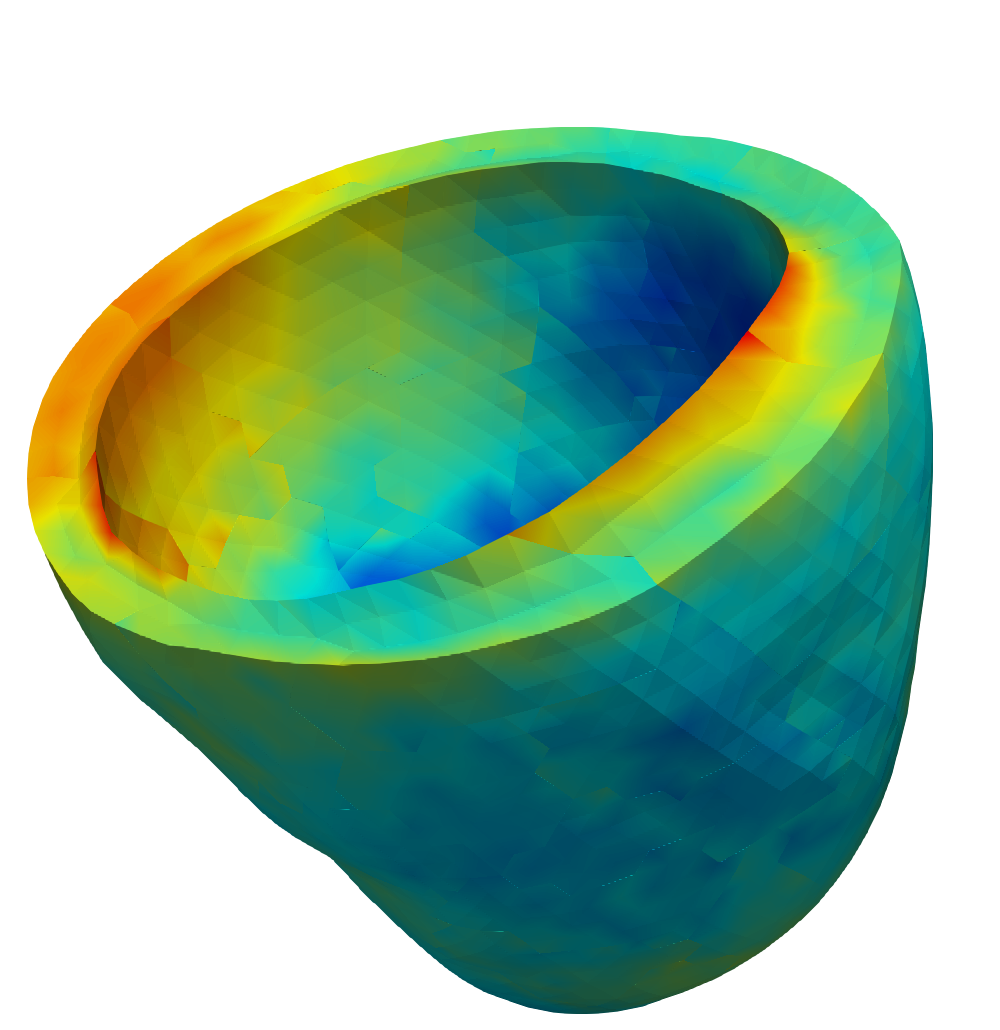}
    }
    \subfloat{
        \includegraphics[width=0.225\linewidth]{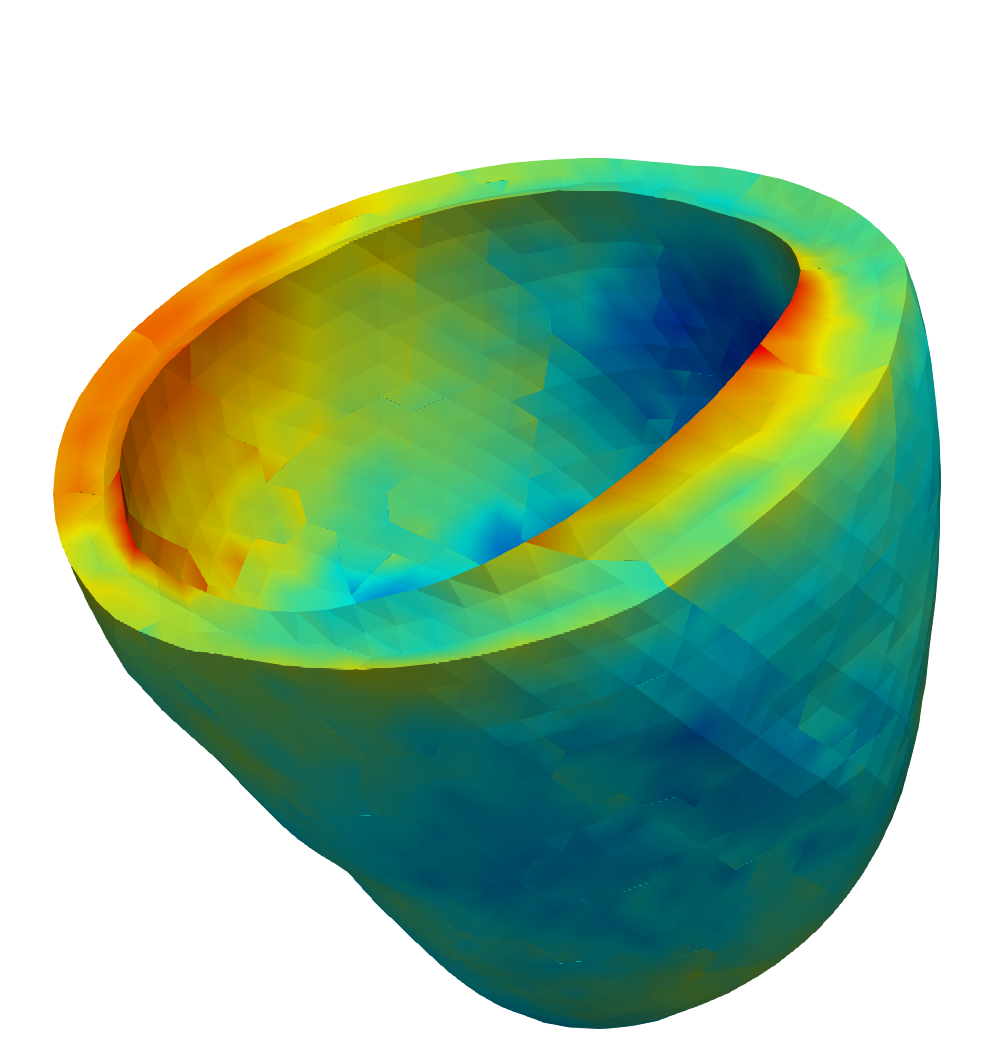}
    }\\
    \subfloat{
        \includegraphics[width=0.98\linewidth]{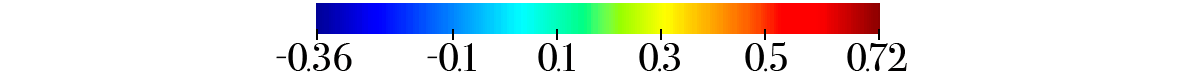}
    }
    \caption{Evolution of the Green-Lagrange strain trace $\operatorname{tr}(\mathbf{E})$ in 
    the left ventricle during the last cardiac cycle.
    Each picture is warped by the 
    displacement vector $\boldsymbol{d}$. Mild, moderate and severe renovascular hypertension with 
    secondary pulmonary hypertension is compared with a healthy individual.
    }
    \label{fig:reno3d-displ}
\end{figure}

In both models, the left ventricular PV loop 
(Figures~\ref{fig:reno0d:lv-loop} and \ref{fig:reno3d:lv-loop})
shows a clear increase in end-systolic pressure due to high systemic vascular 
resistance, leading to an increase in $\text{LV}_\text{SW}$
(Table~\ref{table:reno-out}); 
however, in the 3D--0D model, systolic work is slightly lower in severe hypertension, 
suggesting a more physiologically realistic redistribution of workload, likely due to 
the spatial heterogeneity of myocardial contraction captured by the 3D mechanics. 
As presented in Table~\ref{table:reno-out}, 
$\text{LV}_\text{IEDV}$ increases progressively in both models, but the 3D--0D 
predictions are slightly higher than $(\mathscr{C}_\text{C})$, indicating more 
pronounced ventricular dilatation; $\text{LV}_\text{IESV}$ also 
increases in both models, but the 3D--0D increase is more moderate, reflecting 
relatively preserved systolic emptying compared to the 0D results.

The right ventricle responds to secondary pulmonary hypertension with elevated 
end-diastolic and end-systolic volumes in both models, as shown by 
$\text{RV}_\text{IEDV}$ and $\text{RV}_\text{IESV}$ 
in Table~\ref{table:reno-out}, and by the right ventricular PV loops 
in Figures~\ref{fig:reno0d:rv-loop} and \ref{fig:reno3d:rv-loop}. 
The 3D--0D model generally predicts slightly higher values, in line with increased 
filling in the presence of higher afterload. $\text{RV}_\text{SW}$ and 
$\text{RV}_\text{IESV}$ increase in both 
models, with the 3D--0D simulation showing slightly higher values, reflecting the 
greater energy demand and ventricular filling under increased pulmonary resistance.

Both $p_\text{AR}^\text{SYS}$
(Figures~\ref{fig:reno0d:parsys} and \ref{fig:reno3d:parsys})
and $Q_\text{AR}^\text{SYS}$ 
(Figures~\ref{fig:reno0d:qarsys} and \ref{fig:reno3d:qarsys})
increase in both models, as the circulation adapts to preserve systemic perfusion. 
In $(\mathscr{C}_\text{C})$, these increases appear more uniform throughout the entire 
cardiac cycle, while in the 3D--0D model they are more concentrated around the systolic 
peak: this difference originates from the more detailed spatial representation of the 
left ventricle in the 3D--0D framework, which captures the non-linear interplay between 
ventricular contraction, pressure wave propagation and vascular impedance. As a result, 
the 3D--0D model provides a sharper and more realistic distribution of 
pressures and flows, whereas $(\mathscr{C}_\text{C})$, due to its lumped nature, tends 
to smooth out temporal variations and therefore depicts a more averaged behaviour.
In Figures~\ref{fig:reno0d:parpul} and \ref{fig:reno3d:parpul}, $p_\text{AR}^\text{PUL}$ 
also increases due to secondary pulmonary hypertension, while the 3D--0D model 
highlights how the right ventricle buffers left ventricular preload, maintaining higher 
volumes and systolic work and capturing ventriculo-vascular interactions more accurately 
than the 0D simulation. In Figure~\ref{fig:reno0d-cap}, raised $p_\text{C}^\text{PUL}$ 
and $Q_\text{C}^\text{PUL}$ suggest altered lung microcirculation in the 0D 
representation, risking fluid buildup and impaired gas exchange, while increased 
$p_\text{C}^\text{SYS}$ and $Q_\text{C}^\text{SYS}$ point to systemic microvascular 
dysfunction, leading to tissue oedema and impaired nutrient delivery \cite{capillary-reno}.

In Figure~\ref{fig:reno3d-ta}, active tension markedly increases in maximal mean 
values, enhancing pumping during ejection; maximal and minimal values rise over the 
entire left ventricle, with the upper ventricle showing the greatest increase at 
$t=0.26\,\text{s}$ and $t=0.35\,\text{s}$. 
The trace of $\mathbf{E}$ in Figure~\ref{fig:reno3d-displ} exhibits more pronounced 
positive and negative peaks, reflecting greater volumetric deformation of the 
left ventricle, with stronger dilation and more intense contraction during the cardiac 
cycle. These larger variations are consistent with the altered loading conditions 
imposed by renovascular hypertension, indicating increased mechanical stress on the 
ventricular myocardium.

\section{Conclusions}
\label{conclusions}%
This work developed and applied a computational modeling framework to study hypertension in 
systemic, pulmonary and renovascular forms and in three nuances of severity, using both 
0D models (with and without capillaries) \cite{lumped:capillary} and a 3D--0D coupled model with 
electromechanical representation of the left ventricle \cite{3d-0d_ventr},
within a perspective that moved toward the development of cardiovascular DTs.
The methodology included 
parameter calibration, global and multi-parametric sensitivity analysis and systematic 
comparison between models with and without capillaries.

The calibration approach allowed for a clear distinction between degrees of hypertension, 
reproducing physiological trends: increase in systolic and diastolic blood pressure, 
alterations in valve flows, changes in cardiac volumes and systemic venous and 
pulmonary pressures.
The 0D model with capillaries $(\mathscr{C}_\text{C})$ demonstrated a high degree of 
flexibility in representing a broad spectrum of haemodynamic conditions: it was able to 
incorporate variations in cardiovascular parameters, consistent with the typical 
alterations observed in clinical cases of hypertension. The comparison with the version 
without capillaries $(\mathscr{C}_\text{NC})$ showed that, by adjusting the parameters 
appropriately, comparable results can be obtained, confirming the robustness of the 
simplified model.
The comparison between the 0D model and the 3D--0D model showed substantial differences 
in the ability to reproduce haemodynamic phenomena, especially under pathological 
conditions. The 3D--0D model, thanks to the three-dimensional mechanical description of 
the left ventricle, returned more realistic pressure-volume curves, consistent with 
literature clinical observations. It showed a more physiological behaviour of the cardiac index in 
systemic hypertension, greater sensitivity in the representation of right ventricular 
dysfunction in pulmonary hypertension, and better consistency of venous and arterial 
flows and pressures in the renovascular form. Furthermore, the spatial distribution of 
stresses and strains in the ventricle made it possible to pick up early signs of 
pathological remodeling that the more simplified 0D model tends to overlook.

Overall, the model responded positively to its objective: to represent the 
hypertensive condition in a realistic and robust way 
and to provide a foundation for patient-specific DT frameworks capable of integrating 
both reduced-order and high-fidelity electromechanical descriptions.

We recognize that our study has certain limitations. 
The calibration was conducted on synthetic data generated by 
$(\mathscr{C}_\text{C})$ model: it ensures a controlled comparison, but introduces the 
risk of not fully reflecting the individual variability of real patients.
The lack of direct calibration of the 3D--0D model means that some mechanical 
aspects of the ventricle (particularly electrophysiological parameters and active 
contractility parameters) remain unoptimised 
and potentially less realistic, as certain physiological mechanisms regulating 
contraction may not be fully captured.
3D simulation requires high computational effort, limiting the number of 
scenarios that can be explored and its applicability in real-time clinical settings;
this can be circumvented by means of reduced order models \cite{ldnets}.

A significant extension of the work could arise from the integration of real 
clinical data to calibrate and validate the models. Using specific measurements, such 
as pressure, flow, volume, cardiac geometries, would allow the creation of true  
cardiovascular DTs, potentially useful for personalised diagnostics and monitoring. 
Such an approach requires more robust strategies to handle inter-individual variability 
and implement adaptive cost functions.
Possible avenues of extension include 3D modeling of vascular and cardiac structures 
relevant in hypertension (e.g. \cite{model_hyper:es,model-whole-heart1, model-two-chambers}):
\begin{enumerate}
    \item Systemic hypertension: 3D modeling of the aorta, carotids and femoral 
    arteries to study pressure wave propagation.
    \item Pulmonary hypertension: include right atrium, right ventricle and pulmonary 
    vessels to better validate heart-lung interaction.
    \item Renovascular hypertension with secondary pulmonary hypertension: add a 3D 
    representation of the renal circulation, with dedicated arteries and veins, and 
    the inclusion of equations describing the effect of sodium concentration on 
    pressure control. The kidney contributes to hypertension both through diuretic and 
    natriuretic mechanisms \cite{kidney7,kidney6} and by modulating sympathetic tone through reflexes that 
    increase nerve activity and blood pressure \cite{kidney8}: these aspects are particularly 
    important in salt-sensitive cases and in chronic conditions, particularly in 
    renovascular hypertension.
\end{enumerate}

Finally, the possibility of directly calibrating the 3D--0D model, including 
electromechanical parameters, shall be considered. Although computationally demanding, 
this extension would represent a crucial step toward DT-ready electromechanical models, 
bridging the gap between high-fidelity simulations and real-time personalised applications.

\section*{Acknowledgements}
The authors acknowledge Dr. Andrea Parlangeli for having facilitated the collaboration between Professors Parati and Quarteroni on this interesting scientific topic.

AT, FR, LD and AQ are members of the INdAM group GNCS 
\say{Gruppo Nazionale per il Calcolo Scientifico} 
(National Group for Scientific Computing).

AT, FR and LD acknowledge the INdAM GNCS project CUP E53CE53C24001950001.
FR has received support from the project FIS, MUR, Italy 2025-2028, Project code: 
FIS-2023-02228, CUP: D53C24005440001, 
\say{SYNERGIZE: Synergizing Numerical Methods and Machine Learning for a new generation 
of computational models}.
LD acknowledges the project PRIN2022, MUR, Italy, 2023-2025, 202232A8AN 
\say{Computational modeling of the heart: from eﬀicient numerical solvers to cardiac 
digital twins}.

The present research is part of the activities of 
“Dipartimento di Eccellenza 2023-2027”, MUR, Italy, Dipartimento di Matematica, 
Politecnico di Milano.

\newpage

\appendix
\section{3D--0D Model Equations}
\label{appendix:3d-0d_eq}
The complete 3D--0D model is defined by the set of Equations~\eqref{3d0d-eq}: 
Equations~\eqref{3d0d-eq:electrical}--\eqref{3d0d-eq:mechanic} are solved in the space-time 
domain $\Omega_0\times(0,\text{T}]$, while Equation~\eqref{3d0d-eq:circ} 
in $(0,\text{T}]$, and Equation~\eqref{3d0d-eq:coupling} represents the coupling 
condition between 3D and 0D. 
\begin{subequations}\label{3d0d-eq}
    \begin{align}
    &(\mathscr{E})
    \begin{cases}\label{3d0d-eq:electrical}
        J\chi_\text{m}\left[\text{C}_\text{m}\dfrac{\partial u}{\partial t}+\mathcal{I}_{\text{ion}}(u,\boldsymbol{w},\boldsymbol{z})\right]-&\!\!\!\!\!\nabla\cdot\left(J\mathbf{F}^{-1}\mathbf{D}_\text{M}\mathbf{F}^{-T}\nabla u\right)=J\chi_\text{m}\,\mathcal{I}_\text{app}(t), \\
        \left(J\mathbf{F}^{-1}\mathbf{D}_\text{M}\mathbf{F}^{-T}\nabla u\right)\cdot \mathbf{N}=0 &\text{on }\partial\Omega_0\times(0,\text{T}],\\
        u=u_0 &\text{in }\Omega_0\times \left\{0\right\},
    \end{cases}\\
    &(\mathscr{I})
    \begin{cases}\label{3d0d-eq:ions}
        \dfrac{\partial\boldsymbol{w}}{\partial t}-\mathbf{H}(u,\boldsymbol{w})=\boldsymbol{0},\\[7pt]
        \dfrac{\partial\boldsymbol{z}}{\partial t}-\mathbf{G}(u,\boldsymbol{w},\boldsymbol{z})=\boldsymbol{0},\\
        \boldsymbol{w}=\boldsymbol{w}_0, \,\,\,\, \boldsymbol{z}=\boldsymbol{z}_0 &\quad\text{in }\Omega_0\times \left\{0\right\},
    \end{cases}\\
    &(\mathscr{A})
    \begin{cases}\label{3d0d-eq:active}
        \dfrac{\partial\boldsymbol{s}}{\partial t}=\mathbf{F}_\text{act}\left(\boldsymbol{s},\left[\text{Ca}^{2+}\right]_i,\text{SL}, \dfrac{\partial \text{SL}}{\partial t}\right),\\
        \boldsymbol{s}=\boldsymbol{s}_0 &\text{in }\Omega_0\times \left\{0\right\},
    \end{cases}\\
    &(\mathscr{M})
    \begin{cases}\label{3d0d-eq:mechanic}
        \rho_s\dfrac{\partial^2\boldsymbol{d}}{\partial t^2}-\nabla\cdot\mathbf{P}\left(\boldsymbol{d},T_a(\boldsymbol{s})\right)=\boldsymbol{0},\\
        \mathbf{P}\left(\boldsymbol{d},T_a(\boldsymbol{s})\right)\mathbf{N}+\mathbf{K}^\text{epi}\boldsymbol{d}+\mathbf{C}^\text{epi}\dfrac{\partial \boldsymbol{d}}{\partial t}=\boldsymbol{0} &\text{on }\Gamma_0^\text{epi}\times(0,\text{T}],\\[4pt]
        \mathbf{P}\left(\boldsymbol{d},T_a(\boldsymbol{s})\right)\mathbf{N}=p_\text{LV}(t)\left|J\mathbf{F}^{-T}\mathbf{N}\right|\mathbf{v}^\text{base} &\text{on }\Gamma_0^\text{base}\times(0,\text{T}],\\
        \mathbf{P}\left(\boldsymbol{d},T_a(\boldsymbol{s})\right)\mathbf{N}=-p_\text{LV}(t)J\mathbf{F}^{-T}\mathbf{N} &\text{on }\Gamma_0^\text{endo}\times(0,\text{T}],\\[3pt]
        \boldsymbol{d}=\boldsymbol{d}_0, \quad\dfrac{\partial\boldsymbol{d}}{\partial t}=\dot{\boldsymbol{d}}_0 &\text{in }\Omega_0\times \left\{0\right\},
    \end{cases}\\
    &(\mathscr{C})\,\,\,
    \begin{cases}\label{3d0d-eq:circ}
        \dfrac{d\boldsymbol{c}_1(t)}{d t}=\overset{\sim}{\mathbf{D}}\bigl(t,\boldsymbol{c}_1(t),p_\text{LV}(t)\bigr) &\text{in }(0,\text{T}],\\
        \boldsymbol{c}_1(0)=\boldsymbol{c}_{1,0},
    \end{cases}\\
    &(\mathscr{V})\label{3d0d-eq:coupling}
    \,\,\,\,\,V_\text{LV}^\text{0D}\left(\boldsymbol{c}_1(t)\right)=V_\text{LV}^\text{3D}\left(\boldsymbol{d}(t)\right)\qquad\quad\text{in }(0,\text{T}].
\end{align}
\end{subequations}

\section{Healthy individual parameters}
\label{appendix:param-healthy}%
\begin{table}[H]
    \centering 
    \begin{tabular}{lccccc}
    \hline
    Parameter & $i=\text{LA}$ & $i=\text{LV}$ & $i=\text{RA}$ & $i=\text{RV}$ &Unit\\
    \hline
    $\text{E}^a_i$          & 0.255     & 8.442     &$6\cdot 10^{-2}$       & 0.495 & $\text{mmHg}/\text{mL}$\\
    $\text{E}^p_i$          & 0.1512    &0.126      &$7\cdot 10^{-2}$       &$7\cdot 10^{-2}$ & $\text{mmHg}/\text{mL}$ \\
    $\text{T}_\text{C}^i$   & 0.15      &0.25       &0.1        &0.25 & s \\
    $\text{T}_\text{R}^i$   & 0.8       &0.5        & 0.7       & 0.4 & s \\
    $\text{t}_\text{C}^i$   & 0.75      &0          &0.8        &0 & s \\
    $\text{V}_{0,i}$        & 4         &42         &4          &16& mL \\
    $\text{R}_{\text{min}}$ &$7.5\cdot 10^{-3}$ &$7.5\cdot 10^{-3}$ &$7.5\cdot 10^{-3}$ &$7.5\cdot 10^{-3}$& $\text{mmHg}\cdot\text{s}/\text{mL}$ \\
    $\text{R}_{\text{max}}$ & 75006.2   &75006.2    &75006.2    &75006.2& $\text{mmHg}\cdot\text{s}/\text{mL}$ \\
    \hline
    \end{tabular}
    \caption{Reference parameter values for a healthy individual, used in 
    $(\mathscr{C}_\text{NC})$, describing the heart chambers.}
    \label{table:params-healthy:chambers:nocap}
\end{table}

\begin{table}[H]
    \centering 
    \begin{tabular}{lccc}
    \hline
    Parameter  & $i=\text{SYS}$ & $i=\text{PUL}$ &Unit \\
    \hline
    R$_\text{AR}^i$ & {0.42}               &{0.104} & $\text{mmHg}\cdot\text{s}/\text{mL}$  \\
    L$_\text{AR}^i$ & {$5\cdot 10^{-3}$}   &{$5\cdot 10^{-4}$} & $\text{mmHg}\cdot\text{s}^2/\text{mL}$ \\
    C$_\text{AR}^i$ & {0.96}               &{5}  & $\text{mL}/\text{mmHg}$\\
    R$_\text{VEN}^i$ & {0.352}             &$1.05\cdot 10^{-2}$ & $\text{mmHg}\cdot\text{s}/\text{mL}$  \\
    L$_\text{VEN}^i$ & {$5\cdot 10^{-4}$}  &$5\cdot 10^{-4}$  & $\text{mmHg}\cdot\text{s}^2/\text{mL}$ \\
    C$_\text{VEN}^i$ & {60}                &{16} & $\text{mL}/\text{mmHg}$\\
    \hline
    \end{tabular}
    \caption{Reference parameter values for a healthy individual, used in 
    $(\mathscr{C}_\text{NC})$, describing the circulation.}
    \label{table:params-healthy:RLC:nocap}
\end{table}

\begin{table}[H]
    \centering 
    \begin{tabular}{lccccc}
    \hline
    Parameter  & $i=\text{LA}$  &$i=\text{LV}$      &$i=\text{RA}$  &$i=\text{RV}$ &Unit \\
    \hline
    $\text{E}^a_i$          & 0.38      &2.7        &0.126      & 0.43 & $\text{mmHg}/\text{mL}$ \\
    $\text{E}^p_i$          & 0.27      &$6.9\cdot 10^{-2}$      &0.195      &$4.1264\cdot 10^{-2}$ & $\text{mmHg}/\text{mL}$ \\
    $\text{T}_\text{C}^i$   & 0.1       &0.265      &0.1        &0.3 & s \\
    $\text{T}_\text{R}^i$   & 0.8       &0.4        & 0.7       & 0.4 & s \\
    $\text{t}_\text{C}^i$   & 0.75      &0          &0.8        &0 & s \\
    $\text{V}_{0,i}$        & 4         &3.541      &3.5385     &8.4067 & mL \\
    $\text{R}_{\text{min}}$ & $6.2872\cdot 10^{-3}$    &$6.2872\cdot 10^{-3}$     &$6.2872\cdot 10^{-3}$     &$6.2872\cdot 10^{-3}$ & $\text{mmHg}\cdot\text{s}/\text{mL}$ \\
    $\text{R}_{\text{max}}$ & 94168   &94168    &94168    &94168 & $\text{mmHg}\cdot\text{s}/\text{mL}$ \\
    \hline
    \end{tabular}
    \caption{Reference parameter values for a healthy individual, used in 
    $(\mathscr{C}_\text{C})$, describing the heart chambers.}
    \label{table:params-healthy:chambers:cap}
\end{table}

\begin{table}[H]
    \centering 
    \begin{tabular}{lccc}
    \hline
    Parameter  & $i=\text{SYS}$      & $i=\text{PUL}$ &Unit \\
    \hline
    R$_\text{AR}^i$ &0.5911                 &$7.14\cdot 10^{-2}$ & $\text{mmHg}\cdot\text{s}/\text{mL}$  \\
    L$_\text{AR}^i$ &$2.0643\cdot 10^{-4}$  &$2.0643\cdot 10^{-5}$ & $\text{mmHg}\cdot\text{s}^2/\text{mL}$ \\
    C$_\text{AR}^i$ &1.3315                 &6.0043 & $\text{mL}/\text{mmHg}$ \\
    R$_\text{VEN}^i$ &0.3596                &$3.75\cdot 10^{-2}$ & $\text{mmHg}\cdot\text{s}/\text{mL}$  \\
    L$_\text{VEN}^i$ &$2.0643\cdot 10^{-5}$ &$2.0643\cdot 10^{-5}$ & $\text{mmHg}\cdot\text{s}^2/\text{mL}$ \\
    C$_\text{VEN}^i$ &{75}                &13.181 & $\text{mL}/\text{mmHg}$  \\
    \hline
    \end{tabular}
    \caption{Reference parameter values for a healthy individual, used in 
    $(\mathscr{C}_\text{C})$, describing systemic and pulmonary circulation.}
    \label{table:params-healthy:RLC:cap}
\end{table}

\begin{table}[H]
    \centering 
    \begin{tabular}{lcc}
    \hline
    Parameter  &  &Unit \\
    \hline
    R$_\text{C}^\text{SYS}$ &$2.17\cdot 10^{-2}$ & $\text{mmHg}\cdot\text{s}/\text{mL}$  \\
    C$_\text{C}^\text{SYS}$ &0.27981 & $\text{mL}/\text{mmHg}$ \\
    R$_\text{C}^\text{PUL}$ &$1.7538\cdot 10^{-2}$ & $\text{mmHg}\cdot\text{s}/\text{mL}$  \\
    C$_\text{C}^\text{PUL}$ &5.7803 & $\text{mL}/\text{mmHg}$ \\
    R$_\text{SH}$           &0.35174 & $\text{mmHg}\cdot\text{s}/\text{mL}$  \\
    C$_\text{SH}$           &$4.9043\cdot 10^{-2}$ & $\text{mL}/\text{mmHg}$ \\
    \hline
    \end{tabular}
    \caption{Reference parameter values for a healthy individual, used in 
    $(\mathscr{C}_\text{C})$, describing capillary circulation.}
    \label{table:params-healthy:RLC:cap:capillary}
\end{table}

Note: the parameters presented in 
Tables~\ref{table:params-healthy:chambers:nocap}--\ref{table:params-healthy:RLC:cap:capillary}
are a modification of literature values 
\cite{healthy-params-fonte1,healthy-params-fonte2} in such a way that the
model outputs lie in the healthy individual's ranges, listed in 
\ref{appendix:ranges-healthy}.

\section{Healthy individual's ranges}
\label{appendix:ranges-healthy}%

\begin{table}[H]
    \centering
    \begin{minipage}{0.48\textwidth}
        \centering
        \begin{tabular}{lccc}
            \hline
            Output  & Healthy range  &Unit  & Source \\
            \hline
            $\text{LA}_\text{IVmax}$ & $[24, 57]$ &$\text{mL}/\text{m}^2$ & \cite{healthy-ranges3}   \\
            $\text{LA}_\text{IVmin}$ & $[9, 28]$ &$\text{mL}/\text{m}^2$ & \cite{healthy-ranges3}   \\
            $\text{LA}_\text{IVpreAC}$ & $[15, 46]$ &$\text{mL}/\text{m}^2$ & \cite{healthy-ranges4}   \\
            $\text{LA}_\text{PassEF}$ & $[8, 44]$ &\% & \cite{healthy-ranges4}   \\
            $\text{LA}_\text{ActEF}$ & $[17, 58]$ &\% & \cite{healthy-ranges4}   \\
            $\text{LA}_\text{TotEF}$ & $[37, 70]$ &\% & \cite{healthy-ranges3}   \\
            $\text{LA}_\text{Pmax}$ & $[6, 20]$ &mmHg & \cite{healthy-ranges2}   \\
            $\text{LA}_\text{Pmin}$ & $[-2, 9]$ &mmHg & \cite{healthy-ranges2}   \\
            $\text{LA}_\text{Pmean}$ & $[4, 12]$ &mmHg & \cite{healthy-ranges2}   \\
            \hline
        \end{tabular}
        \caption{Left atrium-related ranges of outputs computed for a healthy individual.}
        \label{table:healthy-ranges:LA}
    \end{minipage}
    \hfill
    \begin{minipage}{0.48\textwidth}
        \centering
        \begin{tabular}{lccc}
            \hline
            Output  & Healthy range  &Unit  & Source \\
            \hline
            $\text{RA}_\text{IVmax}$ & $[28, 76]$ &$\text{mL}/\text{m}^2$ & \cite{healthy-ranges3}   \\
            $\text{RA}_\text{IVmin}$ & $[9, 45]$ &$\text{mL}/\text{m}^2$ & \cite{healthy-ranges3}   \\
            $\text{RA}_\text{IVpreAC}$ & $[19, 61]$ &$\text{mL}/\text{m}^2$ & \cite{healthy-ranges4}   \\
            $\text{RA}_\text{PassEF}$ & $[4, 41]$ &\% & \cite{healthy-ranges4}   \\
            $\text{RA}_\text{ActEF}$ & $[11, 55]$ &\% & \cite{healthy-ranges4}   \\
            $\text{RA}_\text{TotEF}$ & $[29, 68]$ &\% & \cite{healthy-ranges3}   \\
            $\text{RA}_\text{Pmax}$ & $[2, 14]$ &mmHg & \cite{healthy-ranges2}   \\
            $\text{RA}_\text{Pmin}$ & $[-2, 6]$ &mmHg & \cite{healthy-ranges2}   \\
            $\text{RA}_\text{Pmean}$ & $[-1, 8]$ &mmHg & \cite{healthy-ranges2}   \\
            \hline
        \end{tabular}
        \caption{Right atrium-related ranges of outputs computed for a healthy individual.}
        \label{table:healthy-ranges:RA}
    \end{minipage}
\end{table}

\begin{table}[H]
    \centering
    \begin{minipage}{0.48\textwidth}
        \centering
        \begin{tabular}{lccc}
            \hline
            Output  & Healthy range  &Unit  & Source \\
            \hline
            $\text{LV}_\text{ISV}$ & $[30, 66]$ &$\text{mL}/\text{m}^2$ & \cite{healthy-ranges3}   \\
            $\text{LV}_\text{IEDV}$ & $[47, 107]$ &$\text{mL}/\text{m}^2$ & \cite{healthy-ranges3}   \\
            $\text{LV}_\text{IESV}$ & $[11, 47]$ &$\text{mL}/\text{m}^2$ & \cite{healthy-ranges3}   \\
            $\text{LV}_\text{EF}$ & $[51, 76]$ &\% & \cite{healthy-ranges3}   \\
            $\text{LV}_\text{Pmax}$ & $[90, 140]$ &mmHg & \cite{healthy-ranges2}   \\
            $\text{LV}_\text{Pmin}$ & $[4, 12]$ &mmHg & \cite{healthy-ranges2}   \\
            \hline
        \end{tabular}
        \caption{Left ventricle-related ranges of outputs computed for a healthy individual.}
        \label{table:healthy-ranges:LV}
    \end{minipage}
    \hfill
    \begin{minipage}{0.48\textwidth}
        \centering
        \begin{tabular}{lccc}
            \hline
            Output  & Healthy range  &Unit  & Source \\
            \hline
            $\text{RV}_\text{ISV}$ & $[28, 75]$ &$\text{mL}/\text{m}^2$ & \cite{healthy-ranges3}   \\
            $\text{RV}_\text{IEDV}$ & $[53, 123]$ &$\text{mL}/\text{m}^2$ & \cite{healthy-ranges3}   \\
            $\text{RV}_\text{IESV}$ & $[17, 59]$ &$\text{mL}/\text{m}^2$ & \cite{healthy-ranges3}   \\
            $\text{RV}_\text{EF}$ & $[42, 72]$ &\% & \cite{healthy-ranges3}   \\
            $\text{RV}_\text{Pmax}$ & $[15, 28]$ &mmHg & \cite{healthy-ranges2}   \\
            $\text{RV}_\text{Pmin}$ & $[0, 8]$ &mmHg & \cite{healthy-ranges2}   \\
            \hline
        \end{tabular}
        \caption{Right ventricle-related ranges of outputs computed for a healthy individual.}
        \label{table:healthy-ranges:RV}
    \end{minipage}
\end{table}

\begin{table}[H]
    \centering 
    \begin{tabular}{lccc}
    \hline
    Output  & Healthy range  &Unit  & Source \\
    \hline
    $\text{CI}$ & $[2.8, 4.2]$ &$\text{m}^2\cdot \text{L}/\text{min}$ & \cite{healthy-ranges2}   \\
    $\text{SAP}_\text{max}$ & $[0, 140]$ &mmHg & \cite{healthy-ranges1}   \\
    $\text{SAP}_\text{min}$ & $[0, 80]$ &mmHg & \cite{healthy-ranges1}   \\
    $\text{PAP}_\text{max}$ & $[15, 28]$ &mmHg & \cite{healthy-ranges2}   \\
    $\text{PAP}_\text{min}$ & $[5, 16]$ &mmHg & \cite{healthy-ranges2}   \\
    $\text{PAP}_\text{mean}$ & $[10, 22]$ &mmHg & \cite{healthy-ranges2}   \\
    $\text{PWP}_\text{max}$ & $[9, 23]$ &mmHg & \cite{healthy-ranges2}   \\
    $\text{PWP}_\text{min}$ & $[1, 12]$ &mmHg & \cite{healthy-ranges2}   \\
    $\text{PWP}_\text{mean}$ & $[6, 15]$ &mmHg & \cite{healthy-ranges2}   \\
    $\text{SVR}$ & $[11.3, 17.5]$ &$\text{mmHg}\cdot \text{min}/\text{L}$ & \cite{healthy-ranges2}   \\
    $\text{PVR}$ & $[1.9, 3.1]$ &$\text{mmHg}\cdot \text{min}/\text{L}$ & \cite{healthy-ranges2}   \\
    $\text{S}_f$ & $[0,0.05]$ &$1/\text{m}^2$ & \cite{shunt-ranges}   \\
    \hline
    \end{tabular}
    \caption{Further ranges of  
    outputs computed for a healthy 
    individual.}
    \label{table:healthy-ranges:other}
\end{table}

Note: 
the ranges presented in 
Tables~\ref{table:healthy-ranges:LA}--\ref{table:healthy-ranges:other} 
are derived from echocardiography (the ones taken from 
\cite{healthy-ranges1,healthy-ranges2}) and magnetic resonance imaging (MRI) 
(the ones taken 
from \cite{healthy-ranges3,healthy-ranges4}), meaning that using general healthy ranges 
may not always be entirely appropriate; however, 
this is not an issue since, despite the differences between echocardiography 
and MRI, the quantities computed by both $(\mathscr{C}_\text{NC})$ and 
$(\mathscr{C}_\text{C})$ for a healthy individual 
fall within all ranges.

\newpage

\section{Parameter Tuning for Cross-Model Equivalence}
\label{appendix:capillary}%
\begin{algorithm}[t]
    \label{alg:matching}
    \caption{Parameter matching algorithm for capillary and non-capillary 0D models}
    \label{alg:alg1}
    \label{protocol}
    \begin{algorithmic}[1]
    \STATE{\textbf{Input:} $\boldsymbol\theta_0^\text{C}$, $\boldsymbol\theta_0^\text{NC}$, $n$, $k$} 
    \STATE{\textbf{Output:} $\boldsymbol\theta_0$, $\varepsilon$} 
    \STATE{$\boldsymbol\theta_0\longleftarrow\left[\boldsymbol\theta_0^\text{C},\, \boldsymbol\theta_0^\text{NC}\right]$}
    \STATE{$\textsc{MNC}\longleftarrow\textsc{0dModelNoCapillary}(\boldsymbol\theta_0^\text{NC})$}
    \STATE $\textsc{MC}\longleftarrow\textsc{0dModelCapillary}(\boldsymbol\theta_0)$
    \STATE{$\varepsilon\longleftarrow\textsc{ComputeError}(\textsc{MNC}, \textsc{MC})$}
    \FOR{$\theta \in \boldsymbol\theta_0$}
    \FOR{$j=1,\dots,n$}
    \STATE{$\bar{\theta}\longleftarrow\theta\cdot k$}
    \STATE $\textsc{MC}\longleftarrow\textsc{0dModelCapillary}\bigl(\textsc{UpdateParameter}(\bar{\theta}, \boldsymbol\theta_0)\bigr)$
    \STATE{$\delta\longleftarrow\textsc{ComputeError}(\textsc{MNC}, \textsc{MC})$}
    \IF{$\delta<\varepsilon$}
    \STATE{$\boldsymbol\theta_0\longleftarrow\textsc{UpdateParameter}(\bar{\theta}, \boldsymbol\theta_0)$}
    \STATE{$\varepsilon\longleftarrow\delta$}
    \ELSE
    \STATE{\textbf{break}}
    \ENDIF
    \ENDFOR
    \ENDFOR
    \RETURN{$\boldsymbol\theta_0$, $\varepsilon$}
    \end{algorithmic}
\end{algorithm}

\begin{figure}[t]
    \centering
    \includegraphics[width=140 mm]{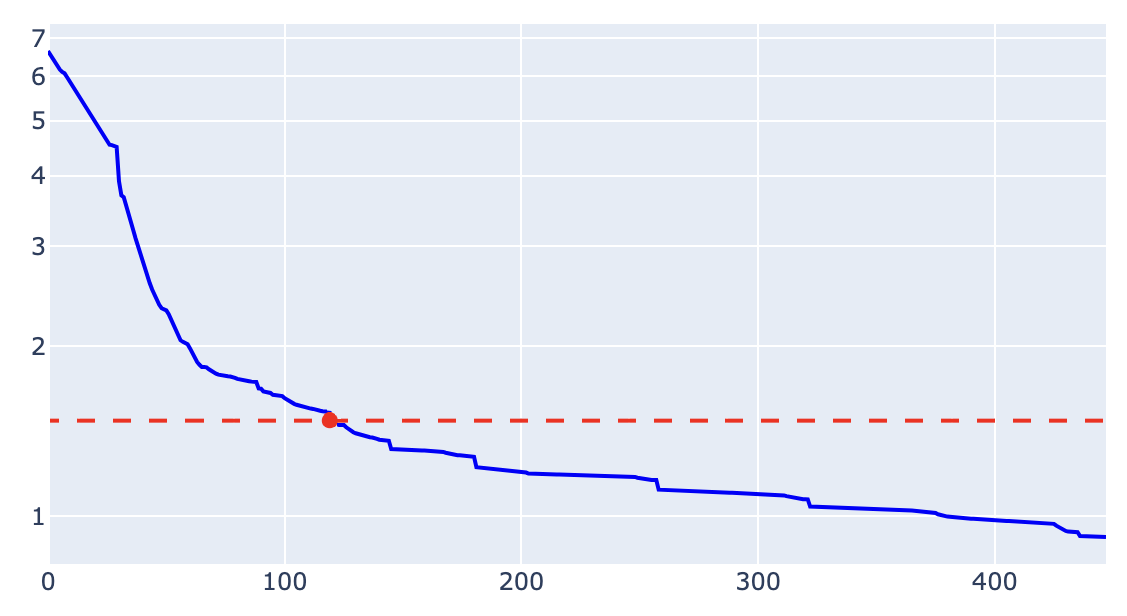}
    \caption{Evolution of the error (in blue) measuring the distance 
    between the output variables of $(\mathscr{C}_{\text{NC}})$ 
    and $(\mathscr{C}_{\text{C}})$, as a function of the number of 
    iterations, until it reaches 0.92. The y-axis is in logarithmic scale. The red 
    point indicates the iteration where the error falls below the threshold of 1.5, 
    which is considered an acceptable value.}
    \label{fig:error_capillary}
\end{figure}

When comparing two sets of output coming from simulations carried on with different 0D models, 
it is essential to assess whether the presence of capillaries 
introduces significant variations that compromise the validity of the comparison.
The idea is to determine the parameters of $(\mathscr{C}_{\text{C}})$ such that its 
simulated outputs are as close as possible to those obtained from 
$(\mathscr{C}_{\text{NC}})$. 
For this purpose, Algorithm~\ref{alg:alg1} has been implemented. Firstly, a set of 
parameters $\boldsymbol\theta_0^\text{NC}$ is chosen for 
$(\mathscr{C}_{\text{NC}})$, and a simulation (MNC) is performed. Then, a 
simulation (MC) is executed by means of $(\mathscr{C}_{\text{C}})$: in this case, 
the same parameter set $\boldsymbol\theta_0^\text{NC}$ is employed, with the addition 
of values corresponding to capillary circulation parameters, denoted by 
$\boldsymbol\theta_0^\text{C}$.
The initial tolerance $\varepsilon$ is calculated as 
the error between the outputs of these two simulations in the following way:
\begin{equation}\label{eq:alg:err}
    \varepsilon=\sum_{\alpha\in A} \dfrac{\delta V_{\alpha}}{\log{\left(2+\overline{V}_{\alpha}^{\text{C}}\right)}},
\end{equation}
where:
\begin{align*}
    \delta V_{\alpha}&=\int_{\text{T}-\text{T}_{\text{HB}}}^{\text{T}} \left|V_{\alpha}^{\text{NC}}(t)-
    V_{\alpha}^{\text{C}}(t)\right|\,dt, & \overline{V}_{\alpha}^{\text{C}}&=
    \dfrac{1}{\text{T}_{\text{HB}}}\int_{\text{T}-\text{T}_{\text{HB}}}^{\text{T}} \left|V_{\alpha}^{\text{C}}(t)\right|\,dt,
\end{align*}
with $\left\{V_\alpha\right\}_{\alpha\in A}$ being the set of variable that contribute to the error computation; 
it comprehends a total of 20 variables: pressures and volumes of the four heart chambers, pressures and fluxes of systemic 
and pulmonary (arterial and venous) circulation, and fluxes through the four heart valves.  
Having defined the initial guess $\boldsymbol\theta_0=[\boldsymbol\theta_0^\text{C},\boldsymbol\theta_0^\text{NC}]$, the parameters of 
$(\mathscr{C}_{\text{C}})$ are adjusted iteratively: each parameter $\theta$ is modified individually using a 
scaling coefficient $k$. After each modification, a new simulation of MC is performed, and the 
error $\delta$ is computed as the distance between the outputs of MNC and the 
current MC, using Equation~\eqref{eq:alg:err}. If the new error is smaller than that of the previous step,
the current parameter $\theta$ is updated (line 13 
in Algorithm~\ref{alg:alg1}; the function $\textsc{UpdateParameter}(\bar{\theta}, \boldsymbol\theta_0)$ 
takes as input $\bar{\theta}$ and $\boldsymbol\theta_0$, and returns as output the new set 
of parameters, in which $\theta$ has been replaced with $\bar{\theta}$), 
and the previous error $\varepsilon$ is replaced with the new one. This procedure is 
repeated $n$ times for each parameter.

In order to reduce the error between the outputs of 
the two simulations below a desired threshold, the whole algorithm has to 
be iterated several times, for different values of $n$ and $k$. Specifically, 
$k$ is chosen within a neighbourhood of 1 ($k \in U(1)$), 
alternating between values greater and less than 1: this ensures that parameters can 
be either increased or decreased as needed. It 
is also possible to apply the algorithm backwards, adapting the parameters of the 
$(\mathscr{C}_{\text{NC}})$ so that its output variables are as close as possible to 
those of the $(\mathscr{C}_{\text{C}})$. In this way, depending on the application to 
be made or the model that is used as reference, the procedure can be flexibly applied 
in either direction.
A strict threshold of 0.92 is chosen, since reaching a lower threshold becomes 
computationally expensive without significantly reducing the error;
furthermore, a second 
threshold of 1.5 is adopted, representing a good balance 
between computational cost and quality of the parameters. 
The whole analysis has been performed on a MacBook Pro
(Apple M1 Pro, 10 cores, 3.2 GHz, 16 GB RAM); typically, the entire process takes only 
a few hundreds of iterations, in which the parameter set $\boldsymbol{\theta}_0$ 
is effectively updated (lines 12--14 in Algorithm~\ref{alg:alg1}): the reason is 
partly because the same starting parameter set, 
$\boldsymbol{\theta}_0^\text{NC}$, is employed for both models $(\mathscr{C}_{\text{NC}})$ and 
$(\mathscr{C}_{\text{C}})$, hence the variables are already sufficiently 
close at the first step.

In order to highlight the versatility of Algorithm~\ref{alg:alg1}, an example is provided
in which the goal is optimizing parameters of $(\mathscr{C}_{\text{C}})$ in such a 
way that it minimizes deviations of its output variables from the reference model 
$(\mathscr{C}_{\text{NC}})$. As shown in Figure~\ref{fig:error_capillary}, to achieve 
an error below the specified 
threshold of 0.92, approximately 76 minutes are 
required; nevertheless, an
acceptable solution, defined as one having an error less than 1.5, is already 
obtained in 16 minutes. Importantly, all time-independent outputs stay within healthy 
ranges (\ref{appendix:ranges-healthy} for more information about the healthy ranges), 
validating the adjusted parameters. To reduce 
the computation time even further, it is possible to remove some variables from the 
matching process, or the tolerance may be increased, depending on the desired accuracy. 
The values of the parameters 
used in this example as initial guess
are provided in \ref{appendix:param-healthy}, while the changes applied to the 
parameters of $(\mathscr{C}_{\text{C}})$ are shown in 
Table~\ref{table:hyper-changes}.

In conclusion, the outputs from $(\mathscr{C}_{\text{C}})$ and $(\mathscr{C}_{\text{NC}})$ 
can be compared reliably since, by adjusting the parameters of one model to the other, 
very similar results can be achieved. 

\newpage

\bibliographystyle{unsrt}
\bibliography{arxiv_celora_et_al_2025.bib}

\end{document}